\documentclass[a4paper,twoside,12pt,DIVcalc]{scrbook}

\usepackage{ngerman,german}
\usepackage[dvips]{hyperref}
\usepackage[latin 1]{inputenc}
\usepackage{graphicx,supertabular}
\usepackage{mathrsfs,amsmath,amsfonts,amssymb,amscd,stmaryrd}
\usepackage{oldgerm}
\usepackage{ifthen}
\usepackage{pifont}
\usepackage{verbatim}
\usepackage{makeidx, showidx, euscript, array}
\usepackage{ntheorem}
\usepackage{paralist}
\usepackage[T1]{fontenc}
\usepackage{scrpage}
\usepackage[DIVcalc]{typearea}

  \usepackage[all,poly,dvips]{xy}
  \usepackage[small]{caption}
  \SelectTips{cm}{}

\theoremstyle{changebreak}

\newtheorem{thm}{Theorem}[section]
\newtheorem{kor}[thm]{Korollar.}

\newtheorem{lem}[thm]{Lemma.}

\newtheorem{prop}[thm]{Proposition.}
\newtheorem{defn}[thm]{Definition.}

\newtheorem{defn+folg}[thm]{Definition und Folgerung.}

\newtheorem{defn+sat}[thm]{Definition und Satz.}
\newtheorem{defn+lem}[thm]{Definition und Lemma.}

\newtheorem{bem}[thm]{Bemerkung.}

\newtheorem{sat}[thm]{Satz.}

\newtheorem{bsp}[thm]{Beispiel.}

\newcommand{\adj}{\operatorname{adj}}

\newcommand{\Char}{char}

\newcommand{\rk}{rk}

\newcommand{\id}{\operatorname{id}}

\newcommand{\Hom}{\operatorname{Hom}}

\newcommand{\trdeg}{\operatorname{trdeg}}

\newcommand{\Gl}{\operatorname{Gl}}

\renewcommand{\Im}{\operatorname{Im}}

\newcommand{\bewub}{\noindent \textbf{Beweis  } }

\newcommand{\Iso}{\operatorname{Iso}}

\newcommand{\Mod}{\texttt{Mod}}

\newcommand{\nat}{\mathbb{N}}

\makeatletter \@addtoreset{equation}{section}

\makeatother

\newcounter{drawingnr}[section]
\newcounter{subdrawingnr}[drawingnr]

\newcounter{Beispnr}[section]

\newcounter{Aufgnr}[section]

\newboolean{showborder}

\newcommand\diagramiso{\hbox{\lower1ex\hbox{$\scriptstyle\sim$}}}
\newcommand\Rep{Rep}

\DeclareMathOperator{\preim}{pre\kern0.13em im} %
\DeclareMathOperator{\preker}{pre\kern0.13em ker} %
\DeclareMathOperator{\precoker}{pre\kern0.13em coker} %

\DeclareMathOperator{\Cat}{Cat} %
\DeclareMathOperator{\mor}{Mor} %
\DeclareMathOperator{\rg}{rg} %
\DeclareMathOperator{\spec}{spec} %
\DeclareMathOperator{\Lie}{Lie}%

\DeclareMathOperator{\supp}{supp} %

\newcommand\set[2][auto]{
     \ifthenelse{\equal{#1}{auto}}{\left\lbrace}{\csname #1\endcsname\lbrace} #2 \ifthenelse{\equal{#1}{auto}}{\right\rbrace}{\csname #1\endcsname\rbrace} }

\newcommand\ob{\operatorname{Obj}}

\newcommand{\indlim}{\operatorname*{\underrightarrow{\lim}}} %
\renewcommand{\projlim}{\operatorname*{\underleftarrow{\lim}}}

\newcommand\longto{\to}

\newcommand\into{\hookrightarrow}

\newcommand{\Aut}{\ensuremath{\mathrm{Aut}}}
\newcommand{\Gal}{\ensuremath{\mathrm{Gal}}}

\newcommand{\Spec}{\ensuremath{\mathrm{Spec}}}

\mathchardef\ordinarycolon\mathcode`\: \mathcode`\:=\string"8000
\begingroup \catcode`\:=\active
  \gdef:{\mathrel{\mathop\ordinarycolon}}

\endgroup

\begin{document}
\hypersetup{pdfauthor={Urs Hackstein}}
\hypersetup{pdftitle={Prinzipalbündel auf p-adischen Kurven und Paralleltransport}}
\hypersetup{pdfsubject={}}

 \thispagestyle{empty}
\begin{center}
\text{}\\
\vspace{7cm}
{\Large \textbf{Urs Hackstein\\\vspace{2cm}
Prinzipalbündel auf $p$"=adischen Kurven und Paralleltransport\\\vspace{1cm}
--2006--}}
\end{center}
\newpage
\thispagestyle{empty}
\text{}
\newpage
\thispagestyle{empty}
 \begin{center}
 {\large{Reine Mathematik}\\
 \vspace{5cm}
\Huge{\textbf{Prinzipalbündel auf p-adischen Kurven und Paralleltransport}}\\
 \vspace{4cm}
 \large Inaugural"=Dissertation zur Erlangung des Doktorgrades der Naturwissenschaften im Fachbereich Mathematik und Informatik der Mathematisch"=Naturwissenschaftlichen Fakultät der Westfälischen Wilhelms"=Universität Münster\\
\vfill
vorgelegt von \\ \vspace{1cm}
  Urs Hackstein\\
  aus Köln\\
  --2006--}
 \end{center}\newpage
\thispagestyle{empty}
\text{}
\vfill
\begin{tabular}{ll}
Dekan: & Prof. Dr. Joachim Cuntz\\
Erster Gutachter: & Prof. Dr. Christopher Deninger\\
Zweite Gutachterin: & Prof. Dr. Annette Werner\\
& (Universität Stuttgart)\\
Tag der mündlichen Prüfung: &15. Dezember 2006\\
Tag der Promotion: &31. Januar 2007\\
\end{tabular}

 \tableofcontents\markboth{}{}\newpage\thispagestyle{empty}

\chapter*{Einleitung}
\markboth{EINLEITUNG}{}\thispagestyle{plain}
\addcontentsline{toc}{chapter}{\numberline{}Einleitung}

Auf einer kompakten Riemannschen Fläche induziert jede endlich"=dimensionale komplexe Darstellung der Fundamentalgruppe ein flaches Vektorbündel und somit ein holomorphes Vektorbündel. Nach einem klassischen Resultat von André Weil wird genauer ein holomorphes Vektorbündel auf einer Riemannschen Fläche genau dann durch eine Darstellung der Fundamentalgruppe gegeben, wenn jede seiner unzerlegbaren Komponenten vom Grad Null ist (siehe \cite{W}).\par
Später konnten Mudumbai S. Narasimham und Conjeeveram S. Seshadri in einem bekannten Theorem zeigen, dass es auf einer Riemannschen Fläche vom Geschlecht $g\geq 2$ eine Äquivalenz von Kategorien zwischen der Kategorie der unitären Darstellungen der Fundamentalgruppe und der Kategorie der polystabilen Vektorbündel vom Grad Null gibt (siehe \cite[Corollary 12.1]{NS}). Außerdem kommt jedes stabile Vektorbündel vom Grad Null von einer irreduziblen unitären Darstellung (\cite[Corollary 12.2]{NS}).\par
Diese Resultate wurden von Annamalai Ramanathan in \cite{R} auf den Fall von $G$"=Prinzipalbündeln auf einer Riemannschen Fläche $X$ ausgedehnt, wobei $G$ eine zusammenhängende reduktive algebraische Gruppe über $\mathbb{C}$ ist. Er zeigt in \cite{R}, wie man ausgehend von gewissen Darstellungen von $\pi_{1}(X-x_{0})$, wobei $x_{0}$ ein Punkt von $X$ ist, $G$"=Prinzipalbündel auf $X$ erhält. Zur genaueren Beschreibung dehnt Ramanathan dabei die bekannte Definition der (Semi)stabilität von Vektorbündeln nach David Mumford (siehe \cite{Mf1}) auf den Fall von $G$"=Prinzipalbündeln aus, indem er definiert, dass ein holomorphes $G$"=Prinzipalbündel $E$ auf $X$ (semi)stabil heiße, wenn für jede maximale parabolische Untergruppe $P$ von $G$ und jeden Schnitt $\sigma\colon X\rightarrow E/P$ gilt, dass $\deg(\sigma^{*}T_{E/P})>0$ (bzw. $\geq 0$) ist, wobei $T_{E/P}$ das relative Tangentialbündel entlang der Fasern von $E/P\to X$ ist.\par
Mittels dieses Stabilitätsbegriffs beweist er, dass ein $G$"=Prinzipalbündel genau dann stabil ist, wenn es von einer irreduziblen unitären Darstellung von $\pi_{1}(X-x_{0})$ herkommt. Eine Darstellung $\rho$ heißt dabei unitär, wenn ihr Bild in einer fixierten maximalen kompakten Untergruppe von $G$ enthalten ist, und sie heißt zusätzlich irreduzibel, wenn die Teilmenge \mbox{$\{Y; \text{ad}h(Y)=Y\,\forall\,h\in\Im\rho\}$} der Lie"=Algebra von $G$ mit deren Zentrum übereinstimmt.\par
In ihrer Arbeit \cite{DW1} entwickelten Christopher Deninger und Annette Werner ein teilweises $p$"=adisches Analogon zu den Ergebnissen von Narasimham und Seshadri. Dazu betrachten sie Vektorbündel $E$ auf $X_{\mathbb{C}_{p}}$, wobei $X$ eine glatte und projektive Kurve über $\overline{\mathbb{Q}_{p}}$ ist, die potentiell streng semistabile Reduktion vom Grad Null besitzen, d.\,h. solche, für die es einen endlichen und étalen Morphismus $\alpha\colon X'\to X$ glatter und projektiver Kurven über $\overline{\mathbb{Q}_{p}}$ gibt, so dass $\alpha_{\mathbb{C}_{p}}^{*}E$ zu einem Vektorbündel $\mathcal{E}$ auf $\mathfrak{X'}_{\mathfrak{o}}=\mathfrak{X'}\otimes\mathfrak{o}$ für ein geeignetes Modell $\mathfrak{X'}$ von $X'$ über $\overline{\mathbb{Z}_{p}}$ ausgedehnt werden kann ($\mathfrak{o}$ sei dabei der Ring der ganzen Zahlen von $\mathbb{C}_{p}$), so dass das Pullback der speziellen Faser $\mathcal{E}_{k}$ von $\mathcal{E}$ zu der Normalisierung jeder irreduziblen Komponente von $\mathfrak{X'}_{k}$ streng semistabil vom Grad Null ist. Zu jedem solchen Vektorbündel $E$ auf $X_{\mathbb{C}_{p}}$ konstruieren sie funktorielle Isomorphismen von "`Paralleltransport"' entlang étaler Wege zwischen den Fasern von $E$ auf $X_{\mathbb{C}_{p}}$, deren Konstruktion mit Tensorprodukten, der Bildung dualer Objekte, inneren Homomorphismen, Pullbacks und äußeren Produkten von Vektorbündeln sowie mit Galoiskonjugation verträglich ist. Insbesondere erhält man eine Darstellung $\rho_{E,x}$ von $\pi_{1}(X,x)$ auf $x^{*}E=:E_{x}$ für jeden Punkt $x$ in $X(\mathbb{C}_{p})$.
In einer weiteren nachfolgenden Arbeit (\cite{DW2}) konnten Deninger und Werner zeigen, dass die Ausgangskategorie der Vektorbündel $E$ auf $X_{\mathbb{C}_{p}}$ mit potentiell streng semistabiler Reduktion vom Grad Null eine neutrale Tannaka"=Kategorie ist. Ebenso wies Gabriel Herz nach, dass die Konstruktion von Deninger und Werner auf Mumfordkurven kompatibel mit einer älteren Konstruktion von Gerd Faltings ist (siehe \cite{He} und \cite{F1}). Außerdem vereinfachte Jilong Tong mit Hilfe der Theorie der TAC"=Kurven in \cite{T} den Beweis der Resultate von Deninger und Werner. Für einen Überblick über diesen Themenkomplex sei auch auf den Überblicksartikel \cite{We} hingewiesen.\par

Es sei ferner in diesem Zusammenhang angemerkt, dass sich die Resultate von Narasimham und Seshadri bzw. Deninger und Werner in eine allgemeinere Theorie einbetten lassen: In \cite{S} bewies Carlos Simpson eine Äquivalenz von Kategorien zwischen sogenannten stabilen Higgs"=Bündeln auf einer Kählermannigfaltigkeit $X$ über $\mathbb{C}$ und irreduziblen lokalen Systemen auf $X$. Ein Higgs"=Bündel ist dabei ein Paar $(E,\theta)$ bestehend aus einem holomorphen Vektorbündel $E$ auf $X$ zusammen mit einer holomorphen Abbildung $\theta\colon E\to  E\otimes\Omega^{1}_{X}$, die $\theta\wedge\theta =0$ erfüllt, während lokale Systeme auf $X$ zu Darstellungen der Fundamentalgruppe von $X$ korrespondieren. Für den Spezialfall $\theta =0$ erhält man aus den Resultaten von Simpson die Theorie von Narasimham und Seshadri.\par
In \cite{F2} konstruierte wiederum Gerd Faltings ein $p$"=adisches Analogon zu dieser Simpson"=Korrespondenz, indem er mit Hilfe seiner Theorie der "`almost étale extensions"' eine Äquivalenz von Kategorien zwischen der Kategorie der Higgs"=Bündel auf einer $p$"=adischen Kurve $X$ und der Kategorie gewisser "`verallgemeinerter Darstellungen"' der étalen Fundamentalgruppe von $X$ zeigen konnte.\\\par 

In dieser Arbeit wird nun die Konstruktion von Deninger und Werner im Fall von Vektorbündeln auf $X_{\mathbb{C}_{p}}$ auf den Fall von $G$"=Torseuren auf $X_{\mathbb{C}_{p}}$ ausgeweitet, wobei $G$ ein zusammenhängendes reduktives Gruppenschema von endlicher Präsentation über $\mathfrak{o}$ ist. Es werden dabei $G$"=Torseure $E$ betrachtet, die potentiell streng semistabile Reduktion vom Grad Null haben. Wir verwenden dabei die oben genannte Definition der Semistabilität von Torseuren nach Ramanathan, die auch über einem beliebigen Körper beliebiger Charakteristik gültig bleibt, und sagen in Analogie zum Fall der Vektorbündel, dass ein $G$"=Torseur $E$ auf $X_{\mathbb{C}_{p}}$ streng semistabile Reduktion vom Grad Null habe, wenn $E$ sich zu einem $G$"=Torseur $\widetilde{E}$ auf $\mathfrak{X}_{\mathfrak{o}}=\mathfrak{X}\otimes\mathfrak{o}$ für ein geeignetes Modell $\mathfrak{X}$ von $X$ ausdehnt, so dass das Pullback der speziellen Faser $\widetilde{E}_{k}$ von $\widetilde{E}$ zu der Normalisierung jeder irreduziblen Komponente von $\mathfrak{X}_{k}$ streng semistabil vom Grad Null ist. Gibt es einen endlichen und étalen Morphismus $\alpha\colon Y\to X$ glatter und projektiver Kurven über $\overline{\mathbb{Q}_{p}}$, so dass $\alpha_{\mathbb{C}_{p}}^{*}E$ streng semistabile Reduktion vom Grad Null hat, so besitzt $E$ potentiell streng semistabile Reduktion vom Grad Null.  Ist $E$ ein derartiger $G$"=Torseur auf $X_{\mathbb{C}_{p}}$, so gibt es funktorielle Isomorphismen von "`Paralleltransport"' entlang étaler Wege der Fasern von $E_{\mathbb{C}_{p}}$ auf $X_{\mathbb{C}_{p}}$. Insbesondere gibt es für jeden solchen $G$"=Torseur $E$ einen stetigen Funktor $\rho_{E}$ von $\Pi_{1}(X)$ nach $\mathcal{P}(G(\mathbb{C}_{p}))$ mit $\rho_{E}(x)=E_{x}=x^{*}E$ für alle $x\in X(\mathbb{C}_{p})$. Dabei ist $\mathcal{P}(G(\mathbb{C}_{p}))$ die Kategorie der topologischen Räume mit einer einfach"=transitiven und stetigen $G(\mathbb{C}_{p})$"=Aktion von rechts, wobei die Morphismen gerade die stetigen $G(\mathbb{C}_{p})$"=äquivarianten Abbildungen sind.
Insgesamt erhält man so einen Funktor $\rho\colon\mathcal{B}_{X_{\mathbb{C}_{p}}}^{ps}\to\Rep_{\Pi_{1}(X)}(G(\mathbb{C}_{p}))$, wobei $\mathcal{B}_{X_{\mathbb{C}_{p}}}^{ps}$ die Kategorie der $G$"=Torseure mit potentiell streng semistabiler Reduktion vom Grad Null und $\Rep_{\Pi_{1}(X)}(G(\mathbb{C}_{p}))$ die Kategorie der stetigen Funktoren von $\Pi_{1}(X)$ nach $\mathcal{P}(G(\mathbb{C}_{p}))$ ist, der sich funktoriell bezüglich Morphismen von Kurven über $\overline{\mathbb{Q}_{p}}$, Morphismen von zusammenhängenden reduktiven Gruppenschemata von endlicher Präsentation über $\mathfrak{o}$ und $\mathbb{Q}_{p}$"=Automorphismen von $\overline{\mathbb{Q}_{p}}$ verhält und verträglich mit den entsprechenden Funktoren im Vektorbündelfall ist.\par
Der Aufbau der Arbeit ist wie folgt:\par
Im ersten Kapitel werden einige wichtige Resultate und Definitionen zu Torseuren, zum étalen Fundamentalgruppoid und zu gewissen Kategorien von Überdeckungen, die in \cite{DW1} definiert wurden, wiederholt und zusammengestellt.\par
Aufbauend auf ersten Skizzen von Deninger (\cite{DW3}) wird dann im zweiten Kapitel ein étaler Paralleltransport zunächst für $G$"=Torseure definiert, die für einen gegebenen Divisor $D$ auf $X$ in der Kategorie $\mathcal{B}_{\mathfrak{X}_{\mathfrak{o}},D}(G)$ bzw. der Kategorie $\mathcal{B}_{X_{\mathbb{C}_{p}},D}(G)$ liegen. Für diese Konstruktionen genügt es, an $G$ lediglich zu fordern, dass es sich um ein affines, glattes Gruppenschema von endlicher Präsentation über $\mathfrak{o}$ handelt, also nicht notwendigerweise auch um ein zusammenhängendes reduktives Gruppenschema. In der Definition von Deninger besteht dabei die Kategorie $\mathcal{B}_{\mathfrak{X}_{\mathfrak{o}},D}(G)$ für ein fest gewähltes Modell von $X$ aus allen $G$"=Torseuren $E$ auf $\mathfrak{X}_{\mathfrak{o}}$, so dass für jedes $n\geq 1$ eine Überdeckung $\pi\colon\mathcal{Y}\to\mathfrak{X}$ von $\mathfrak{X}$ existiert, die zu der in \cite{DW1} definierten Kategorie von Überdeckungen $S_{\mathfrak{X},D}$ gehört, so dass $\pi^{*}E$ trivial modulo $p^{n}$ ist. Die Kategorie $\mathcal{B}_{X_{\mathbb{C}_{p}},D}(G)$ wiederherum besteht aus den $G$"=Torseuren, die isomorph zu der generischen Faser eines $G$"=Torseurs $\widetilde{E}$ aus $\mathcal{B}_{\mathfrak{X}_{\mathfrak{o}},D}(G)$ für ein Modell $\mathfrak{X}$ von $X$ sind. \par
Die Konstruktion des Paralleltransport verläuft dann, den Skizzen von Deninger folgend, weitgehend analog zu der Konstruktion des Paralleltransports im Fall der Vektorbündel: Für einen gegebenen $G$"=Torseur $E\in\mathcal{B}_{\mathfrak{X}_{\mathfrak{o}},D}(G)$ und ein festes $n\geq 1$ sei $\pi\colon\mathcal{Y}\to\mathfrak{X}$ eine Überdeckung in $S_{\mathfrak{X},D}$, so dass $\pi_{n}^{*}E_{n}$ ein trivialer $G$"=Torseur auf $\mathcal{Y}_{n}$ ist. Hierbei bezeichne der Index $n$ stets die Reduktion modulo $p^{n}$. Man betrachte nun zwei Punkte $x$ und $x'$ in $X(\mathbb{C}_{p})=\mathfrak{X}(\mathfrak{o})$ und man wähle einen Punkt $y$ in $Y=\mathcal{Y}_{\mathbb{C}_{p}}$ über $x$. Für einen étalen Weg $\gamma$ von $x$ nach $x'$, d.\,h. einen Isomorphismus von Faserfunktoren zwischen den Faserfunktoren $F_{x}$ und $F_{x'}$, sei $\gamma y$ der korrespondierende Punkt über $x'$. Liegt $\pi$ in der in \cite{DW1} definierten vollen Unterkategorie $S^{good}_{\mathfrak{X},D}$ von $S_{\mathfrak{X},D}$, so hat man Isomorphismen \[E_{x_{n}}:=(x_{n}^{*}E)(\mathfrak{o}_{n})\stackrel{y_{n}^{*}}{\stackrel{\sim}{\longleftarrow}}\Gamma(\mathcal{Y}_{n},\pi_{n}^{*}E_{n})\stackrel{(\gamma y)_{n}^{*}}{\stackrel{\sim}{\longrightarrow}}((x'_{n})^{*}E)(\mathfrak{o}_{n})=:E_{x'_{n}}.\] Die so gegebenen Abbildungen $\rho_{E,n}(\gamma)=(\gamma y)_{n}^{*}\circ(y_{n}^{*})^{-1}$ bilden ein projektives System und wir definieren den Paralleltransport $\rho_{E}(\gamma)\colon E_{x}\xrightarrow{\sim}E_{x'}$ als ihren projektiven Limes. 
Setzt man $\rho_{E}(x):=E_{x}$ für alle $x\in (X-D)(\mathbb{C}_{p})$, so ergibt dies einen stetigen Funktor von $\mathcal{B}_{\mathfrak{X}_{\mathfrak{o}},D}(G)$ in die Kategorie der stetigen Funktoren vom étalen Fundamentalgruppoid $\Pi_{1}(X-D)$ in die oben definierte Kategorie $\mathcal{P}(G(\mathfrak{o}))$.
Dieser Funktor wird anschließend zu einem stetigen Funktor von $\mathcal{B}_{X_{\mathbb{C}_{p}},D}(G)$ in die Kategorie der stetigen Funktoren vom étalen Fundamentalgruppoid $\Pi_{1}(X-D)$ in die analog zu $\mathcal{P}(G(\mathfrak{o}))$ definierte Kategorie $\mathcal{P}(G(\mathbb{C}_{p}))$ ausgedehnt.\par
Zum Abschluss des Kapitels wird gezeigt, dass diese Konstruktion sich funktoriell unter Morphismen glatter und projektiver Kurven über $\overline{\mathbb{Q}_{p}}$, Morphismen afffiner, glatter Gruppenschemata von endlicher Präsentation über $\mathfrak{o}$ und unter Galoiskonjugation verhält und kompatibel mit den im Falle der Vektorbündel konnstruierten Funktoren ist.\par
Um den Beweis des eingangs erwähnten Theorems zu ermöglichen, wird dann im dritten Kapitel der Arbeit eine handlichere Charakterisierung der Kategorie $\mathcal{B}_{\mathfrak{X}_{\mathfrak{o}},D}(G)$ erarbeitet, indem analoge Resultate zu den Theoremen 17,18 und 20 von \cite{DW1} bewiesen werden. Als erstes wird gezeigt, dass es genügt zu wissen, dass $\pi_{k}^{*}E_{k}$ für ein geeignetes $\pi\in S_{\mathfrak{X},D}$ trivial ist, um $E\in\mathcal{B}_{\mathfrak{X}_{\mathfrak{o}},D}(G)$ nachzuweisen. Dies gilt unter der Voraussetzung, dass $G$ ein affines, glattes Gruppenschema über $\mathfrak{o}$ ist. Setzt man anschließend stets voraus, dass $G$ sogar ein zusammenhängendes reduktives Gruppenschema von endlicher Präsentation über $\mathfrak{o}$ ist, so erhält man in Analogie zu \cite[Theorem 17]{DW1}, dass ein $G$"=Torseur auf $\mathfrak{X}_{o}$ genau dann in $\mathcal{B}_{\mathfrak{X}_{\mathfrak{o}},D}(G)$ für einen gewissen Divisor $D$ liegt, wenn er streng semistabile Reduktion vom Grad Null besitzt. Entscheidend dafür ist ähnlich wie im Fall der Vektorbündel die Charakterisierung der $G$"=Torseure auf einem rein"=eindimensionalen eigentlichen Schema über einem endlichen Körper $\mathbb{F}_{q}$, die streng semistabile Reduktion vom Grad Null haben, d.\,h. deren Pullback zu der Normalisierung jeder irreduziblen Komponente streng semistabil vom Grad Null ist. Es wird gezeigt, dass dies gerade die Torseure sind, deren Pullback entlang eines endlichen surjektiven Morphismus zu einem rein"=eindimensionalen eigentlichen $\mathbb{F}_{q}$"=Schema trivial ist. Die Rolle des im Vektorbündelfall zugrundeliegenden Resultats von Herbert Lange und Ulrich Stuhler, das die Charakterisierung für den Fall einer glatten und projektiven Kurve lieferte, übernimmt hier ein Resultat von Deligne (siehe Satz \ref{Laszlo}) in \cite{Las}, aus dem sich ein Zusammenhang zwischen Trivialisierbarkeit des gegebenen Torseurs und der Tatsache, dass er isomorph zu seinem Pullback unter Frobenius ist, ergibt. Daraus folgt wiederum wie in \cite{LS} die gegebene Beschreibung für den Fall einer glatten und projektiven Kurve.\par
Das letzte Kapitel besteht schließlich aus dem Beweis der beiden Hauptresultate der Arbeit unter Vewendung der im vorherigen Kapitel erzielten Charakerisierung der Kategorie $\mathcal{B}_{\mathfrak{X}_{\mathfrak{o}},D}(G)$.\\\par

Zuallererst möchte ich Prof. Christopher Deninger für die gute Betreuung meiner Arbeit und die Ideen, wie man im Falle von $G$"=Torseuren einen Paralleltransport definieren könne, sowie für die Erlaubnis, seine ersten Resultate zu einer solchen Definition des Paralleltransports und zu einigen Funktorialitätseigenschaften der Konstruktion (vgl. \cite{DW3}), die den Ausgangspunkt dieser Arbeit darstellen, verwenden zu dürfen, herzlichst danken. Ebenso danke ich ihm für viele weiterführende Diskussionen und sein Vertrauen. Außerdem gebührt mein Dank dem Marie Curie Research Training Network "`Arithmetic Algebraic Geometry"' (MRTN"=CT"=2003"=504917) der Europäischen Union, das mir einen sechsmonatigen Aufenthalt am Laboratoire des Mathématiques der Université de Paris-Sud in Orsay ermöglichte, ohne den diese Arbeit nicht hätte entstehen können und während dessen ich vom dortigen Koordinator des Netzwerkes, Prof. Étienne Fouvry, Unterstützung bei allen denkbaren administrativen Fragen erhalten habe. Für viele fruchtbare Diskussionen möchte ich ferner Pham Ngoc Duy, Philippe Gille, Sylvain Maugeais, Michel Raynaud, Jean-Marc Fontaine, Alexander Steinmetz, Gabriel Herz, Stefan Wiech, Thomas Ludsteck, Jilong Tong, Vikram Metha, Yogish Holla und Markus Reineke herzlichst danken. Yves Laszlo danke ich für den Hinweis auf das oben erwähnte Resultat von Deligne in seinem Artikel \cite{Las}, das sich als sehr wertvoll für die vorliegende Arbeit herausgestellt hat. Nicht zu vergessen sind in diesem Themenbereich auch mehrere Diskussionen mit Philippe Gille und Pham Ngoc Duy, die sehr hilfreich für ein besseres Verständnis waren.\par
Schließlich wäre die vorliegende Arbeit auch nicht möglich gewesen ohne alle Menschen und Institutionen wie dem Collegium musicum instrumentale der Universität Münster, dem Ensemble des violoncelles d'Orsay, dem COGE (Choeurs et orchestres des grandes écoles, Paris), dem Orchestre du campus d'Orsay und dem Unihockeyclub Münster, die mich immer wieder auf den Boden zurückgeholt und auf andere Gedanken gebracht haben und so neu motiviert und inspiriert haben, und nicht zuletzt alle, die an mich geglaubt haben, und auch die, die es nicht getan haben.\par
Außer durch das Netzwerk "`Arithmetic Algebraic Geometry"' wurde die vorliegende Dissertation auch teilweise durch die Deutsche Forschungsgemeinschaft über den SFB 478 "`Geometrische Strukturen in der Mathematik"' unterstützt.
\chapter{Präliminarien}
 
In diesem Kapitel sollen einige für die Arbeit wichtige Begriffe zusammen gestellt werden.\\\par

\section{Torseure}

Es sei $S$ ein Schema. Dann ist ein $S$"=Gruppenschema $G$ ein $S$"=Schema $G$ zusammen mit einem Kompositionsgesetz $\gamma\colon h_{G}\times h_{G}\to h_{G}$, das ein Gruppengesetz ist. Dabei sei $h_{G}$ der Punktfunktor \[h_{G}\colon(Sch/S)\to(Ens), h_{G}(T):=G(T)=\Hom(T,G).\] Das Gruppengesetz $\gamma$ besteht aus einer Menge von Abbildungen \[\gamma_{T}\colon h_{G}(T)\times h_{G}(T)\to h_{G}(T)\] für $T\in(Sch/S)$, die verträglich mit Abbildungen $u\colon T'\to T$ sind, d.\,h. für die das Diagramm \begin{equation*}
        \xymatrix@=3em{%
        h_{G}(T)\times h_{G}(T) \ar[r]^{\gamma_{T}} \ar[d]^{h_{G}(u)\times h_{G}(u)} & h_{G}(T) \ar[d]^{h_{G}(u)} \\
        h_{G}(T')\times h_{G}(T') \ar[r]^{\gamma_{T'}} & h_{G}(T')}
           \end{equation*}
kommutiert, und es ist $h_{G}(T)$ eine Gruppe unter $\gamma_{T}$ für alle $S$"=Schemata $T$.\par
Diese Daten sind äquivalent zu einem Tupel $(G,m,\varepsilon,i)$ bestehend aus einem $S$"=Schema $G$, einer Multiplikation $m\colon G\times_{S}G\to G$, einem Einsschnitt \mbox{$\varepsilon\colon S\to G$} und einer inversen Abbildung $i\colon G\to G$, die jeweils Morphismen von $S$"=Schemata sind, so dass die folgenden Diagramme kommutieren:\par
\begin{itemize}
\item[a)] Assoziativität: \begin{equation*}
        \xymatrix@=3em{%
        G\times_{S}G\times_{S}G \ar[r]^{m\times\id_{G}} \ar[d]^{\id_{G}\times m} & G\times_{S}G \ar[d]^{m} \\
        G\times_{S}G \ar[r]^{m} & G}
           \end{equation*}
\item[b)] Existenz einer Linksidentität:
					\begin{equation*}
        \xymatrix@=3em{%
        G \ar[r]^{(p,\id_{G})} \ar[drr]^{\id_{G}} & S\times_{S}G \ar[r]^{\varepsilon\times\id_{G}} & G\times_{S}G \ar[d]^{m} \\
        & & G}
           \end{equation*}
\item[c)] Existenz eines Linksinversen:
	\begin{equation*}
        \xymatrix@=3em{%
        G \ar[r]^{(i,\id_{G})} \ar[d]^{p} & G\times_{S}G \ar[d]^{m} \\
        S \ar[r]^{\varepsilon} & G}
           \end{equation*}\\
\end{itemize}\par

Für eine beliebige gefaserte Kategorie $\mathcal{C}$ definiert man ein Gruppenobjekt vollkommen analog, man ersetze lediglich $(Sch/S)$ stets durch $\mathcal{C}$ und \mbox{"`$S$"=Schema"'} durch "`$S$"=Objekt"'.\\\par
Es sei nun wieder $S$ ein Schema und $G$ ein $S$"=Gruppenschema, das treuflach und lokal von endlicher Präsentation über $S$ ist. Dann definiert der Multiplikationsmorphismus $G\times_{S}G\to G$ von $G$ eine Aktion von $G$ auf sich selbst. Morphismen, Monomorphismen, Epimorphismen, Isomorphismen etc. von $S$"=Schemata mit einer $G$"=Aktion sind zudem in offensichtlicher Weise definiert.\par
Damit können wir nun den Begriff eines \emph{$G$"=Torseurs} definieren:\\\par

\begin{defn}\label{torsor}
Es sei $G$ ein Gruppenschema, das treuflach und lokal von endlicher Präsentation über dem Basisschema $S$ ist. Dann bezeichnet man als einen Rechts"=$G$"=Torseur $P$ bezüglich der fppf"=Topologie auf $S$ ein $S$"=Schema $P$, auf dem $G$ durch einen Morphismus $f\colon P\times_{S}G\to P, (x,g)\mapsto xg$ operiert, so dass die folgende Bedingung erfüllt ist:\par
Es gibt eine Überdeckung $(U_{i}\to S)_{i\in I}$ für die fppf"=Topologie auf $S$, so dass $P\times_{S}U_{i}$ mit seiner $G\times_{S}U_{i}$"=Aktion isomorph zu $G\times_{S}U_{i}$, versehen mit der durch den Multiplikationsmorphismus gegebenen Aktion von $G\times_{S}U_{i}$ auf sich selbst, für alle $i\in I$ ist.\\
\end{defn}\par

\begin{bsp}
Beispielsweise ist $G$ zusammen mit der durch den Multiplikationsmorphismus gegebenen Aktion auf sich selbst ein $G$"=Torseur (für die fppf"=Topologie).\par
Einen $G$"=Torseur $P$, der als $G$"=Torseur isomorph zu $G$ ist, bezeichnet man als trivialen $G$"=Torseur.\end{bsp}

Nach \cite[Proposition III.4.1]{M} hat man die folgende äquivalente Beschreibung eines $G$"=Torseurs für die fppf"=Topologie:\\\par

\begin{prop}
Es sei $G$ ein $S$"=Gruppenschema, das treuflach und lokal von endlicher Präsentation über $S$ ist. Ist $P$ ein $S$"=Schema, auf dem $G$ durch einen Morphismus $f\colon P\times_{S}G\to P, (x,g)\mapsto xg$ operiert, so sind die folgenden Bedingungen äquivalent:\par
\begin{itemize}
\item[a)] $P$ ist treuflach und lokal von endlicher Präsentation über $S$ und der Morphismus $P\times_{S}G\xrightarrow{(id_{P},f)}P\times_{S}P, (x,g)\mapsto(x,xg)$ ist ein Isomorphismus.
\item[b)] Es gibt eine Überdeckung $(U_{i}\to S)_{i\in I}$ für die fppf"=Topologie auf $S$, so dass $P\times_{S}U_{i}$ für jedes $i$ mit seiner Aktion von $G\times_{S}U_{i}$ isomorph zu dem trivialen $G\times_{S}U_{i}$"=Torseur $G\times_{S}U_{i}$ ist.
\end{itemize}
\end{prop}

Ebenso definiert man Garben-Torseure für die fppf-Topologie:\\\par

\begin{defn}
Es sei $G$ ein Gruppenobjekt in der Kategorie der Garben bezüglich der fppf"=Topologie auf $S$. Dann bezeichnet man als einen Rechts"=$G$"=(Garben-)Torseur $P$ bezüglich der fppf"=Topologie auf $S$ eine Garbe $P$ von Mengen bezüglich der fppf"=Topologie auf $S$, auf der $G$ durch einen Morphismus \[f\colon P\times_{S}G\to P, (x,g)\mapsto xg\] operiert, so dass die folgende Bedingung erfüllt ist:\par
Es gibt eine Überdeckung $(U_{i}\to S)_{i\in I}$ für die fppf"=Topologie auf $S$, so dass $P\times_{S}U_{i}$ für jedes $i$ mit seiner Aktion von $G\times_{S}U_{i}$ isomorph zu $G\times_{S}U_{i}$ ist, wobei $G\times_{S}U_{i}$ bzw. $P\times_{S}U_{i}$ die von $G$ bzw. $P$ induzierten Garben auf $(U_{i})_{fppf}$ bezeichnen und $G\times_{S}U_{i}$ zusammen mit der durch die Gruppenoperation induzierten Operation von $G\times_{S}U_{i}$ auf sich selbst betrachtet werde.\\
Als trivialen Torseur bezeichnet man einen $G$"=Torseur $P$, der als $G$"=Torseur isomorph zu $G$ ist, wobei $G$ die durch die Gruppenoperation gegebene Struktur eines $G$"=Torseurs trägt.\\
\end{defn}\par

Es ist offensichtlich, dass ein Garbentorseur $P$ für die fppf"=Topologie ein Torseur im Sinne von Definition \ref{torsor} ist, falls $P$ und $G$ durch Schemata darstellbar sind. Aus dem Yoneda"=Lemma folgt zudem, dass zwei Schemata $P$ und $P'$ isomorph als $G$"=Torseure sind, falls sie es als $G$"=Garbentorseure sind.\\\par

\begin{bem}\label{affin1}
Nicht jeder Garben"=Torseur $P$ in der fppf"=Topologie ist durch ein Schema $P$ darstellbar. Nach \cite[Theorem III.4.3a)]{M} gilt jedoch: Ist $G$ ein $S$"=Gruppenschema, das affin über $S$ ist, so ist jeder $G$"=Garben"=Torseur $P$ durch ein Schema $P$ darstellbar.\\  
\end{bem}\par

Torseure lassen sich auch allgemeiner für beliebige Topoi definieren:\\\par

\begin{defn}[\protect{\cite[Definition III.1.4.1, III.1.4.1.1]{G}}]
Es sei $T$ ein Topos und $S$ ein Objekt von $T$. Ein (Rechts-)Torseur von $T$ über $S$ ist ein Tupel $(P,G,m)$ bestehend aus einem $S$"=Objekt $P$, einem Gruppenobjekt $G$ von $T$ und einer Familie von Abbildungen \[m(S')\colon P(S')\times G(S')\to P(S')\] für alle $S$"=Objekte $S'$, so dass die folgenden Bedingungen erfüllt sind:\par 
\begin{itemize}
\item[i)] Der Morphismus $p\colon P\to S$ ist ein Epimorphismus, d.\,h. es existiert eine Familie von Epimorphismen $\{S_{i}\to S\}_{i\in I}$, so dass die Mengen \mbox{$P(S_{i})=\Hom_{S}(S_{i},P)$} nicht leer sind.
\item[ii)] Der Morphismus $u\colon P\times_{S}G\to P\times_{S}P,(p,g)\to(p,m(p,g))$ ist ein Isomorphismus.\\\end{itemize}
\end{defn}\par

\begin{bem}
Es seien $T$ ein Topos, $S$ ein Objekt von $T$ und $G$ ein Gruppenobjekt von $T$. Es sei $G_{d}$ das $G$"=Objekt, das man erhält, wenn man $G$ auf sich selbst durch Rechts\-translationen operieren lässt, d.\,h. es sei $G$ mit der von der Gruppenoperation $G\times_{S}G\to G$ induzierten Operation von $G$ auf sich selbst versehen.\par
Da nach Definition eines Gruppenobjekts der Morphismus $G_{d}\to S$ einen Schnitt besitzt, ist der Morphismus $G_{d}\to S$ ein Epimorphismus. Die zweite Bedingung für einen Torseur ist trivialerweise erfüllt, so dass also $G_{d}$ ein $G$"=Torseur ist.\par
Man bezeichnet $G_{d}$ als den \emph{trivialen Torseur}.\\
\end{bem}\par

Für den assoziierten Topos zu dem étalen Situs erhält man insbesondere:\\\par
 
\begin{defn}
Es sei $G$ ein Gruppenobjekt in der Kategorie der Garben bezüglich der étalen Topologie $T$ auf einem Schema $S$. Dann bezeichnet man als einen Rechts"=$G$"=(Garben-)Torseur $P$ bezüglich der étalen Topologie auf $S$ eine étale Garbe $P$ auf $S$, auf der $G$ durch einen Morphismus $f\colon P\times_{S}G\to P, (x,g)\mapsto xg$ operiert, so dass die folgende Bedingung erfüllt ist:\par
Es gibt eine Überdeckung $(U_{i}\to S)_{i\in I}$ für die étale Topologie auf $S$, so dass $P\times_{S}U_{i}$ für jedes $i$ mit seiner Aktion von $G\times_{S}U_{i}$ isomorph zu dem trivialen Torseur $G\times_{S}U_{i}$ ist, wobei $G\times_{S}U_{i}$ bzw. $P\times_{S}U_{i}$ die von $G$ bzw. $P$ induzierten Garben auf $(U_{i})_{\acute{e}t}$ bezeichnen.\\\end{defn}\par

Für beliebige Topoi gilt außerdem nach \cite[III.1.4.1.3]{G}:\\\par

\begin{lem}\label{pullback}
Es sei $f\colon T\to T'$ ein Morphismus von Topoi und $S$ ein Objekt von $T$ sowie $S'=f(S)$ das Bild von $S$ in $T'$ unter $f$. Ferner sei $G$ ein $S$"=Gruppenobjekt des Topos $T$ und $G'$ ein $S'$"=Gruppenobjekt des Topos $T'$. Ist dann $P$ ein $G$"=Torseur über $S$ in $T$, so ist $f_{*}(P)$ ein $f_{*}(G)$"=Torseur über $S'$ in $T'$. Ist umgekehrt $P'$ ein $G'$"=Torseur über $S'$ in $T'$, so ist $f^{*}P'$ ein  $f^{*}G'$"=Torseur über $S$ in $T$.\par
Diese Konstruktionen sind funktoriell.
\end{lem}
\bewub{\cite[III.1.4.1.3]{G}$\hfill\Box$}\\\par

Wichtig ist für uns noch der Zusammenhang zwischen Torseuren und Kohomologie:\\\par

\begin{sat}\label{koh}
Es seien $E$ ein Situs, $T$ der zu $E$ assoziierte Topos von abelschen Garben, $G$ ein Gruppenobjekt von $T$ und $S$ ein Objekt von $T$.\par
Dann definiert jeder $G$"=Torseur $P$ von $T$ über $S$ ein eindeutig bestimmtes Element $c(P)$ von $H^{1}(E,G)$, wobei die durch die Klasse des trivialen Torseurs punktierte Menge $H^{1}(E,*)$ die erste Rechtsableitung des Funktors \[H^{0}(E,*)\colon T\to Ens\] ist. Dabei ist $H^{0}(E,\mathscr{F})$ für eine Garbe $\mathscr{F}$ auf $E$ als \[H^{0}(E,\mathscr{F}):=\projlim_{U\in\Cat\,E^{0}}\mathscr{F}(U)\] (vgl. \cite[Notation II.2.2]{A}) definiert. Diese Zuordnung ist eine Bijektion.
\end{sat}
\bewub{\cite[Définition III.2.4.2, Remarque III.3.5.4]{G}$\hfill\Box$}\\\par

Aus diesem Satz folgt insbesondere:\\\par

\begin{kor}[vgl. auch \protect{\cite[Proposition III.4.6., Corollary III.4.7. ff]{M}}]
Ist $G$ ein flaches Gruppenschema von endlichem Typ über einem Schema $X$, so kann man die Menge der Isomorphieklassen der $G$"=Garbentorseure bezüglich der fppf"=Topologie (bzw. der Zariski"=Topologie bzw. der étalen Topologie) auf $X$ mit $\check{H}^{1}(X_{fppf},G)$ (bzw. $\check{H}^{1}(X,G)$ bzw. $\check{H}^{1}_{\acute{e}t}(S,G)$) identifizieren.\\
Ist $G$ zusätzlich affin, so kann man die Menge der Isomorphieklassen der $G$"=Torseure bezüglich der fppf"=Topologie (bzw. der Zariski"=Topologie bzw. der étalen Topologie) auf $X$ mit $\check{H}^{1}(X_{fppf},G)$ (bzw. $ \check{H}^{1}(X,G)$ bzw. $\check{H}^{1}_{\acute{e}t}(S,G)$) identifizieren.\\
\end{kor}
\par

Abschließend sei noch erwähnt, dass man für die fppf"=Topologie auf einem Schema $S$ den Begriff eines $G$"=Torseurs auch für beliebige abstrakte endliche Gruppen definieren kann:\\\par

\begin{defn}
Es sei $S$ ein Schema und $G$ eine abstrakte endliche Gruppe. Dann bezeichnet man ein Schema $P$ über $S$ als einen $G$"=Torseur bezüglich der fppf"=Topologie auf $S$, wenn die folgenden Bedingungen erfüllt sind:\par
\begin{itemize}
\item[a)] Der Strukturmorphismus $p\colon P\to S$ ist treuflach und lokal von endlicher Präsentation.
\item[b)] Der Morphismus $P\times_{S}G_{P}\xrightarrow{(id_{P},f)}P\times_{S}P, (x,g)\mapsto(x,xg)$ ist ein Isomorphismus. Dabei bezeichnet $G_{P}$ das Schema $G_{P}=\coprod_{g\in G}P_{g}$ mit $P_{g}=P$ für alle $g\in G$.\\
\end{itemize}
\end{defn}

Es sei außerdem noch an den Begriff einer Galoisüberlagerung erinnert:\\\par

\begin{defn}[\protect{\cite[Remark I.5.4]{M}}]
Es sei $G$ eine endliche Gruppe, die auf einem Schema $Y$ über $X$ operiert. Dann nennt man $Y$ eine Galoisüberlagerung von $X$ mit Gruppe $G$, falls $Y$ treuflach über $X$ ist und der kanonische Morphismus \[\varphi\colon G_{Y}\to Y\times_{X}Y,\] der durch \[\varphi\mathop{|}\nolimits_{Y_{\sigma}}=(y\mapsto(y,y\sigma))\]  für alle $g\in G$ gegeben ist, ein Isomorphismus ist.\\
\end{defn}\par

\begin{bem}\label{galo}
Nach \cite[Lemma 4.4.1.8]{Mu} wird jeder endliche étale Morphismus \mbox{$f\colon Y\to X$} durch eine Galoisüberdeckung $Y'\to X$ dominiert.\\
\end{bem}\par

Damit erhalten wir die folgende Charakterisierung von $G$"=Torseuren:\\\par

\begin{lem}
Es sei $k$ ein algebraisch abgeschlossener Körper, $G$ ein Gruppenschema über $k$ und $X$ ein beliebiges Schema über $k$. Dann entsprechen die Isomorphieklassen der $G$"=Torseure bezüglich der étalen Topologie auf $X$, die nach Pullback unter einer fest gewählten Galoisüberdeckung $X'\to X$ mit Gruppe $\mathfrak{g}$ trivial werden, bijektiv den Elementen von $H^{1}(\mathfrak{g},\Gamma(X',G))$.
\end{lem}
\bewub{Die \v{C}ech"=Kohomologie für die étale Topologie liefert die folgende exakte Sequenz:\par
\begin{equation*}
\xymatrix@=3em{%
1 \ar[r] & \check{H}^{1}(X'/X,G) \ar[r] & H^{1}_{\acute{e}t}(X,G) \ar[r] & H^{1}_{\acute{e}t}(X',G) \ar[r] & 1.}\end{equation*}
Nach \cite[Example III.2.6]{M} ist aber $\check{H}^{1}(X'/X,G)$ isomorph zu $H^{1}(\mathfrak{g},\Gamma(X',G))$, so dass die Behauptung folgt. \hfill$\Box$
     }\\\par

\section{Verschiedene Kategorien von Überdeckungen}

Für die Definition einer gewissen Kategorie $\mathscr{B}_{\mathfrak{X}_{\mathfrak{o}},D}(G)$ von $G$"=Torseuren auf einem Schema $\mathfrak{X}_{\mathfrak{o}}$ und ihre Charakterisierung im zweiten und dritten Kapitel benötigen wir verschiedene Kategorien von Überdeckungen aus \cite{DW1}.\par
Wie in \cite{DW1} verstehen wir hier und in der gesamten Arbeit unter einer Varietät über einem Körper $k$ ein geometrisch irreduzibles und geometrisch reduziertes separiertes Schema von endlichem Typ über $k$. Eine Kurve ist eine eindimensionale Varietät.\par
Es sei im folgenden $R$ ein Bewertungsring mit Quotientenkörper $Q$, dessen Charakteristik 0 sei.\par
Zu einer glatten projektiven Kurve $X$ über $Q$ betrachten wir jeweils ein Modell $\mathfrak{X}$ von $X$ über $R$. Darunter sei stets ein endlich präsentiertes, flaches und eigentliches Schema $\mathfrak{X}$ über $\Spec\,R$ verstanden, so dass $X$ isomorph zu $\mathfrak{X}\otimes_{R}Q$ ist. Außerdem schreiben wir für einen Divisor $D$ auf $X$ immer kurz $X-D$ anstelle von $X-\supp\,D$.\\\par

\begin{defn}[\protect{\cite[S. 4,5]{DW1}}]
Die Kategorie $S_{\mathfrak{X},D}$ ist wie folgt definiert:\par
Die Objekte der Kategorie sind endlich präsentierte, eigentliche $R$"=Morphismen \mbox{$\pi\colon\mathcal{Y}\to\mathfrak{X}$}, deren generische Faser $\pi_{Q}\colon\mathcal{Y}_{Q}\to X$ endlich ist und für welche die induzierte Abbildung $\pi_{Q}\colon\pi_{Q}^{-1}(X-D)\to X-D$ étale ist.\par
Ein Morphismus von $\pi_{1}\colon\mathcal{Y}_{1}\to\mathfrak{X}$ nach $\pi_{2}\colon\mathcal{Y}_{2}\to\mathfrak{X}$ ist ein Morphismus $\varphi\colon\mathcal{Y}_{1}\to\mathcal{Y}_{2}$, so dass $\pi_{1}=\pi_{2}\circ\varphi$ ist.\par
Existiert solch ein Morphismus von $\pi_{1}\colon\mathcal{Y}_{1}\to\mathfrak{X}$ nach $\pi_{2}\colon\mathcal{Y}_{2}\to\mathfrak{X}$, so sagt man, dass $\pi_{1}$ den Morphismus $\pi_{2}$ dominiert. Induziert zusätzlich $\varphi_{Q}$ einen Isomorphismus der lokalen Ringe in zwei generischen Punkten, so sagt man, dass $\pi_{1}$ den Morphismus $\pi_{2}$ strikt dominiert.\\
\end{defn}\par

Ist $R'$ ein Bewertungsring (mit Quotientenkörper $Q'$ der Charakteristik Null), der $R$ enthält, so hat man offensichtlich einen Basiswechselfunktor \[S_{\mathfrak{X},D}\to S_{\mathfrak{X}\otimes_{R}R',D'},\] wobei $D'$ das Urbild von $D$ in $\mathfrak{X}\otimes_{R}R'$ ist. Klar ist auch, dass endliche Produkte und endliche gefaserte Produkte in $S_{\mathfrak{X},D}$ existieren. Außerdem halten wir fest:\\\par

\begin{bem}
Es sei $\varphi\colon(\pi_{1}\colon\mathcal{Y}_{1}\to\mathfrak{X})\to(\pi_{2}\colon\mathcal{Y}_{2}\to\mathfrak{X})$ ein Morphismus in $S_{\mathfrak{X},D}$. Dann gilt:\par
\begin{itemize}
\item[i)] $\varphi\colon\mathcal{Y}_{1}\to\mathcal{Y}_{2}$ ist eigentlich nach \cite[Proposition 3.3.16]{L1}, da $\pi_{1}$ und $\pi_{2}$ jeweils insbesondere eigentlich sind.
\item[ii)] $\varphi\colon\mathcal{Y}_{1}\to\mathcal{Y}_{2}$ ist von endlicher Präsentation nach \cite[Proposition 1.6.2.v)]{EGA IV}, da $\pi_{2}$ per Definition insbesondere quasi"=separiert und von endlicher Präsentation ist.
\item[iii)] $\varphi_{Q}$ ist endlich über $X-D$ nach \cite[Lemma 3.3.15]{L1}, da nach Voraussetzung $\pi_{2}$ und damit $(\pi_{2})_{Q}$ eigentlich ist.
\item[iv)]  $\varphi_{Q}$ ist étale über $X-D$ nach \cite[Corollary I.3.6]{M}, da nach Voraussetzung $(\pi_{1})_{Q}$ und $(\pi_{2})_{Q}$ étale über $X-D$ sind.
\end{itemize}
\end{bem}

Wir benötigen außerdem folgende volle Unterkategorie $S_{\mathfrak{X},D}^{good}$ von $S_{\mathfrak{X},D}$:\\\par

\begin{defn}
Die volle Unterkategorie $S_{\mathfrak{X},D}^{good}$ von $S_{\mathfrak{X},D}$ besteht aus den Elementen $\pi\colon\mathcal{Y}\to\mathfrak{X}$ von $S_{\mathfrak{X},D}$ mit der Eigenschaft, dass der Strukturmorphismus $\lambda\colon\mathcal{Y}\to\Spec\,R$ flach ist, $\lambda_{*}\mathcal{O}_{\mathcal{Y}}=\mathcal{O}_{\Spec\,R}$ universell erfüllt ist und ihre generische Faser \linebreak\mbox{$\lambda_{Q}\colon\mathcal{Y}_{Q}\to\Spec\,Q$} glatt ist.\\
\end{defn}\par

\begin{bem}\label{zusammen}
Es sei $(\pi\colon\mathcal{Y}\to\mathfrak{X})\in S_{\mathfrak{X},D}^{good}$. Dann ist $\mathcal{Y}_{Q}$ geometrisch zusammenhängend.
\end{bem}
\bewub{Es sei $(\pi\colon\mathcal{Y}\to\mathfrak{X})\in S_{\mathfrak{X},D}^{good}$ und $\overline{Q}$ ein algebraischer Abschluss von $Q$. Da $\mathcal{Y}_{\overline{Q}}$ glatt über $\overline{Q}$ ist, ist $\mathcal{Y}_{\overline{Q}}$ genau dann zusammenhängend, wenn $\mathcal{Y}_{\overline{Q}}$ irreduzibel ist.\par
Es genügt also zu zeigen, dass $\mathcal{Y}_{\overline{Q}}$ irreduzibel ist.\par
Dafür genügt es zu zeigen, dass $H^{0}(\mathcal{Y}_{\overline{Q}},\mathcal{O}_{\mathcal{Y}_{\overline{Q}}})$ eindimensional über $\overline{Q}$ ist: Denn wäre $\mathcal{Y}_{\overline{Q}}$ nicht irreduzibel, also etwa $\mathcal{Y}_{\overline{Q}}=A\coprod B$ mit nichtleeren $A$ und $B$, so wäre $H^{0}(\mathcal{Y}_{\overline{Q}},\mathcal{O}_{\mathcal{Y}_{\overline{Q}}})=H^{0}(A,\mathcal{O}_{\mathcal{Y}_{\overline{Q}}})\oplus H^{0}(B,\mathcal{O}_{\mathcal{Y}_{\overline{Q}}})$. Da $H^{0}(A,\mathcal{O}_{\mathcal{Y}_{\overline{Q}}})$ und $H^{0}(B,\mathcal{O}_{\mathcal{Y}_{\overline{Q}}})$ jeweils mindestens eindimensional über $\overline{Q}$ sind, wäre dann $H^{0}(\mathcal{Y}_{\overline{Q}},\mathcal{O}_{\mathcal{Y}_{\overline{Q}}})$ mindestens zweidimensional über $\overline{Q}$.\par
Wir zeigen also, dass $H^{0}(\mathcal{Y}_{\overline{Q}},\mathcal{O}_{\overline{Q}})$ eindimensional über $\overline{Q}$ ist:\par
Da $\lambda_{*}\mathcal{O}_{\mathcal{Y}}=\mathcal{O}_{\Spec\,R}$ nach Definition von $S_{\mathfrak{X},D}^{good}$ universell gilt, gilt \linebreak \mbox{$(\lambda_{\overline{Q}})_{*}\mathcal{O}_{\mathcal{Y}_{\overline{Q}}}=\mathcal{O}_{\Spec\,\overline{Q}}$}. Daraus folgt, dass \[H^{0}(\mathcal{Y}_{\overline{Q}},\mathcal{O}_{\mathcal{Y}_{\overline{Q}}})\cong H^{0}(\Spec\,\overline{Q},(\lambda_{\overline{Q}})_{*}\mathcal{O}_{\mathcal{Y}_{\overline{Q}}})\cong H^{0}(\Spec\,\overline{Q}, \mathcal{O}_{\Spec\,\overline{Q}})=\overline{Q}\] ist,
also eindimensional über $\overline{Q}$ ist. $\hfill\Box$}\\\par

\begin{kor}\label{irred}
Es sei $(\pi\colon\mathcal{Y}\to\mathfrak{X})\in S_{\mathfrak{X},D}^{good}$. Dann ist $\mathcal{Y}_{Q}$ eine glatte und projektive Kurve und also $\mathcal{Y}$ irreduzibel und reduziert.
\end{kor}
\bewub{Es sei $(\pi\colon\mathcal{Y}\to\mathfrak{X})\in S_{\mathfrak{X},D}^{good}$. Dann ist nach Bemerkung \ref{zusammen} $\mathcal{Y}_{Q}$ geometrisch zusammenhängend und daher, wie in deren Beweis gesehen, geometrisch irreduzibel.\par
Da $X$ nach Voraussetzung eine Kurve über $Q$ ist, ist $X$ insbesondere separiert über $Q$. Außerdem ist $\mathcal{Y}_{Q}$ nach Definition der Kategorie $S_{\mathfrak{X},D}$ endlich über $X$, also affin nach \cite[6.1.2]{EGA II} und damit separiert nach \cite[Proposition 1.2.4]{EGA II}, so dass insgesamt der Morphismus $\mathcal{Y}_{Q}\to Q$ als Komposition zweier separierter Morphismen separiert ist.\par
Die Glattheit von $\mathcal{Y}_{Q}$ über $Q$ ergibt sich direkt aus der Definition von $S_{\mathfrak{X},D}^{good}$.\par
Da $\overline{Q}$ als ein Körper nach \cite[17.E]{Ma} ein regulärer lokaler Ring ist und damit auch $\Spec\,\overline{Q}$ regulär ist (man benutze dabei \cite[Theorem 4.2.16]{L1}) sowie $\mathcal{Y}_{\overline{Q}}$ glatt über $\overline{Q}$ ist, ist nach \cite[Theorem 4.3.36]{L1} auch $\mathcal{Y}_{\overline{Q}}$ regulär. Also sind alle lokalen Ringe von $\mathcal{Y}_{\overline{Q}}$ regulär und damit nach \cite[Proposition 4.2.11]{L1} auch reduziert. Damit ist $\mathcal{Y}_{\overline{Q}}$ reduziert und also $\mathcal{Y}_{Q}$ geometrisch reduziert.\par
Weil der Morphismus $\pi\colon\mathcal{Y}\to\mathfrak{X}$ nach Definition von $S_{\mathfrak{X},D}$ von endlicher Präsentation ist, ist auch $\pi_{Q}\colon\mathcal{Y}_{Q}\to\mathfrak{X}_{Q}$ von endlicher Präsentation. Da $\mathfrak{X}$ ein Modell von $X$ und damit insbesondere von endlicher Präsentation über $R$ ist, ist außerdem $\mathfrak{X}_{Q}$ von endlicher Präsentation über $Q$. Also ist der Morphismus $\mathcal{Y}_{Q}\to Q$ als Komposition zweier Morphismen von endlicher Präsentation ebenfalls von endlicher Präsentation, insbesondere damit von endlichem Typ.\par
Insgesamt ist damit $\mathcal{Y}_{Q}$ eine glatte Varietät über $Q$.\par
Desweiteren ist $\mathcal{Y}_{Q}$ nach dem bisher gezeigten geometrisch irreduzibel, nach \cite[Remark 3.2.11]{L1} also irreduzibel. Aus der Endlichkeit von $\pi_{Q}$ folgt dann, dass der Funktionenkörper $K(\mathcal{Y}_{Q})$ eine endliche algebraische Erweiterung von $K(X)$ ist, so dass $\trdeg_{Q}K(\mathcal{Y}_{Q})=1$ und also nach \cite[Proposition 7.4.1]{EGA II} $\mathcal{Y}_{Q}$ eindimensional ist.\par
Also ist $\mathcal{Y}_{Q}$ eine glatte Kurve über $Q$.\par
Weil $\pi\in S_{\mathfrak{X},D}$ ist, ist ferner der Morphismus $\pi_{Q}$ endlich und damit nach \cite[Corollaire 6.11]{EGA II} projektiv. Da nach Voraussetzung $X$ projektiv über $Q$ ist, ist also der Strukturmorphismus $\mathcal{Y}_{Q}\to Q$ als Komposition zweier projektiver Morphismen projektiv.\par
Es bleibt noch zu zeigen, dass $\mathcal{Y}$ irreduzibel und reduziert ist:\par
Da nach dem bisher gezeigten $\mathcal{Y}_{Q}$ geometrisch reduziert ist, ist nach \cite[Remark 3.2.11]{L1} $\mathcal{Y}_{Q}$ auch reduziert; dass $\mathcal{Y}_{Q}$ irreduzibel ist, wurde bereits gezeigt. Da außerdem $(\pi\colon\mathcal{Y}\to\mathfrak{X})\in S_{\mathfrak{X},D}^{good}$ ist, ist $\mathcal{Y}$ flach und von endlicher Präsentation über $R$, so dass nach \cite[Proposition 4.3.8]{L1} $\mathcal{Y}$ irreduzibel und reduziert ist. \hfill$\Box$}\\\par

Ein im folgenden wichtiges Resultat ist das folgende Theorem:\\\par
\newpage

\begin{thm}[\protect{\cite[Theorem 1, 3)-5)]{DW1}}]\label{Kgood}
Es sei $R$ ein diskreter Bewertungsring (mit Quotientenkörper $Q$ der Charakteristik Null). Dann gilt:\par
\begin{enumerate}
\item[(i)] Für jeden diskreten Bewertungsring $R'$, der über $R$ liegt, setze man \linebreak \mbox{$\mathfrak{X'}:=\mathfrak{X}\otimes_{R}R'$} und definiere $D'$ als das Urbild von $D$ in $X'$. Dann bildet der kanonische Basiswechselfunktor $S_{\mathfrak{X},D}\to S_{\mathfrak{X'},D'}$ die Unterkategorie $S_{\mathfrak{X},D}^{good}$ auf $S_{\mathfrak{X'},D'}^{good}$ ab.
\item[(ii)] Für je endlich viele Objekte $\pi_{i}\colon\mathcal{Y}_{i}\to\mathfrak{X}$ aus $S_{\mathfrak{X},D}$ existiert eine endliche Erweiterung $Q'$ von $Q$, so dass es ein Objekt $\pi'\in S_{\mathfrak{X'},D'}^{good}$ gibt, das alle Objekte $\pi_{i}\otimes_{R}R'\in S_{\mathfrak{X'},D'}$  dominiert. Dabei ist $R'$ ein diskreter Bewertungsring in $Q'$, der $R$ dominiert, $\mathfrak{X'}:=\mathfrak{X}\otimes_{R}R'$ und $D'$ das Urbild von $D$ in $X'$.
\item[(iii)] Zu jedem Objekt $\pi\colon\mathcal{Y}\to\mathfrak{X}$ aus $S_{\mathfrak{X},D}$ existiert eine Erweiterung von diskreten Bewertungsringen $R'/R$ wie in $(ii)$, so dass $\pi\otimes_{R}R'$ von 
einem Objekt aus $S_{\mathfrak{X'},D'}^{good}$ strikt dominiert wird.
\end{enumerate}
\end{thm}
\bewub{\begin{enumerate} \item[(i)] Dies ist klar, da Glattheit und Flachheit stabil unter Basiswechsel sind.
\item[(ii)] Weil endliche Produkte in der Kategorie $S_{\mathfrak{X},D}$ existieren, folgt Aussage $(ii)$ aus Aussage $(iii)$.
\item[(iii)] Nach \cite[Proposition 3.2.20]{L1} gibt es einen $Q^{sep}$"=wertiger Punkt von \mbox{$X-D$}. Da $\pi^{-1}_{Q}(X-D)$ endlich und étale über $X-D$ ist, existiert somit nach \cite[4.2]{Mu} ein $Q^{sep}$"=rationaler Punkt von $\pi_{Q}^{-1}(X-D)$ über dem gegebenen Punkt von $X-D$. Mittels noetherschem Descent (vgl. \cite[§8]{EGA IV}) erhält man daraus eine endliche Erweiterung $Q'$ von $Q$, so dass $\mathcal{Y}$ einen $Q'$"=rationalen Punkt über $X-D$ besitzt, wobei man $Q'$ so wählen kann, dass die irreduziblen Komponenten von $\mathcal{Y}_{Q'}$ geometrisch irreduzibel sind.\par
Es sei $R'$ ein diskreter Bewertungsring in $Q'$ über $R$ und $\mathcal{Y}_{R'}:=\mathcal{Y}\otimes_{R}R'$. Wählt man nun eine irreduzible Komponente von $\mathcal{Y}_{Q'}$, die einen $Q'$"=wertigen Punkt über $X-D$ besitzt, und definiert $\mathcal{Y}^{*}$ als deren Abschluss in $\mathcal{Y}_{R'}$, der mit der Struktur eines reduzierten Schemas versehen sei, so ist $\mathcal{Y}^{*}$ reduziert und irreduzibel sowie nach \cite[§7.8]{EGA IV} auch exzellent, so dass seine Normalisierung $\widetilde{\mathcal{Y}}$ nach \cite[Proposition 7.8.6 (ii)]{EGA IV} endlich über $\mathcal{Y}^{*}$ ist. Außerdem ist $\widetilde{\mathcal{Y}}$ ein eigentliches und flaches $R'$"=Schema, da $\mathcal{Y}_{R'}$ dies ist. Da die Bildung der Normalisierung mit der des Faserprodukts vertauscht, ist $\widetilde{\mathcal{Y}}\otimes_{R'}Q'$ die Normalisierung von $\mathcal{Y}^{*}_{Q}$. Also besitzt $\widetilde{\mathcal{Y}}\otimes_{R'}Q'$ einen $Q'$"=rationalen Punkt. Man beachte, dass nach Konstruktion $\widetilde{\mathcal{Y}}$ von Dimension 2 ist. Nach Lipmans Auflösung von Singularitäten (siehe \cite[S. 362]{L1}) existiert dann ein irreduzibles und reguläres $R'$"=Schema $\mathcal{Y}^{V}$ zusammen mit einem eigentlichen $R'$"=Morphismus $\mathcal{Y}^{V}\to\widetilde{\mathcal{Y}}$, der ein Isomorphismus auf der generischen Faser ist. Das Schema $\mathcal{Y}^{V}$ erhält man dabei durch wiederholte Aufblasung des singulären Orts gefolgt von der Bildung der Normalisierung. Dieser Prozess wird nach endlich vielen Schritten stationär,
so dass man ein reguläres irreduzibles Schema $\mathcal{Y}^{V}$ erhält, das eigentlich und flach über $R'$ ist, zusammen mir einem eigentlichen Morphismus $\mathcal{Y}^{V}\to\mathfrak{X'}$, der nach Konstruktion den Morphismus $\pi\otimes_{R}\id_{R'}$ strikt dominiert.\par
Dies alles impliziert, dass der Morphismus $\psi\colon\mathcal{Y}^{V}\to\mathfrak{X'}$ in der Kategorie $S_{\mathfrak{X'},D'}$ liegt und auch der Strukturmorphismus $\lambda^{V}\colon\mathcal{Y}^{V}\to\Spec\,R'$ flach sowie dessen generische Faser glatt ist.\par
Um zu sehen, dass $\psi$ sogar in $S_{\mathfrak{X'},D'}^{good}$ liegt, bleibt also nur noch zu zeigen, dass $\lambda^{V}_{*}\mathcal{O}_{\mathcal{Y}^{V}}=\mathcal{O}_{\Spec\,R'}$ universell gilt.\par
Aufgrund der Eigentlichkeit von $\mathcal{Y}^{V}$ induziert aber der $Q'$"=rationale Punkt der generischen Faser von $\mathcal{Y}^{V}$ einen Schnitt des Strukturmorphismus $\lambda^{V}\colon\mathcal{Y}^{V}\to\Spec\,R'$. Damit folgt nun aus einem Theorem von Raynaud (siehe \cite[Théorème 8.2.1, $(ii)\Rightarrow(iv)$]{Ray} oder \cite[Corollary 9.1.24., Corollary 9.1.32.]{L1}), dass $\mathcal{Y}^{V}$ kohomologisch flach in Dimension 0 ist, also insbesondere $\lambda^{V}_{*}\mathcal{O}_{\mathcal{Y}^{V}}=\mathcal{O}_{\Spec\,R'}$ universell gilt. \hfill$\Box$
\end{enumerate}}\par
\text{}\\\par

Für den konkreten Fall $R=\overline{\mathbb{Z}}_{p}$ und $Q=\overline{\mathbb{Q}}_{p}$ gilt:\\\par

\begin{kor}[\protect{\cite[Corollary 3,1)]{DW1}}]\label{kgood}
Es sei $X$ eine glatte und projektive Kurve über $\overline{\mathbb{Q}}_{p}$, $D$ ein Divisor auf $X$ sowie $\mathfrak{X}$ ein Modell von $X$ über $\Spec\,\overline{\mathbb{Z}}_{p}$. Außerdem sei eine endliche Anzahl von Objekten $\pi_{i}\colon\mathcal{Y}_{i}\to\mathfrak{X}$ aus $S_{\mathfrak{X},D}$ gegeben.\par
Dann existieren eine endliche Erweiterung $K$ von $\mathbb{Q}_{p}$ und eine glatte und projektive Kurve $X_{K}$ über $K$ mit Modell $\mathfrak{X}_{\mathfrak{o}_{K}}$ über $\mathfrak{o}_{K}$ sowie ein Divisor $D_{K}$ von $X_{K}$, so dass folgendes gilt: Es ist $X=X_{K}\otimes_{K}\overline{K}$ und $D=D_{K}\otimes_{K}\overline{K}$ sowie $\mathfrak{X}=\mathfrak{X}_{\mathfrak{o}_{K}}\otimes_{\mathfrak{o}_{K}}\overline{\mathbb{Z}}_{p}$. Außerdem existiert ein Element $\pi_{\mathfrak{o}_{K}}\colon\mathcal{Y}_{\mathfrak{o}_{K}}\to\mathfrak{X}_{\mathfrak{o}_{K}}$ von $S_{\mathfrak{X}_{\mathfrak{o}_{K}},D_{K}}^{good}$, so dass der Morphismus $\pi=\pi_{\mathfrak{o}_{K}}\otimes_{\mathfrak{o}_{K}}\overline{\mathbb{Z}}_{p}$ alle Morphismen $\pi_{i}$ dominiert.
\end{kor}
\bewub{Da $\mathfrak{X}$ als Modell von $X$ insbesondere von endlicher Präsentation über $\overline{\mathbb{Z}}_{p}$ ist sowie nach \cite[8.1.2 b)]{EGA IV} $\indlim\limits_{\mathbb{Q}_{p}\subset K\,\text{endlich}}K=\overline{\mathbb{Q}}_{p}$ und damit auch $\indlim\limits_{\mathbb{Q}_{p}\subset K\,\text{endlich}}\mathfrak{o}_{K}=\overline{\mathbb{Z}}_{p}$ gilt, existieren nach \cite[Proposition 8.9.1 (i)]{EGA IV} eine endliche Erweiterung $K_{0}$ von $\mathbb{Q}_{p}$ sowie ein Schema $\mathfrak{X'}_{0}$ von endlichem Typ über $\mathfrak{o}_{K_{0}}$, so dass es einen $\overline{\mathbb{Z}}_{p}$"=Isomorphismus $\mathfrak{X'}_{0}\otimes_{\mathfrak{o}_{K_{0}}}\overline{\mathbb{Z}}_{p}\xrightarrow{\sim}\mathfrak{X}$ gibt. Da $\mathfrak{o}_{K_{0}}$ noethersch ist, ist nach \cite[Définition 1.6.1]{EGA IV} $\mathfrak{X'}_{0}$ auch von endlicher Präsentation über $K_{0}$. Nach \cite[Théorème 8.10.5 (xii)]{EGA IV} und \cite[Théorème 11.2.6 (ii)]{EGA IV} existiert dann eine endliche Erweiterung $K_{1}$ von $\mathbb{Q}_{p}$, die $K_{0}$ enthält, so dass $\mathfrak{X}_{1}:=\mathfrak{X'}_{0}\otimes_{\mathfrak{o}_{K_{o}}}\mathfrak{o}_{K_{1}}$ von endlicher Präsentation, flach und eigentlich über $\mathfrak{o}_{K_{1}}$ ist. Setzt man nun $X_{1}:=\mathfrak{X}_{1}\otimes_{\mathfrak{o}_{K_{1}}}K_{1}$, so folgt aus \cite[Scholie 8.8.3]{EGA IV}, dass $X$ zu $X_{1}$ absteigt. Nach \cite[Proposition 8.9.1]{EGA IV}, \cite[Théorème 8.10.5 (v),(xiii)]{EGA IV} und \cite[Proposition 17.7.8]{EGA IV} sowie \cite[Proposition 3.2.7 (a)]{L1} gibt es also eine endliche Erweiterung $K_{2}$ von $\mathbb{Q}_{p}$ mit $K_{1}\subset K_{2}$, so dass $X_{2}:=X_{1}\otimes_{K_{1}}K_{2}$ eine glatte und projektive Kurve über $K_{2}$ ist.
Da die Eigenschaften von $\mathfrak{X}_{1}$ stabil unter Basiswechsel sind, ist dann $\mathfrak{X}_{2}:=\mathfrak{X}_{1}\otimes_{\mathfrak{o}_{K_{1}}}\mathfrak{o}_{K_{2}}$ ein Modell von $X_{2}$.  \par
Nach \cite[Proposition 8.5.5]{EGA IV} gibt es außerdem eine endliche Erweiterung $K_{3}$ von $K_{2}$, so dass der Divisor $D$ zu einem Divisor $D_{3}$ von $X_{3}$ über $\Spec\,K_{3}$ absteigt, wobei man dazu den zu $D$ assoziierten $\mathcal{O}_{X}$"=Modul $\mathcal{O}_{X}(D)$ betrachte.\par
Da es endlich viele Morphismen $\pi_{i}$ sind, die gegeben sind, gibt es nach \cite[Proposition 8.9.1, Théorème 8.10.5 (xii)]{EGA IV} eine endliche Erweiterung $K_{4}$ von $K_{3}$, so dass alle Morphismen \mbox{$\pi_{i}\colon\mathcal{Y}_{i}\to\mathfrak{X}$} zu endlich"=präsentierten, eigentlichen Morphismen $\pi_{i4}\colon\mathcal{Y}\to\mathfrak{X}_{4}$ absteigen. Da noetherscher Descent verträglich mit der Bildung der generischen Faser ist, existiert dann nach \cite[Théorème 8.10.5 (x)]{EGA IV} und \cite[Proposition 17.7.8]{EGA IV} eine endliche Erweiterung $K_{5}$ von $K_{4}$, so dass alle $(\pi_{i5})_{Q}$ endlich sind und zudem étale über $X_{5}-D_{5}$ sind. Damit liegen alle Morphismen $\pi_{i5}$ in der Kategorie $S_{\mathfrak{X}_{5},D_{5}}$.\par
Nach Theorem \ref{Kgood} existiert also eine Erweiterung $K'_{5}$ von $K_{5}$, so dass die Objekte $\pi_{i5}\otimes_{\mathfrak{o}_{K_{5}}}R'_{5}$ alle von einem einzigen Objekt $\pi_{5}\in S_{\mathfrak{X}_{5}\otimes_{\mathfrak{o}_{K_{5}}}R'_{5},D_{5}\otimes_{\mathfrak{o}_{K_{5}}}R'_{5}}^{good}$ dominiert werden, wobei $R'_{5}$ ein Bewertungsring in $Q'_{5}$ ist, der $\mathfrak{o}_{K_{5}}$ dominiert. Setzt man nun $K:=K'_{5}$, so besitzen $X_{K}$, $\mathfrak{X}_{\mathfrak{o}_{K}}:=\mathfrak{X}_{5}\otimes_{\mathfrak{o}_{K}}R'_{5}$, $D_{K}$ und $\pi_{\mathfrak{o}_{K}}:=\pi_{5}$ alle gewünschten Eigenschaften. 
\hfill$\Box$}
\\\par

Daraus folgt:\\\par

\begin{kor}[\protect{\cite[Corollary 3, 3)]{DW1}}]\label{dom}
Es sei $X$ eine glatte und projektive Kurve über $\overline{\mathbb{Q}}_{p}$ und $D$ ein Divisor auf $X$ sowie $\mathfrak{X}$ ein Modell von $X$ über $\Spec\,\overline{\mathbb{Z}}_{p}$. Dann gilt:\par
Jede endliche Anzahl von Objekten $\pi_{i}\colon\mathcal{Y}_{i}\to\mathfrak{X}$ in $S_{\mathfrak{X},D}$ wird von einem gemeinsamen Objekt $(\pi\colon\mathcal{Y}\to\mathfrak{X})\in S_{\mathfrak{X},D}^{good}$ dominiert.
\end{kor}
\bewub{Nach Korollar \ref{kgood} existieren eine endliche Erweiterung $K$ von $\mathbb{Q}_{p}$, eine glatte und projektive Kurve $X_{K}$ über $K$ mit einem Modell $\mathfrak{X}_{\mathfrak{o}_{K}}$ über $\mathfrak{o}_{K}$ und ein Divisor $D_{K}$ von $X_{K}$, so dass $X=X_{K}\otimes\overline{K}$, $D=D_{K}\otimes\overline{K}$ und $\mathfrak{X}=\mathfrak{X}_{\mathfrak{o}_{K}}\otimes_{\mathfrak{o}_{K}}\overline{\mathbb{Z}_{p}}$ gilt, sowie zusätzlich ein Element $\pi_{\mathfrak{o}_{K}}\in S^{good}_{X_{\mathfrak{o}_{K}},D_{K}}$, so dass der Morphismus $\pi:=\pi_{\mathfrak{o}_{K}}\otimes_{\mathfrak{o}_{K}}\overline{\mathbb{Z}_{p}}$ alle Morphismen $\pi_{i}$ dominiert. Dieser Morphismus $\pi$ liegt in der Kategorie $S_{\mathfrak{X},D}^{good}$, da die Definition der vollen Unterkategorie $S_{\mathfrak{X},D}^{good}$ nach Theorem \ref{Kgood}(i) stabil unter der Erweiterung des zugrundeliegenden Bewertungsrings ist. 
\hfill$\Box$}\\\par

Außer der vollen Unterkategorie $S_{\mathfrak{X},D}^{good}$ betrachten wir noch die folgende volle Unterkategorie von $S_{\mathfrak{X},D}$:\\\par

\begin{defn}[vgl. \protect{\cite[S. 5]{DW1}}]
Die volle Unterkategorie $S_{\mathfrak{X},D}^{ss}$ von $S_{\mathfrak{X},D}$ besteht aus den Elementen $\pi\colon\mathcal{Y}\to\mathfrak{X}$ von $S_{\mathfrak{X},D}$, für die $\lambda\colon\mathcal{Y}\to\Spec\,R$ eine semistabile Kurve ist, deren generische Faser $\mathcal{Y}_{Q}$ eine glatte und projektive Kurve über $Q$ ist.\par
Eine Kurve $\lambda\colon\mathcal{Y}\to\Spec\,R$ heißt dabei nach \cite[§ 10.3]{L1} genau dann semistabil, wenn $\lambda$ flach ist und für alle $s\in\Spec\,R$ die geometrische Faser $\mathcal{Y}_{\overline{s}}$ reduziert ist und lediglich gewöhnliche Doppelpunkte als Singularitäten besitzt.\\ 
\end{defn}\par

Analog zu dem Fall der Unterkategorie $S_{\mathfrak{X},D}^{good}$ gilt auch für den Fall der Unterkategorie $S_{\mathfrak{X},D}^{ss}$:\\\par

\begin{thm}[\protect{\cite[Theorem 1, 3)-5)]{DW1}}]\label{thmss}
Es sei $R$ ein diskreter Bewertungsring. Dann gilt:\par
\begin{enumerate}
\item[(i)] Für jeden diskreten Bewertungsring $R'$, der über $R$ liegt, setze man \linebreak \mbox{$\mathfrak{X'}:=\mathfrak{X}\otimes_{R}R'$} und definiere $D'$ als das Urbild von $D$ in $X'$. Dann bildet der kanonische Basiswechselfunktor $S_{\mathfrak{X},D}\to S_{\mathfrak{X'},D'}$ die Unterkategorie $S_{\mathfrak{X},D}^{ss}$ auf $S_{\mathfrak{X'},D'}^{ss}$ ab.
\item[(ii)] Für je endlich viele Objekte $\pi_{i}\colon\mathcal{Y}_{i}\to\mathfrak{X}$ aus $S_{\mathfrak{X},D}$ existiert eine endliche Erweiterung $Q'$ von $Q$, so dass es ein Objekt $\pi'\in S_{\mathfrak{X'},D'}^{ss}$ gibt, das alle Objekte $\pi_{i}\otimes_{R}R'\in S_{\mathfrak{X'},D'}$  dominiert. Dabei ist $R'$ ein diskreter Bewertungsring in $Q'$, der $R$ dominiert, $\mathfrak{X'}:=\mathfrak{X}\otimes_{R}R'$ und $D'$ das Urbild von $D$ in $X'$.
\item[(iii)] Zu jedem Objekt $\pi\colon\mathcal{Y}\to\mathfrak{X}$ aus $S_{\mathfrak{X},D}$ existiert eine Erweiterung von diskreten Bewertungsringen $R'/R$ wie in $(ii)$, so dass $\pi\otimes_{R}R'$ von 
einem Objekt aus $S_{\mathfrak{X'},D'}^{ss}$ strikt dominiert wird.
\end{enumerate}
\end{thm}
\bewub{\begin{enumerate} \item[(i)] Dies ist klar, da Glattheit und Flachheit sowie Semistabilität (siehe \cite[Proposition 10.3.15(a)]{L1}) stabil unter Basiswechsel sind.
\item[(ii)] Weil endliche Produkte in der Kategorie $S_{\mathfrak{X},D}$ existieren, folgt Aussage $(ii)$ aus Aussage $(iii)$.
\item[(iii)] Es sei $\mathcal{Y}^{V}\to\mathfrak{X'}$ wie im Beweis von Theorem \ref{Kgood}$(iii)$ konstruiert. Weil $\mathcal{Y}^{V}$, wie oben gesehen, irreduzibel und regulär, eigentlich und flach über $R'$ ist, ist nach einem Resultat von Lichtenbaum (\cite[Theorem 2.8]{Li}) $\mathcal{Y}^{V}$ projektiv über $R'$. Nach \cite[Theorem 0.2]{L2} existieren also eine endliche Erweiterung $Q^{\dagger}$ von $Q'$, ein diskreter Bewertungsring $R^{\dagger}$ in $Q^{\dagger}$ über $R'$ und ein semistabiles Modell $\mathcal{Y}^{\dagger}$ von $\mathcal{Y}^{V}\otimes_{R'}Q^{\dagger}$ zusammen mit einem Morphismus $\varphi\colon\mathcal{Y}^{\dagger}\to\mathcal{Y}^{V}\otimes_{R'}R^{\dagger}$ über $\Spec\,R^{\dagger}$. Nach Konstruktion ist dann die Komposition \[\mathcal{Y}^{\dagger}\xrightarrow{\varphi}\mathcal{Y}^{V}\otimes_{R'}R^{\dagger}\to\mathfrak{X}^{\dagger}=\mathfrak{X}\otimes_{R}R^{\dagger}\]
ein Element der Kategorie $S^{ss}_{\mathfrak{X}^{\dagger},D^{\dagger}}$, das den Morphismus $\pi^{\dagger}=\pi\otimes_{R}R^{\dagger}$ strikt dominiert.\hfill$\Box$\end{enumerate}}\par
\text{}\par

Für den konkreten Fall $R=\overline{\mathbb{Z}}_{p}$ und $Q=\overline{\mathbb{Q}}_{p}$ gilt wie im Falle von $S_{\mathfrak{X},D}^{good}$:\\\par

\begin{kor}[\protect{\cite[Corollary 3,1)]{DW1}}]\label{kss}
Es sei $X$ eine glatte und projektive Kurve über $\overline{\mathbb{Q}}_{p}$ und $D$ ein Divisor auf $X$ sowie $\mathfrak{X}$ ein Modell von $X$ über $\Spec\,\overline{\mathbb{Z}}_{p}$. Außerdem sei eine endliche Anzahl von Objekten $\pi_{i}\colon\mathcal{Y}_{i}\to\mathfrak{X}$ aus $S_{\mathfrak{X},D}$ gegeben.\par
Dann existieren eine endliche Erweiterung $K$ von $\mathbb{Q}_{p}$ und eine glatte und projektive Kurve $X_{K}$ über $K$ mit Modell $\mathfrak{X}_{\mathfrak{o}_{K}}$ über $\mathfrak{o}_{K}$ sowie ein Divisor $D_{K}$ von $X_{K}$, so dass folgendes gilt: Es ist $X=X_{K}\otimes_{K}\overline{K}$ und $D=D_{K}\otimes_{K}\overline{K}$ sowie $\mathfrak{X}=\mathfrak{X}_{\mathfrak{o}_{K}}\otimes_{\mathfrak{o}_{K}}\overline{\mathbb{Z}}_{p}$. Außerdem existiert ein Element $\pi_{\mathfrak{o}_{K}}\colon\mathcal{Y}_{\mathfrak{o}_{K}}\to\mathfrak{X}_{\mathfrak{o}_{K}}$ von $S_{\mathfrak{X}_{\mathfrak{o}_{K}},D_{K}}^{ss}$, so dass der Morphismus $\pi=\pi_{\mathfrak{o}_{K}}\otimes_{\mathfrak{o}_{K}}\overline{\mathbb{Z}}_{p}$ alle Morphismen $\pi_{i}$ dominiert.
\end{kor}
\bewub{Der Beweis ist wortwörtlich derselbe wie der für Korollar \ref{kgood}. Man ersetzt lediglich den Verweis auf Theorem \ref{Kgood} durch den auf Theorem \ref{thmss} und stets "`good"' durch "`ss"'. \hfill$\Box$}\\\par

Daraus folgt:\\\par

\begin{kor}[\protect{\cite[Corollary 3, 3)]{DW1}}]\label{doms}
Es sei $X$ eine glatte und projektive Kurve über $\overline{\mathbb{Q}}_{p}$ und $D$ ein Divisor auf $X$ sowie $\mathfrak{X}$ ein Modell von $X$ über $\Spec\,\overline{\mathbb{Z}}_{p}$. Dann gilt:\par
Jede endliche Anzahl von Objekten $\pi_{i}\colon\mathcal{Y}_{i}\to\mathfrak{X}$ in $S_{\mathfrak{X},D}$ wird von einem gemeinsamen Objekt $(\pi\colon\mathcal{Y}\to\mathfrak{X})\in S_{\mathfrak{X},D}^{ss}$ dominiert.
\end{kor}
\bewub{Dies folgt direkt aus Korollar \ref{kss} zusammen mit der Tatsache, dass die Definition der Kategorie $S_{\mathfrak{X},D}^{ss}$ nach Theorem \ref{thmss}$(i)$ veträglich mit Erweiterungen des zugrundeliegenden Bewertungsrings ist. \hfill$\Box$}\\\par

Neben den Kategorien $S_{\mathfrak{X},D}, S_{\mathfrak{X},D}^{ss}$ und $S_{\mathfrak{X},D}^{good}$ benötigen wir noch die folgenden Kategorien:\\\par

\begin{defn}[vgl. \protect{\cite[S.10f]{DW1}}]
\begin{itemize}
\item[i)] Die Objekte der Kategorie $\mathfrak{T}_{\mathfrak{X},D}$ sind endlich präsentierte eigentliche $G$"=äquivariante $\overline{\mathbb{Z}}_{p}$"=Morphismen $\pi\colon\mathcal{Y}\to\mathfrak{X}$, wobei $G$ eine endliche abstrakte Gruppe ist, die $\overline{\mathbb{Z}}_{p}$"=linear von links auf $\mathcal{Y}$ und trivial auf $\mathfrak{X}$ operiert. Außerdem sei verlangt, dass die generischen Fasern endlich sind und die Einschränkung $(\mathcal{Y}_{\overline{\mathbb{Q}}_{p}}-\pi^{*}D)\to (X-D)$ ein étaler $G$"=Torseur ist.\par
Ein Morphismus von dem $G_{1}$"=äquivarianten Morphismus $\pi_{1}\colon\mathcal{Y}_{1}\to\mathfrak{X}$ zu dem $G_{2}$"=äquivarianten Morphismus $\pi_{2}\colon\mathcal{Y}_{2}\to\mathfrak{X}$ besteht aus einem Morphismus $\varphi\colon\mathcal{Y}_{1}\to\mathcal{Y}_{2}$, so dass $\pi_{1}=\pi_{2}\circ\varphi$ gilt, zusammen mit einem Gruppenhomomorphismus $\gamma\colon G_{1}\to G_{2}$, so dass $\varphi$ ein $G_{1}$"=äquivarianter Morphismus ist, wenn $G_{1}$ vermöge $\gamma$ auf $\mathcal{Y}_{2}$ operiert.
\item[ii)] Offensichtlicherweise existiert ein Vergissfunktor $\mathfrak{T}_{\mathfrak{X},D}\to S_{\mathfrak{X},D}$. Die volle Unterkategorie $\mathfrak{T}_{\mathfrak{X},D}^{good}$ von $\mathfrak{T}_{\mathfrak{X},D}$ besteht aus den Elementen von $\mathfrak{T}_{\mathfrak{X},D}$, die auf Objekte in $S_{\mathfrak{X},D}^{good}$ abgebildet werden.\\
\end{itemize}
\end{defn}\par

Ein wichtiges Resultat aus \cite{DW1} ist, dass jedes Objekt von $S_{\mathfrak{X},D}$ durch das Bild eines Objektes von $\mathfrak{T}_{\mathfrak{X},D}^{good}$ dominiert wird:\\\par

\begin{thm}[\protect{\cite[Theorem 4]{DW1}}]\label{tgood}
Für jedes Objekt $\pi\colon\mathcal{Y}\to\mathfrak{X}$ aus $S_{\mathfrak{X},D}$ existieren eine endliche Gruppe $G$ und ein $G$"=äquivarianter Morphismus $\pi'\colon\mathcal{Y'}\to\mathfrak{X}$, die ein Objekt von $\mathfrak{T}_{\mathfrak{X},D}^{good}$ definieren, für das ein Morphismus $\varphi\colon\mathcal{Y'}\to\mathcal{Y}$ existiert, so dass $\pi\circ\varphi=\pi'$ ist. 
\end{thm}
\bewub{Für den genauen, etwas längeren Beweis siehe \cite[Theorem 4]{DW1}. Er benutzt wieder noetherschen Descent, die Kontruktionen der Beweise der Theoreme \ref{Kgood} und \ref{thmss} sowie verschiedene Aussagen über die $G$"=Äquivarianz eines Morphismus aus \cite{L2} und \cite{EGA II}. \hfill$\Box$}\\\par

\section{Das étale Fundamentalgruppoid und Faserfunktoren}

Im weiteren Verlauf benötigen wir den Begriff des étalen Fundamentalgruppoids und den von Faserfunktoren. Als Literatur sei hierzu auf \cite{SGA 1} oder auch \cite[S.577 ff]{DW1} verwiesen.\par
Es sei $Z$ eine Varietät über $\overline{\mathbb{Q}}_{p}$. Dann definiert man das étale Fundamentalgruppoid $\Pi_{1}(Z)$ von $Z$ wie folgt:\\\par

\begin{defn}
Es sei $Z$ eine Varietät über $\overline{\mathbb{Q}}_{p}$ und $z\in Z(\mathbb{C}_{p})$ ein geometrischer Punkt.\par 
Dann ist der zugehörige Faserfunktor $F_{z}$ definiert als der Funktor \[F_{z}=Mor_{Z}(z,\_)\] von der Kategorie der endlichen étalen Überlagerungen $Z'$ von $Z$ in die Kategorie der endlichen Mengen, der jedem solchen $Z'$ die Menge der über $z$ liegenden $\mathbb{C}_{p}$"=wertigen Punkte von $Z'$ zuordnet.\par
Außerdem sei die topologische Kategorie $\Pi_{1}(Z)$ wie folgt definiert:\par
Die Objekte der Kategorie sind die geometrischen Punkte von $Z$, d.h. die Menge der Objekte ist $Z(\mathbb{C}_{p})$.
Ein Morphismus in $\Pi_{1}(Z)$ zwischen zwei $\mathbb{C}_{p}$"=wertigen Punkten $z$ und $z^{*}$ ist ein Isomorphismus der assoziierten Faserfunktoren. Solch ein Isomorphismus der Faserfunktoren wird als étaler Weg von $z$ nach $z^{*}$ bezeichnet.\par
Da nach \cite[Exposé V, 4. bzw. 7.]{SGA 1} jeder Faserfunktor pro"=repräsentier\-bar ist, ist $Mor_{\Pi_{1}(Z)}(z,z^{*})$ eine proendliche Menge und damit ein kompakter total unzusammenhängender Hausdorffraum (vgl. \cite[S. 577]{DW1}). Ferner induziert die Komposition von Morphismen eine stetige Abbildung \[Mor_{\Pi_{1}(Z)}(z,z^{*})\times Mor_{\Pi_{1}(Z)}(z^{*},z^{**})\to Mor_{\Pi_{1}(Z)}(z,z^{**}).\]
Die so definierte topologische Kategorie $\Pi_{1}(Z)$ wird als étales Fundamentalgruppoid von $Z$ bezeichnet.\\
\end{defn}\par

Desweiteren sei an einige Funktorialitätseigenschaften von $\Pi_{1}$ erinnert (vgl. dazu \cite[S. 578f]{DW1}, die wir hier weitestgehend übernehmen):\par
Ist $\alpha\colon Z_{1}\to Z_{2}$ ein Morphismus von Varietäten über $\overline{\mathbb{Q}_{p}}$, so wird wie folgt ein stetiger Funktor $\alpha_{*}\colon\Pi_{1}(Z_{1})\to\Pi_{1}(Z_{2})$ induziert:\par Auf Ebene der Objekte ist $\alpha_{*}$ die Abbildung $\alpha\colon Z_{1}(\mathbb{C}_{p})\to Z_{2}(\mathbb{C}_{p})$. Für zwei Punkte $z,z'\in Z_{1}(\mathbb{C}_{p})$ definiert man außerdem eine stetige Abbldung $\alpha_{*}\colon\Iso(F_{z},F_{z'})\to\Iso(F_{\alpha(z)},F_{\alpha(z')})$ folgendermaßen:\par  
Für einen endlichen étalen Morphismus $Y_{2}\to Z_{2}$ betrachte man den Basiswechsel $Y_{1}=Y_{2}\times_{Z_{2}}Z_{1}\to Z_{1}$. Dann hat man kanonische Bijektionen $F_{z}(Y_{1})\cong F_{\alpha(z)}(Y_{2})$ und $F_{z'}(Y_{1})\cong F_{\alpha(z')}(Y_{2})$. Dies erlaubt es, für einen étalen Weg \mbox{$\gamma\in\Iso(F_{z},F_{z'})$} den Morphimus $\alpha_{*}(\gamma)(Y_{2})$ als die Verkettung \[\alpha_{*}(\gamma)(Y_{2})\colon F_{\alpha(z)}(Y_{2})\cong F_{z} (Y_{1})\stackrel{\gamma(Y_{1})}{\xrightarrow{\sim}} F_{z'}(Y_{1})\cong F_{\alpha(z')}(Y_{2})\] zu definieren. Dies definiert einen Isomorphismus $\alpha_{*}(\gamma)\colon F_{\alpha(z)}\xrightarrow{\sim}F_{\alpha(z')}$ von Faserfunktoren. Nach Konstruktion ist es klar, dass die Abbildung $\gamma\mapsto\alpha_{*}(\gamma)$ stetig ist und insgesamt $\alpha_{*}$ mit dieser Definition ein Funktor ist.\par
Ist $\beta\colon Z_{2}\to Z_{3}$ ein weiterer Morphismus von Varietäten über $\overline{\mathbb{Q}_{p}}$, so erhält man $(\beta\circ\alpha)_{*}=\beta_{*}\circ\alpha_{*}\colon\Pi_{1}(Z_{1})\to\Pi_{1}(Z_{3})$ und es ist offensichtlicherweise $\id_{*}=\id$.\\\par
Als nächstes werde das Verhalten von étalen Fundamentalgruppoiden unter Galoiskonjugation untersucht:\par
Für ein Schema $Y$ über $\overline{\mathbb{Q}_{p}}$ und einen $\mathbb{Q}_{p}$"=Automorphismus $\sigma$ von $\overline{\mathbb{Q}_{p}}$ setze man ${}^{\sigma}Y=Y\otimes_{\overline{\mathbb{Q}_{p}},\sigma}\overline{\mathbb{Q}_{p}}$ und schreibe $\sigma\colon Y\to{}^{\sigma}Y$ für das Inverse der Pro\-jek\-tions\-ab\-bil\-dung. Dann lässt sich wie folgt ein stetiger Funktor \[\sigma_{*}\colon\Pi_{1}(Z)\to\Pi_{1}({}^{\sigma}Z)\] definieren:\par 
Auf Ebene der Objekte sei $\sigma_{*}$ als die Abbildung von $z\in Z(\mathbb{C}_{p})$ auf \linebreak \mbox{${}^{\sigma}z=\sigma\circ z\circ\sigma^{-1}$} in ${}^{\sigma}Z(\mathbb{C}_{p})$ definiert. Außerdem definiert man die stetige Abbildung \[\sigma_{*}\colon\Iso(F_{z},F_{z'})\xrightarrow{\sim}\Iso(F_{{}^{\sigma}z},F_{{}^{\sigma}z'})\] zwischen den Räumen von Morphismen wie folgt:\par
Es ist klar, dass jede endliche étale Überlagerung von ${}^{\sigma}Z$ von der Gestalt ${}^{\sigma}Y$  für eine endliche étale Überlagerung $Y$ von $Z$ ist und dass in kanonischer Weise $F_{{}^{\sigma}z}({}^{\sigma}Y)\cong F_{z}(Y)$ für jeden Punkt $z\in Z(\mathbb{C}_{p})$ gilt. Damit kann man $\sigma_{*}(\gamma)({}^{\sigma}Y)$ als die Komposition \[\sigma_{*}(\gamma)({}^{\sigma}Y)\colon F_{{}^{\sigma}z}({}^{\sigma}Y)\cong F_{z}(Y)\xrightarrow{\gamma}F_{z'}(Y)\cong F_{{}^{\sigma}z'}({}^{\sigma}Y)\] definieren, wodurch ein Isomorphismus von Faserfunktoren $\sigma_{*}(\gamma)$ definiert wird. Die Abbildung  $\gamma\mapsto\sigma_{*}(\gamma)$ ist wieder stetig und man erhält durch diese Definitionen einen wohldefinierten Funktor $\sigma_{*}$. Es ist zudem offensichtlich, dass $(\sigma\tau)_{*}=\sigma_{*}\circ\tau_{*}$ als Funktoren von $\Pi_{1}(Z)$ nach $\Pi_{1}({}^{\sigma\tau}Z)=\Pi_{1}({}^{\sigma}({}^{\tau} Z))$ gilt.\\\par

Ist ferner $Z$ bereits über einer endlichen Erweiterung $\mathbb{Q}_{p}\subset K\subset\overline{\mathbb{Q}_{p}}$ von $\mathbb{Q}_{p}$ definiert, d.\,h. ist also $Z=Z_{K}\otimes_{K}\overline{\mathbb{Q}_{p}}$ für eine Varietät $Z_{K}$ über $K$, so liefert für jedes Element $\sigma$ der Galoisgruppe $G_{K}=\Gal(\overline{\mathbb{Q}_{p}}/K)$ die Abbildung $\id\times_{\Spec\,K}\Spec(\sigma^{-1})$ einen $\overline{\mathbb{Q}_{p}}$"=linearen Isomorphismus ${}^{\sigma}Z\xrightarrow{\sim}Z$, mittels dessen man ${}^{\sigma}Z$ mit $Z$ identifizieren kann. Es folgt, dass für ein derartiges $Z$ die Gruppe $G_{K}$ von links durch stetige Automorphismen auf der Kategorie $\Pi_{1}(Z)$ operiert.\\\par
\newpage
\thispagestyle{empty}

\chapter{Paralleltransport für Prin\-zi\-pal\-bün\-del}
\section{Konstruktion des Funktors $\mathscr{B}_{\mathfrak{X}_{\mathfrak{o}},D}(G)\to \Rep_{\pi_{1}(U)}(G(\mathfrak{o}))$}

Es sei $\mathfrak{o}$ der Ring der ganzen Zahlen in $\mathbb{C}_{p}$ und $S:=\Spec\,\mathfrak{o}$. Außerdem sei $G:=\Spec\,A$ ein glattes und affines Gruppenschema über $\mathfrak{o}$ von endlicher Präsentation. Es folgt sofort, dass $G$ auch treuflach über $\mathfrak{o}$ ist, da die Glattheit von $G$ impliziert, dass $G$ flach über $\mathfrak{o}$ ist, und jedes flache Gruppenschema von endlicher Präsentation wiederum nach \cite[Exposé $VI_{B}$, Proposition 9.2. (xi),(xii)]{SGA3} treuflach ist.\par
Unter einem $G$"=Torseur (oder auch \emph{$G$"=Prinzipalbündel}) $P$ auf einem $\mathfrak{o}$"=Schema $\xi$ werde im folgenden stets ein darstellbarer Rechts"=$G_{\xi}=G\times_{\Spec\mathfrak{o}}\xi$"=Torseur für die fppf"=Topologie auf $\xi$ verstanden.\par
Ein solcher $G$"=Torseur $P$ auf einem $\mathfrak{o}$"=Schema $\xi$ ist stets glatt über $\xi$:\\\par

\begin{lem}\label{descent}
Es sei $P$ ein $G$"=Torseur auf einem $\mathfrak{o}$"=Schema $\xi$. Dann ist $P$ glatt über $\xi$.
\end{lem}
\bewub{Da nach Voraussetzung $G$ glatt über $\mathfrak{o}$ ist, ist $G_{\xi}$ glatt über $\xi$ und damit auch $P\times_{\xi}G_{\xi}$ glatt über $P$. Da nach Definition eines Torseurs $P\times_{\xi}G_{\xi}$ isomorph zu $P\times_{\xi}P$ ist, ist dann auch $P\times_{\xi}P$ glatt über $P$. Weil $P$ über $\xi$ nach Definition treuflach und quasi"=kompakt ist, folgt mittels Descenttheorie aus dem kartesischen Diagramm \begin{equation*}
        \xymatrix@=3em{%
        P\times_{\xi}P \ar[r]\ar[d]& P \ar[d]\\
        P \ar[r]& \xi,}
           \end{equation*} dass $P$ glatt über $\xi$ ist.$\hfill\Box$}\\\par
           
Ferner halten wir fest:\\\par

\begin{sat}\label{etale}
Es seien $G$, $P$, $\xi$ wie bisher. Dann ist $P$ sogar ein $G$"=Torseur für die étale Topologie auf $\xi$. 
\end{sat}
\bewub{Nach Satz \ref{koh} können wir $P$ mit einer Klasse $c(P)\in H^{1}_{fppf}(\xi,G)$ identifizieren. Nach \cite[Théorème 11.7.]{Gr} ist unter den gegebenen Voraussetzungen der kanonische Morphismus $H^{1}_{\acute{e}t}(\xi,G_{\acute{e}t})\to H^{1}_{fppf}(\xi,G)$, der durch den kanonischen Morphismus von Siten $j\colon\xi_{fppf}\to\xi_{\acute{e}t}$ induziert wird, ein Isomorphismus, wobei $G_{\acute{e}t}$ die durch $G$ induzierte étale Garbe ist. Also können wir $c(P)$ in eindeutiger Weise mit einer Klasse $c'\in H^{1}_{\acute{e}t}(\xi,G)$ identifizieren. Weil $H^{1}_{\acute{e}t}(X,G)=\check{H}^{1}_{\acute{e}t}(X,G)$ gilt, ist diese Klasse $c'$ nach Satz \ref{koh} durch einen Garben"=$G$"=Torseur $P'$ für die étale Topologie auf $\xi$ gegeben, d.\,h. es ist $P=j^{*}P'$.  
Aufgrund der Isomorphismen und der Darstellbarkeit von $P$ folgt also, dass $P'$ durch $P$ darstellbar ist und $P$ einen $G$"=Torseur für die étale Topologie definiert. 
\hfill$\Box$
}\\\par

Für einen $G$"=Torseur $P$ auf $\xi$ definiert man $\Gamma(\xi, P)$ als die Menge der Schnitte von $P$. Fasst man $P$ als Garbe auf, so liefert dies gerade die gewöhnliche Definition von $\Gamma(\xi, P)$. Da jedes Element $(\alpha\colon\mathfrak{o}\to G)\in G(\mathfrak{o})$ ein Element $(\alpha_{\xi}\colon\xi\to G_{\xi})\in G_{\xi}(\xi)$ induziert und $G_{\xi}$ als Garbe auf $P$ als Garbe operiert, operiert $G_{\xi}(\xi)$ auf $\Gamma(\xi, P)$, so dass $\Gamma(\xi, P)$ ein Rechts"=$G(\mathfrak{o})$"=Modul ist.\\\par

In Analogie zu dem in \cite{DW1} dargestellten Fall von Vektorbündeln soll im folgenden ein étaler Paralleltransport für Prinzipalbündel definiert werden. Dazu wird zunächst wie in \cite{DW1} ein Funktor $\rho$ von einer gewissen Kategorie $\mathscr{B}_{\mathfrak{X}_{\mathfrak{o}},D}(G)$ von Prinzipalbündeln in die Kategorie $\Rep_{\Pi_{1}(X-D)}(G(\mathfrak{o}))$ der stetigen Funktoren vom étalen Fundamentalgruppoid $\pi_{1}(X-D)$ in die Kategorie $\mathcal{P}(G(\mathfrak{o}))$ definiert. Diese Kategorie $\mathcal{P}(G(\mathfrak{o}))$ ist wiederum die Kategorie der topologischen Räume mit einer einfach"=transitiven und stetigen Rechts"=$G(\mathfrak{o})$"=Operation, wobei die Morphismen in dieser Kategorie gerade deren stetige $G(\mathfrak{o})$"=äquivariante Abbildungen sind. Die Definition dieses Funktors $\rho\colon\mathscr{B}_{\mathfrak{X}_{\mathfrak{o}},D}(G)\to\Rep_{\Pi_{1}(X-D)}G(\mathfrak{o})$ wurde bereits in \cite{DW3} skizziert.\\\par

Für einen Morphismus $f\colon\xi_{1}\to\xi_{2}$ von $\mathfrak{o}$"=Schemata und einen $G$"=Torseur $P$ auf $\xi_{2}$ ist nach Bemerkung \ref{affin1} und Lemma \ref{pullback} $f^{*}P=\xi_{1}\times_{\xi_{2}}P$ ein $G_{\xi_{1}}$"=Torseur auf $\xi_{1}$. Ebenso folgt, wenn man die Torseure jeweils als Garben betrachtet, aus Lemma \ref{pullback} und dem Adjunktionsmorphismus $P\to f_{*}f^{*}P$ sowie der Definition von $f_{*}$, dass es eine natürliche Rechts"=$G(\mathfrak{o})$"=äquivariante Abbildung $f^{*}\colon\Gamma(\xi_{2},P)\to\Gamma(\xi_{1},f^{*}P)$ gibt.\\\par

Als nächstes führen wir einige Bezeichnungen ein:\\\par
Für $n\geq 1$ sei $\mathfrak{o}_{n}:=\mathfrak{o}/p^{n}\mathfrak{o}=\overline{\mathbb{Z}}_{p}/p^{n}\overline{\mathbb{Z}}_{p}$. Dann ist $\mathfrak{o}=\projlim_{n}\mathfrak{o}_{n}$. Ebenso setze man $A_{n}:=A\otimes_{\mathfrak{o}}\mathfrak{o}_{n}$. Dann ist $G_{n}:=\Spec\,A_{n}$ die Reduktion von $G$ modulo $p^{n}$. Für jedes $\mathfrak{o}$"=Schema $\xi$ setze man ferner $\xi_{n}:=\xi\otimes_{\mathfrak{o}}\mathfrak{o}_{n}$. Für einen $G$"=Torseur $P$ auf $\xi$ schließlich sei $P_{n}$ die Reduktion von $P$ zu einem $G$"= bzw. $G_{n}$"=Torseur auf $\xi_{n}$, d.\,h. es sei $P_{n}=i_{n}^{*}P$, wobei $i_{n}\colon\xi_{n}\to\xi$ der kanonische Morphismus ist. \\\par

Wir versehen nun $G(\mathfrak{o})$ mit einer prodiskreten Topologie, so dass $G(\mathfrak{o})$ zu einer topologischen Gruppe wird:\par Es gibt einen kanonischen Isomorphismus \[G(\mathfrak{o})\cong Hom_{\mathfrak{o}"=Alg.}(A,\mathfrak{o})\cong\projlim_{n}Hom_{\mathfrak{o}"=Alg.}(A,\mathfrak{o}_{n})\cong\projlim_{n}G(\mathfrak{o}_{n}).\] Versieht man $G(\mathfrak{o}_{n})$ mit der diskreten Topologie, so erhält man auf $\prod_{n\geq 1}G(\mathfrak{o}_{n})$ die induzierte Produkttopologie und auf $\projlim_{n}G(\mathfrak{o}_{n})$ die induzierte Teilraumtopologie. Mittels des obigen Isomorphismus liefert das eine Topologie auf $G(\mathfrak{o})$.\par
Da $G$ nach Voraussetzung glatt und damit insbesondere formell glatt über $\mathfrak{o}$ ist, ist ferner die Abbildung \[Hom_{\Spec\,\mathfrak{o}}(\Spec\,\mathfrak{o}_{n+1},G)\to Hom_{\Spec\,\mathfrak{o}}(\Spec\,\mathfrak{o}_{n},G)\] surjektiv, so dass alle Übergangsabbildungen $G(\mathfrak{o}_{n+1})\to G(\mathfrak{o}_{n})$ surjektiv sind. Nach Konstruktion des projektiven Limes sind damit alle Morphismen \linebreak\mbox{$G(\mathfrak{o})\to G(\mathfrak{o}_{n})$} ebenfalls surjektiv.\\\par

Mit diesen Vorbereitungen können wir nun eine Kategorie $\mathscr{B}_{\mathfrak{X}_{\mathfrak{o}},D}(G)$ definieren:\\\par

\begin{defn}
Für ein Modell $\mathfrak{X}$ über $\overline{\mathbb{Z}}_{p}$ einer glatten und projektiven Kurve $X$ über $\overline{\mathbb{Q}}_{p}$ und einen Divisor $D$ auf $X$ definiere man $\mathscr{B}_{\mathfrak{X}_{\mathfrak{o}},D}(G)$ als die Kategorie der $G$"=Torseure $P$ auf $\mathfrak{X}_{\mathfrak{o}}=\mathfrak{X}\otimes_{\overline{\mathbb{Z}}_{p}}\mathfrak{o}$ mit der folgenden Eigenschaft:\par
Für alle $n\geq 1$ existiert ein Objekt $\pi\colon\mathcal{Y}\to\mathfrak{X}$ von $S_{\mathfrak{X},D}$, so dass $\pi_{n}^{*}P_{n}$ trivial ist, d.\,h. als $G_{n}$"=Torseur isomorph zu dem trivialen Torseur $\mathcal{Y}_{n}\otimes_{\mathfrak{o}_{n}}G_{n}$ ist.\par
Dabei sei $\pi_{n}\colon\mathcal{Y}_{n}\to\mathfrak{X}_{n}$ der von $\pi$ induzierte Morphismus.
\end{defn}\par

Im folgenden seien nun stets $X$ eine glatte projektive Kurve über $\overline{\mathbb{Q}}_{p}$, $\mathfrak{X}$ ein Modell von $X$ über $\overline{\mathbb{Z}_{p}}$ und $D$ ein Divisor auf $X$.\\\par

\begin{bem}\label{pullback3}
Für einen Morphismus $f\colon\mathfrak{X}\to\mathfrak{X'}$ von Modellen über $\overline{\mathbb{Z}_{p}}$ der glatten und projektiven Kurven $X$ und $X'$ über $\overline{\mathbb{Q}_{p}}$ und einen Divisor $D'$ auf $X'$ induziert der Pullback"=Funktor $f^{*}$ für $G$"=Torseure einen Funktor \[f^{*}\colon\mathscr{B}_{\mathfrak{X'}_{\mathfrak{o}},D'}(G)\to\mathscr{B}_{\mathfrak{X}_{\mathfrak{o}},f^{*}D'}(G).\]
\end{bem}
\bewub{Es sei $P'\in\mathscr{B}_{\mathfrak{X'}_{\mathfrak{o}},D'}(G)$ und $n\in\nat, n\geq 1$ beliebig. Dann existiert nach Definition von $\mathscr{B}_{\mathfrak{X'}_{\mathfrak{o}},D'}(G)$ ein Objekt $\pi'\colon\mathcal{Y'}\to\mathfrak{X'}$ von $S_{\mathfrak{X'},D'}$, so dass $\pi'^{*}_{n}P'_{n}$ trivial ist.\par
Setzt man $\mathcal{Y}:=\mathcal{Y'}\times_{\mathfrak{X'}}\mathfrak{X}$ und bezeichnet $\pi$ den von $\pi'$ induzierten Morphismus $\pi\colon\mathcal{Y}\to\mathfrak{X}$, so gilt für den $G$"=Torseur $P:=f^{*}P'$ auf $X$, dass $\pi^{*}_{n}P_{n}$ trivial ist. Da die Definition der Kategorie $S_{\mathfrak{X},D}$ nach dem ersten Kapitel der Arbeit verträglich mit Basiswechsel ist, ist $\pi\in S_{X,D}$, wobei $D$ das Urbild von $D'$ auf $X$ ist, so dass $P\in\mathscr{B}_{\mathfrak{X}_{\mathfrak{o}},f^{*}D'}(G)$ gilt. \hfill$\Box$}\\\par

Es sei nun $x\colon\Spec\,\mathbb{C}_{p}\to X$ ein $\mathbb{C}_{p}$"=wertiger Punkt von $X$ und $P$ ein $G$"=Torseur über $\mathfrak{X}_{\mathfrak{o}}$. Dann induziert $x$ einen Morphismus $\widetilde{x}\colon\Spec\,\mathbb{C}_{p}\to X\to\mathfrak{X}$. Da $\mathfrak{X}$ eigentlich über $\overline{\mathbb{Z}}_{p}$ ist, existiert ein eindeutig bestimmter Morphismus \[x'=\colon\Spec\,\mathfrak{o}\to\mathfrak{X},\] so dass das Diagramm \begin{equation*}
        \xymatrix@=3em{%
        \Spec\,\mathbb{C}_{p} \ar[r]^{\widetilde{x}} \ar[d] & \mathfrak{X} \ar[d]\\
        \Spec\,\mathfrak{o} \ar[r] \ar[ur]^{x'} & \Spec\,\overline{\mathbb{Z}}_{p} }
           \end{equation*} kommutiert. Dieser $\mathfrak{o}$"=wertige Punkt $x'$ von $\mathfrak{X}$ induziert wiederum einen eindeutig bestimmten Morphismus \[x_{\mathfrak{o}}\colon\Spec\mathfrak{o}\to\mathfrak{X}_{\mathfrak{o}}.\] Für $n\in\nat,n\geq1$ setze man ferner \mbox{$x_{n}:=x_{\mathfrak{o}}\circ(\Spec\,\mathfrak{o}_{n}\to\Spec\,\mathfrak{o})$}, so dass sich Gruppen $P_{x_{\mathfrak{o}}}$ und $P_{x_{n}}$ durch $P_{x_{\mathfrak{o}}}:=(x_{\mathfrak{o}}^{*}P)(\mathfrak{o})$ und \mbox{$P_{x_{n}}:=(x_{n}^{*}P)(\mathfrak{o}_{n})$} definieren lassen.\par
           
 Um auf $P_{x_{\mathfrak{o}}}$ bzw. $P_{x_{n}}$ Topologien zu definieren, benötigen wir das folgende Lemma:\\\par
 
 \begin{lem}\label{affin}
 Der $G$"=Torseur $x_{\mathfrak{o}}^{*}P$ über der affinen Basis $\Spec\,\mathfrak{o}$ ist affin.
 \end{lem}
 \bewub{Da nach Voraussetzung $G$ affin über $\mathfrak{o}$ ist, ist auch $G_{\mathfrak{X}_{\mathfrak{o}}}$ affin über $\mathfrak{X}_{\mathfrak{o}}$. Also ist auch der Morphismus $P\times_{\mathfrak{X}_{\mathfrak{o}}}G_{\mathfrak{X}_{\mathfrak{o}}}\to P$ affin. Da nach Definition eines Torseurs $P\times_{\mathfrak{X}_{\mathfrak{o}}}G_{\mathfrak{X}_{\mathfrak{o}}}$ isomorph zu $P\times_{\mathfrak{X}_{\mathfrak{o}}}P$ ist, ist dann auch $P\times_{\mathfrak{X}_{\mathfrak{o}}}P$ affin über $P$. Weil $P$ nach Definition treuflach und quasi"=kompakt über $\mathfrak{X}_{\mathfrak{o}}$ ist, folgt mittels Descenttheorie aus dem kartesischen Diagramm \begin{equation*}
        \xymatrix@=3em{%
        P\times_{\mathfrak{X}_{\mathfrak{o}}}P \ar[r]\ar[d]& P \ar[d]\\
        P \ar[r]& \mathfrak{X}_{\mathfrak{o}}, }
           \end{equation*} dass $P$ affin über $\mathfrak{X}_{\mathfrak{o}}$ ist.$\hfill\Box$}\\\par
           
Nach Definition von $x_{\mathfrak{o}}$ und $x_{n}$ gilt ferner $(x_{\mathfrak{o}}^{*}P)_{n}=x_{n}^{*}P$. Daraus folgt:\\\par

\begin{lem}
Es gilt $P_{x_{\mathfrak{o}}}=\projlim_{n}P_{x_{n}}$.
\end{lem}
\bewub{Nach Definition ist $P_{x_{\mathfrak{o}}}=(x_{\mathfrak{o}}^{*}P)(\mathfrak{o})$. Da nach Lemma \ref{affin} $x_{\mathfrak{o}}^{*}P$ affin über $\mathfrak{o}$ ist, etwa $x_{\mathfrak{o}}^{*}P=\Spec\,B$, ist dies isomorph zu \[Hom_{\mathfrak{o}-Alg}(B,\mathfrak{o})=\projlim_{n}Hom_{\mathfrak{o}-Alg}(B,\mathfrak{o}_{n}).\] Da $f(p^{n}B)=0$ für jeden Morphismus von $\mathfrak{o}$"=Algebren $f\colon B\to\mathfrak{o}_{n}$ gilt, faktorisiert jeder solcher Morphismus $f$ über $B/p^{n}B$ $=B\otimes_{\mathfrak{o}}\mathfrak{o}_{n}$. Also ist \[Hom_{\mathfrak{o}-Alg}(B,\mathfrak{o}_{n})\cong Hom_{\mathfrak{0}_{n}-Alg}(B\otimes_{\mathfrak{o}}\mathfrak{o}_{n},\mathfrak{o}_{n}).\] Damit gilt unter Verwendung von $(x_{\mathfrak{o}}^{*}P)_{n}=x_{n}^{*}P$:
\begin{align*}
P_{x_{\mathfrak{o}}} & \cong\projlim_{n}Hom_{\mathfrak{0}_{n}-Alg}(B\otimes_{\mathfrak{o}}\mathfrak{o}_{n},\mathfrak{o}_{n})\\
& \cong\projlim_{n}(i_{n}^{*}(x_{\mathfrak{o}}^{*}P))(\mathfrak{o}_{n})\cong\projlim_{n}(x_{n}^{*}P)(\mathfrak{o}_{n})=\projlim_{n}P_{x_{n}}
\end{align*}$\hfill\Box$}\\\par

Versieht man nun $P_{x_{n}}$ für alle $n$ jeweils mit der diskreten Topologie, so erhält man auf $P_{x_{\mathfrak{o}}}$ eine prodiskrete Topologie. Dadurch wird $P_{x_{\mathfrak{o}}}$ zu einem topologischen Raum mit stetiger $G(\mathfrak{o})$"=Aktion. Dies liegt daran, dass, wie bereits gesehen, in gleicher Weise $G(\mathfrak{o})=\projlim_{n}G(\mathfrak{o}_{n})$ gilt, wobei $G(\mathfrak{o}_{n})$ bzw. $G(\mathfrak{o})$ mit der diskreten bzw. der induzierten prodiskreten Topologie versehen wurden, sowie an der Tatsache, dass die Operationen von $G(\mathfrak{o}_{n})$ auf $P_{x_{n}}$ für die einzelnen $n$ miteinander verträglich und jeweils auch mit der Operation von $G(\mathfrak{o})$ auf $P_{x_{\mathfrak{o}}}$ verträglich sind.\par
Die Operation von $G(\mathfrak{o}_{n})$ bzw. $G(\mathfrak{o})$ auf $P_{x_{n}}$ bzw. $P_{x_{\mathfrak{o}}}$ lässt sich mit Hilfe des folgenden Lemmas genauer beschreiben:\\\par

\begin{lem}
$x_{\mathfrak{o}}^{*}P$ ist ein trivialer Torseur für die étale Topologie und damit auch für die fppf"=Topologie.
\end{lem}
\bewub{Da $\mathfrak{o}$ ein vollständiger lokaler Ring ist, ist nach \cite[Propostion I.4.5]{M} $\mathfrak{o}$ und damit $\Spec\,\mathfrak{o}$ strikt henselsch.\par
Außerdem folgt für die Wahl $\xi=\mathfrak{X}_{\mathfrak{o}}$ aus Satz \ref{etale}, dass $P$ auch ein Torseur für die étale Topologie ist. Nach Lemma \ref{pullback} ist dann auch $x_{\mathfrak{o}}^{*}P$ ein Torseur für die étale Topologie und definiert also in eindeutiger Weise nach Satz \ref{koh} ein Element $c(x_{\mathfrak{o}}^{*}P)\in H^{1}_{\acute{e}t}(\Spec\,\mathfrak{o},G)$. Nach \cite[Corollary I.1.19]{FK} ist aber $H^{1}_{\acute{e}t}(\Spec\,\mathfrak{o},G)=0$, also $c(x_{\mathfrak{o}}^{*}P)=0$ und damit $x_{\mathfrak{o}}^{*}P$ trivial als étaler Torseur, also auch trivial als Torseur für die fppf"=Topologie. $\hfill\Box$}\\\par

Damit operieren also $G(\mathfrak{o})$ bzw. $G(\mathfrak{o}_{n})$ jeweils einfach"=transitiv auf $P_{x_{\mathfrak{o}}}$ bzw. $P_{x_{n}}$.\par

Es gilt nun:\\\par

\begin{lem}\label{glatt}
 Der $G$"=Torseur $x_{\mathfrak{o}}^{*}P$ ist glatt über $\Spec\,\mathfrak{o}$.
 \end{lem}
 \bewub{Da nach Voraussetzung $G$ glatt über $\mathfrak{o}$ ist, ist auch der Morphismus $x_{\mathfrak{o}}^{*}P\times_{\Spec\,\mathfrak{o}}G\to x_{\mathfrak{o}}^{*}P$ glatt. Da nach Definition eines Torseurs $x_{\mathfrak{o}}^{*}P\times_{\Spec\,\mathfrak{o}}G$ isomorph zu $x_{\mathfrak{o}}^{*}P\times_{\Spec\,\mathfrak{o}}x_{\mathfrak{o}}^{*}P$ ist, ist dann auch $x_{\mathfrak{o}}^{*}P\times_{\Spec\,\mathfrak{o}}x_{\mathfrak{o}}^{*}P$ glatt über $x_{\mathfrak{o}}^{*}P$. Weil $x_{\mathfrak{o}}^{*}P$ über $\Spec\,\mathfrak{o}$ nach Definition treuflach und quasi"=kompakt ist, folgt mittels Descenttheorie aus dem kartesischen Diagramm \begin{equation*}
        \xymatrix@=3em{%
        x_{\mathfrak{o}}^{*}P\times_{\mathfrak{o}}x_{\mathfrak{o}}^{*}P \ar[r]\ar[d]&  x_{\mathfrak{o}}^{*}P \ar[d]\\
         x_{\mathfrak{o}}^{*}P \ar[r]& \Spec\,\mathfrak{o}, }
           \end{equation*} dass $x_{\mathfrak{o}}P$ glatt über $\mathfrak{o}$ ist.$\hfill\Box$}\\\par

Da $x_{\mathfrak{o}}^{*}P$ also glatt und damit insbesondere formell glatt über $\mathfrak{o}$ ist, ist die Abbildung \[Hom_{\Spec\,\mathfrak{o}}(\Spec\,\mathfrak{o}_{n+1},x_{\mathfrak{o}}^{*}P)\to Hom_{\Spec\,\mathfrak{o}}(\Spec\,\mathfrak{o}_{n},x_{\mathfrak{o}}^{*}P)\] surjektiv, so dass alle Übergangsabbildungen $x_{\mathfrak{o}}^{*}P(\mathfrak{o}_{n+1})\to x_{\mathfrak{o}}^{*}P (\mathfrak{o_{n}})$ surjektiv sind. Nach Konstruktion des projektiven Limes sind damit alle Morphismen $x_{\mathfrak{o}}^{*}P(\mathfrak{o})\to x_{\mathfrak{o}}^{*}P(\mathfrak{o}_{n})$ ebenfalls surjektiv.\\\par

Desweiteren seien die folgenden Kategorien definiert:\\\par

\begin{defn}
Die Kategorie $\mathcal{P}(G(\mathfrak{o}_{n}))$ sei für $n\in\nat, n\geq 1$ wie folgt definiert:\par
Die Objekte sind Mengen mit einer einfach"=transitiven Rechts"=$G(\mathfrak{o}_{n})$"=Aktion; die Morphismen sind $G(\mathfrak{o}_{n})$"=äquivariante Abbildungen.\par
Außerdem definiert man die Kategorie $\mathcal{P}(G(\mathfrak{o}))$, indem man als Objekte topologische Räume mit einer einfach"=transitiven und stetigen Rechts"=$G(\mathfrak{o})$"=Aktion und als Morphismen deren stetige $G(\mathfrak{o})$"=äquivariante Abbildungen wählt.\\
\end{defn}\par

\begin{lem}
Es seien $P_{1},P_{2}\in\mathcal{P}(G(\mathfrak{o}))$. Dann bestimmt die Wahl zweier Elemente $p_{1}\in P_{1}$ und $p_{2}\in P_{2}$ eine Bijektion \begin{align*}{\Phi\colon G(\mathfrak{o}) & \xrightarrow{\sim} Mor_{\mathcal{P}(G(\mathfrak{o}))}(P_{1},P_{2})\\
g & \longmapsto(\varphi_{g}\colon P_{1}\to P_{2},p_{1}h\mapsto p_{2}gh). }\end{align*}
\end{lem}
\bewub{Für alle \mbox{$g\in G(\mathfrak{o})$} ist $\varphi_{g}$ ein wohldefiniertes Element von \linebreak $Mor_{\mathcal{P}(G(\mathfrak{o}))}(P_{1},P_{2})$, da $G(\mathfrak{o})$ einfach"=transitiv auf $P_{1}$ und $P_{2}$ operiert.\par
Ist $\alpha\in Mor_{\mathcal{P}(G(\mathfrak{o}))}(P_{1},P_{2})$, so folgt aus der Tatsache, dass $G(\mathfrak{o})$ einfach"=transitiv auf $P_{1}$ und $P_{2}$ operiert, dass $\alpha=\varphi_{\alpha(p_{1})}$ ist. Damit ist die Abbildung $\Phi$ surjektiv.\par
Ist ferner $\varphi_{g}=\varphi_{g'}$ für zwei Elemente $g\in G(\mathfrak{o})$, so ist $p_{2}gh=p_{2}g'h$ für alle $h\in G(\mathfrak{o})$ und damit $g=g'$, so dass $\Phi$ injektiv und also insgesamt bijektiv ist.\hfill$\Box$}\\\par

\begin{bem}
Die Topologie auf $Mor_{\mathcal{P}(G(\mathfrak{o}))}(P_{1},P_{2})$, die diesen Morphismus zu einem Homöomorphismus macht, ist unabhängig von der Wahl von $p_{1}$ und $p_{2}$. Ersetzt man beispielsweise $p_{2}$ durch $p'_{2}$, so existiert ein Element $g'\in G$, so dass $p_{2}=p'_{2}g'$ ist. Ist dann $M\subset Mor_{\mathcal{P}(G(\mathfrak{o}))}(P_{1},P_{2})$ offen bezüglich der Topologie für den Fall der Wahl von $p_{1}$ und $p_{2}$, so ist $\Phi'^{-1}(M)=\Phi^{-1}(M)g'$, wenn $\Phi'$ die durch $p_{1}$ und $p'_{2}$ induzierte Bijektion $\Phi'\colon G(\mathfrak{o})\xrightarrow{\sim} Mor_{\mathcal{P}(G(\mathfrak{o}))}(P_{1}, P_{2})$ ist, und damit auch $\Phi'^{-1}(M)$ eine offene Teilmenge von $G(\mathfrak{o})$, so dass $M$ auch für den Fall der Wahl von $p_{1}$ und $p'_{2}$ eine offene Teilmenge von $Mor_{\mathcal{P}(G(\mathfrak{o}))}(P_{1}, P_{2})$ ist.\par
Auf analoge Weise sieht man, dass die auf $Mor_{\mathcal{P}(G(\mathfrak{o}))}(P_{1}, P_{2})$ induzierte Topologie unabhängig von der Wahl von $p_{1}$ ist.\par
Damit wird $\mathcal{P}(G(\mathfrak{o}))$ zu einer topologischen Kategorie, d.\,h. für alle $P_{1}$ und $P_{2}$ aus $\mathcal{P}(G(\mathfrak{o}))$ ist $Mor_{\mathcal{P}(G(\mathfrak{o}))}(P_{1},P_{2})$ ein topologischer Raum und die Komposition von Morphismen ist stetig.\\
\end{bem}\par

Es sei nun $P$ ein $G$"=Torseur in $\mathscr{B}_{\mathfrak{X}_{\mathfrak{o}},D}(G)$. Für alle $n\geq 1$ wähle man dann ein Objekt $(\pi\colon\mathcal{Y}\to\mathfrak{X})\in S_{\mathfrak{X},D}^{good}$, so dass $\pi_{n}^{*}P_{n}$ trivial auf $\mathcal{Y}_{n}$ ist. Dies geht nach folgendem Satz:\\\par

\begin{sat}\label{good}
Ist $P$ ein $G$"=Torseur in $\mathscr{B}_{\mathfrak{X}_{\mathfrak{o}},D}(G)$, so existiert für alle $n\geq 1$ ein Objekt $(\pi\colon\mathcal{Y}\to\mathfrak{X})\in S_{\mathfrak{X},D}^{good}$, so dass $\pi_{n}^{*}P_{n}$ trivial auf $\mathcal{Y}_{n}$ ist.
\end{sat}
\bewub{Es sei $P$ ein $G$"=Torseur in $\mathscr{B}_{\mathfrak{X}_{\mathfrak{o}},D}(G)$. Nach Definition von $\mathscr{B}_{\mathfrak{X}_{\mathfrak{o}},D}(G)$ existiert dann für alle $n\geq 1$ ein Objekt $(\pi\colon\mathcal{Y}\to\mathfrak{X})\in S_{\mathfrak{X},D}$, so dass $\pi_{n}^{*}P_{n}$ trivial auf $\mathcal{Y}_{n}$ ist. Nach Korollar \ref{dom} existiert ein Objekt \mbox{$(\pi'\colon\widetilde{\mathcal{Y}}\to\mathfrak{X})\in S_{\mathfrak{X},D}^{good}$}, welches $\pi$ dominiert, d.\,h. es gibt einen Morphismus $\varphi\colon\widetilde{\mathcal{Y}}\to\mathcal{Y}$, so dass $\pi\circ\varphi=\pi'$ ist. Da $\pi_{n}^{*}P_{n}$ trivial auf $\mathcal{Y}_{n}$ ist, ist dann auch $\varphi^{*}\pi_{n}^{*}P_{n}$ trivial auf $\widetilde{\mathcal{Y}}_{n}$, so dass wegen $\pi'=\varphi\circ\pi$ damit $(\pi')_{n}^{*}P_{n}$ trivial auf $\widetilde{\mathcal{Y}}_{n}$ ist.\hfill$\Box$}\\\par

Man setze nun $U=X-D$ und $V=\mathcal{Y}\otimes_{\overline{\mathbb{Z}}_{p}}\overline{\mathbb{Q}}_{p}-\pi^{*}D$. Dann ist $V$ eine endliche étale Überlagerung von $U$. Nach \cite[4.2]{Mu} existiert somit ein Punkt $y\in V(\mathbb{C}_{p})$ über $x\in U(\mathbb{C}_{p})$ und es gilt:\\\par

\begin{lem}
Die Pullback"=Abbildung \[y_{n}^{*}\colon\Gamma(\mathcal{Y}_{n},\pi_{n}^{*}P_{n})\to\Gamma(\Spec\,\mathfrak{o}_{n},y_{n}^{*}\pi_{n}^{*}P_{n})=P_{x_{n}}\] ist ein Isomorphismus von Rechts"=$G(\mathfrak{o})$"=Mengen. 
\end{lem}
\bewub{Es genügt zu zeigen, dass die Abbildung \[y_{n}^{*}\colon\Gamma(\mathcal{Y}_{n},\mathcal{Y}_{n}\times_{\Spec\,\mathfrak{o}_{n}}G_{n})\to\Gamma(\Spec\,\mathfrak{o}_{n},y_{n}^{*}(\mathcal{Y}_{n}\times_{\Spec\,\mathfrak{o}_{n}}G_{n}))=\Gamma(\Spec\,\mathfrak{o}_{n},G_{n})=G_{n}(\mathfrak{o}_{n})\] bijektiv ist.\par
Aufgrund der universellen Eigenschaft des Faserprodukts hat man die folgenden kanonischen Bijektionen:
\[\Gamma(\mathcal{Y}_{n},\mathcal{Y}_{n}\times_{\Spec\,\mathfrak{o}_{n}}G_{n})\cong Mor_{\Spec\,\mathfrak{o}_{n}}(\mathcal{Y}_{n},G_{n})\cong\Hom_{\mathfrak{o}_{n}-Alg}(A_{n},\Gamma(\mathcal{Y}_{n},\mathcal{O}_{\mathcal{Y}_{n}}))\]
Da $\pi\in S_{\mathfrak{X},D}^{good}$ ist und damit \[\Gamma(\mathcal{Y}_{n},\mathcal{O}_{\mathcal{Y}_{n}})=(\lambda_{n})_{*}\mathcal{O}_{\mathcal{Y}_{n}}(\Spec\,\mathfrak{o}_{n})=\mathcal{O}_{\Spec\,\mathfrak{o}_{n}}(\Spec\,\mathfrak{o}_{n})=\mathfrak{o}_{n}\] gilt, ist dies kanonisch isomorph zu $\Hom_{\mathfrak{o}_{n}-Alg}(A_{n},\mathfrak{o}_{n})=G_{n}(\mathfrak{o}_{n})$. Identifiziert man $\Gamma(\mathcal{Y}_{n},\mathcal{Y}_{n}\times_{\Spec\,\mathfrak{o}_{n}}G_{n})$ mit $G_{n}(\mathfrak{o}_{n})$, so wird die Abbildung $y_{n}^{*}$ zur Identität. \hfill$\Box$}\\\par

Damit kann nun für alle $P\in\mathscr{B}_{\mathfrak{X}_{\mathfrak{o}},D}(G)$ und für alle $n\in\nat, n\geq 1$ ein Funktor \[\rho_{P,n}\colon\Pi_{1}(U)\to\mathcal{P}(G(\mathfrak{o}))\] konstruiert werden, wobei wie oben $U=X-D$ sei:\\\par
Für jeden $\mathbb{C}_{p}$"=wertigen Punkt von $U$ setze man $\rho_{P,n}(x):=P_{x_{n}}$. Einem étalen Weg $\gamma$ von $x\in U(\mathbb{C}_{p})$ nach $x'\in U(\mathbb{C}_{p})$ ordnet man ferner wie folgt einen Morphismus $\rho_{P, n}(\gamma)$ in $\mathcal{P}(G(\mathfrak{o}))$ zu:\par
Man wählt zunächst ein Objekt $(\pi\colon\mathcal{Y}\to\mathfrak{X})\in S_{\mathfrak{X},D}^{good}$, so dass $\pi_{n}^{*}P_{n}$ trivial auf $\mathcal{Y}_{n}$ ist. Dies ist möglich nach Satz \ref{good}. Außerdem wähle man einen Punkt $y\in V(\mathbb{C}_{p})$ über $x\in U(\mathbb{C}_{p})$, wobei wie oben $V=\mathcal{Y}\otimes_{\overline{\mathbb{Z}}_{p}}\overline{\mathbb{Q}}_{p}-\pi^{*}_{\overline{\mathbb{Q}}_{p}}D$ sei. Ein solcher Punkt existiert nach \cite[4.2]{Mu}, da $V$ eine endliche étale Überlagerung von $X$ ist. Dann ist $\gamma y$ ein $\mathbb{C}_{p}$"=wertiger Punkt von $V$ über $x'$, wobei $\gamma y$ das Bild von $y$ unter der Abbildung $\gamma_{V}\colon F_{x}(V)\to F_{x'}(V)$ ist. Mit Hilfe dieser Wahlen lässt sich dann $\rho_{P,n}(\gamma)$ als die Komposition \[\rho_{P,n}(\gamma)\colon P_{x_{n}}\xrightarrow{(y_{n}^{*})^{-1}}\Gamma(\mathcal{Y}_{n},\pi_{n}^{*}G_{n})\xrightarrow{(\gamma y)_{n}^{*}} P_{x'_{n}}\] definieren.\\\par

Ähnlich wie in \cite[§3]{DW1} zeigt man:\\\par

\begin{lem}
Der Morphismus $\rho_{P,n}(\gamma)$ ist unabhängig von allen getroffenen Wahlen und definiert einen Morphismus in $\mathcal{P}(G(\mathfrak{o}_{n}))$.\par
Insgesamt erhält man also einen wohldefinierten stetigen Funktor $\rho_{P,n}$ von $\Pi_{1}(X-D)$ in die Kategorie $\mathcal{P}(G(\mathfrak{o}_{n}))$.
\end{lem}
\bewub{\begin{itemize}{\item[(i)] Als erstes soll gezeigt werden, dass $\rho_{P,n}(\gamma)$ unabhängig von der Wahl des Punktes $y$ über $x$ ist:\par
Dazu sei $z\in V(\mathbb{C}_{p})$ ein weiterer Punkt über $x$.\par
Nach Theorem \ref{tgood} existieren nun eine endliche Gruppe $G$ und ein $G$"=äquivarianter Morphismus $\widetilde{\pi}\colon\widetilde{\mathcal{Y}}\to\mathfrak{X}$, die ein Objekt in $\mathfrak{T}_{\mathfrak{X},D}^{good}$ definieren, zusammen mit einem Morphismus $\varphi\colon\widetilde{\mathcal{Y}}\to\mathcal{Y}$, so dass $\widetilde{\pi}=\pi\circ\varphi$ ist und also $\widetilde{\pi}$ das in der Konstruktion von $\rho_{P,n}(\gamma)$ gewählte Element \mbox{$\pi\in S_{\mathfrak{X},D}^{good}$} dominiert. Insbesondere ist  $\widetilde{V}:=\widetilde{Y}-\widetilde{\pi}^{*}D$ ein $G$"=Torseur auf $U$ und damit eine Galoisüberlagerung von $U$ mit Gruppe $G$, wobei $\widetilde{Y}$ die generische Faser von $\widetilde{\mathcal{Y}}$ bezeichnet. Wählt man nun zwei $\mathbb{C}_{p}$"=wertige Punkte $\widetilde{y}$ und $\widetilde{z}$ von $V$ über $y$ bzw. $z$ (die Punkte $\widetilde{y}$ und $\widetilde{z}$ existieren nach \cite[4.2]{Mu}), so liegen die Punkte $\gamma\widetilde{y}$ und $\gamma\widetilde{z}$ über $\gamma y$ bzw. $\gamma z$. Weil sowohl $\widetilde{y}$ als auch $\widetilde{z}$ über $x$ liegen, existiert ein eindeutig bestimmtes Element $\sigma\in G$, so dass $\sigma\widetilde{y}=\widetilde{z}$ ist und damit auch $\sigma\widetilde{y}_{\mathfrak{o}}=\widetilde{z}_{\mathfrak{o}}$ und $\sigma\widetilde{y}_{n}=\widetilde{z}_{n}$  für alle $n\geq1$ gilt.\par
Nach Konstruktion ist außerdem das folgende Diagramm kommutativ:\par
\begin{equation*}
        \xymatrix@=3em{%
        P_{x_{n}} \ar@{=}[d]& \Gamma(\mathcal{Y}_{n},\pi_{n}^{*}P_{n}) \ar[l]_{y_{n}^{*}}^{\sim} \ar[d]^{\sim}_{\varphi_{n}^{*}} \ar[r]_{\sim}^{(\gamma y)^{*}_{n}} & P_{x'_{n}} \ar@{=}[d]\\
        P_{x_{n}} & \Gamma(\widetilde{\mathcal{Y}}_{n},\widetilde{\pi}_{n}^{*}P_{n}) \ar[l]^{\sim}_{\widetilde{y}_{n}^{*}} \ar[r]_{\sim}^{(\gamma\widetilde{y})_{n}^{*}} & P_{x'_{n}}\\}
           \end{equation*}
Also gilt:
\[(\gamma y)_{n}^{*}\circ(y_{n}^{*})^{-1}=(\gamma\widetilde{y})^{*}_{n}\circ(\widetilde{y}_{n}^{*})^{-1}\]
Analog sieht man: \[(\gamma z)_{n}^{*}\circ(z_{n}^{*})^{-1}=(\gamma\widetilde{z})^{*}_{n}\circ(\widetilde{z}_{n}^{*})^{-1}\] Wegen \mbox{$\widetilde{z}=\sigma\widetilde{y}$} ist nun $\widetilde{z}_{n}^{*}$ gleich der Komposition \[\Gamma(\widetilde{\mathcal{Y}}_{n},\widetilde{\pi}_{n}^{*}P_{n})\xrightarrow{\sigma^{*}}\Gamma(\widetilde{\mathcal{Y}}_{n},\widetilde{\pi}_{n}^{*}P_{n})\xrightarrow{\widetilde{y}_{n}^{*}}P_{x_{n}}.\] Dabei ist $\sigma^{*}$ der von $\sigma$ induzierte Morphismus.\par
Aufgrund der Natürlichkeit von $\gamma$ ist desweiteren $\gamma\widetilde{z}=\sigma\gamma\widetilde{y}$ und also \[(\gamma\widetilde{z})^{*}_{n}=(\gamma\widetilde{y})^{*}_{n}\circ\sigma^{*}.\] Also folgt: 
\[(\gamma z)_{n}^{*}\circ(z_{n}^{*})^{-1}=(\gamma\widetilde{z})^{*}_{n}\circ(\widetilde{z}_{n}^{*})^{-1}=(\gamma\widetilde{y})^{*}_{n}\circ\sigma^{*}\circ(\sigma^{*})^{-1}\circ(\widetilde{y}^{*}_{n})^{-1}=(\gamma y)_{n}^{*}\circ(y_{n}^{*})^{-1}.\] Dies zeigt, dass $\rho_{P,n}$ unabhängig von der Wahl des Punktes $y$ über $x$ ist.
\item[(ii)] Als nächstes soll gezeigt werden, dass $\rho_{P,n}$ nicht von der Wahl der Überdeckung $\pi\colon\mathcal{Y}\to\mathfrak{X}$ abhängt:\par
Es sei daher $\widetilde{\pi}\colon\widetilde{\mathcal{Y}}\to\mathfrak{X}$ ein weiteres Objekt aus $S_{\mathfrak{X},D}$, so dass $\widetilde{\pi}_{n}^{*}P_{n}$ ein trivialer $G$"=Torseur ist. Nach Korollar \ref{dom} existiert dann ein Element $\pi'\in S_{\mathfrak{X},D}^{good}$, das sowohl $\pi$ als auch $\widetilde{\pi}$ dominiert. Ohne Beschränkung der Allgemeinheit darf man somit annehmen, dass $\widetilde{\pi}$ den Morphismus $\pi$ dominiert. Mit den Notationen von oben wähle man einen Punkt \mbox{$\widetilde{y}\in\widetilde{V}(\mathbb{C}_{p})$} über $x$ und man setze $y=\varphi_{\overline{\mathbb{Q}}_{p}}(\widetilde{y})$, wobei \mbox{$\varphi_{\overline{\mathbb{Q}}_{p}}\colon\widetilde{Y}\to Y$} die induzierte Abbildung auf den generischen Fasern ist. Es folgt \mbox{$\varphi_{\overline{\mathbb{Q}}_{p}}(\gamma\widetilde{y})=\gamma y$} und aufgrund der Eigentlichkeit von $\mathcal{Y}$ und $\widetilde{\mathcal{Y}}$ über $\Spec\,\overline{\mathbb{Z}_{p}}$ gilt \mbox{$y_{\mathfrak{o}}=\varphi(\widetilde{y}_{\mathfrak{o}})$} und $(\gamma y)_{\mathfrak{o}}=\varphi((\gamma\widetilde{y})_{\mathfrak{o}})$.\par
Damit erhält man dasselbe Diagramm wie oben. Also folgt: \[(\gamma y)_{n}^{*}\circ(y_{n}^{*})^{-1}=(\gamma\widetilde{y})^{*}_{n}\circ(\widetilde{y}_{n}^{*})^{-1}\]
Dies zeigt, dass $\rho_{P,n}$ unabhängig von der Wahl der trivialisierenden Überdeckung $\pi\in S^{good}_{\mathfrak{X},D}$ ist.\par
\item[(iii)] Es ist klar, dass $\rho_{P,n}(\id)=\id$ für den trivialen Weg \mbox{$\id\in\Iso(F_{x},F_{x})$} gilt. Für Wege $\gamma\in\Iso(F_{x},F_{x'})$ und $\gamma'\in\Iso(F_{x'},F_{x''})$ sowie einen gewählten Punkt $y\in V(\mathbb{C}_{p})$ über $x$ liegt der Punkt $\gamma y$ über $x'$ und es gilt daher $\rho_{P,n}(\gamma)=(\gamma y)_{n}^{*}\circ(y_{n}^{*})^{-1}$ und $\rho_{P,n}(\gamma')=(\gamma'(\gamma y))_{n}^{*}\circ((\gamma y)_{n}^{*})^{-1}$. Daraus folgt: $\rho_{P,n}(\gamma')\circ\rho_{P,n}(\gamma)=((\gamma'\circ\gamma)(y))_{n}^{*}\circ(y_{n}^{*})^{-1}=\rho_{P,n}(\gamma'\circ\gamma)$. \hfill$\Box$\\}\end{itemize}}\par

Mit Hilfe dieser Funktoren $\rho_{P,n}$ können wir nun einen Funktor \[\rho_{P}\colon\Pi_{1}(U)\to\mathcal{P}(G(\mathfrak{o}))\] definieren:\par
Für jeden $\mathbb{C}_{p}$"=wertigen Punkt von $U$ setzt man $\rho_{P}(x)=P_{x_{\mathfrak{o}}}$; für einen étalen Weg $\gamma$ definiert man $\rho_{P}(\gamma)$ als $\rho_{P}(\gamma)=\projlim_{n}\rho_{P,n}(\gamma)\colon P_{x_{\mathfrak{o}}}\to P_{x'_{\mathfrak{o}}}$.\par
Dies ergibt einen wohldefinierten Funktor $\rho_{P}$, da man ähnlich wie im Beweis von \cite[Theorem 22]{DW1} zeigt:\\\par

\begin{lem}
Die Familie von Abbildungen $(\rho_{P,n}(\gamma))_{n\geq 1}$ definiert einen Morphismus zwischen den projektiven Systemen $(P_{x_{n}})_{n\geq 1}$ und $(P_{x'_{n}})_{n\geq 1}$.
\end{lem}
\bewub{Es genügt zu zeigen, dass die Abbildungen \[\rho_{P,n}\colon\Iso(F_{x},F_{x'})\to\Hom_{\mathfrak{o}_{n}}(P_{x_{n}},P_{x'_{n}})\] ein projektives System bezüglich der kanonischen Projektionen \[\lambda_{n+1}\colon\Hom_{\mathfrak{o}_{n+1}}(P_{x_{n+1}},P_{x'_{n+1}})\to\Hom_{\mathfrak{o}_{n+1}}(P_{x_{n+1}},P_{x'_{n+1}})\otimes_{\mathfrak{o}_{n+1}}\mathfrak{o}_{n}=\Hom_{\mathfrak{o}_{n}}(P_{x_{n}},P_{x'_{n}})\] bilden, so dass $\lambda_{n+1}\circ\rho_{P,n+1}=\rho_{P,n}$ ist.\\
Dazu wählt man für ein gegebenes $n\geq 1$ ein Element $(\pi\colon\mathcal{Y}\to\mathfrak{X})\in S_{\mathfrak{X},D}^{good}$, so dass $\pi_{n+1}^{*}P_{n+1}$ ein trivialer $G$"=Torseur ist. Dann ist auch $\pi_{n}^{*}P_{n}$ trivial. Für einen Punkt \mbox{$y\in V(\mathbb{C}_{p})$} über $x$ und einen étalen Weg \mbox{$\gamma\in\Iso(F_{x},F_{x'})$} betrachte man das folgende kommutierende Diagramm, wobei $a$ und $b$ die kanonischen Abbildungen sind:\par
\begin{equation*}
        \xymatrix@=3em{%
        \Spec\,\mathfrak{o}_{n} \ar[r]^{y_{n}} \ar[d]^{b} & \mathcal{Y}_{n} \ar[d]^{a} & \Spec\,\mathfrak{o}_{n} \ar[d]^{b} \ar[l]_{(\gamma y)_{n}} \\
        \Spec\,\mathfrak{o}_{n+1} \ar[r]^{y_{n+1}} & \mathcal{Y}_{n+1} & \Spec\,\mathfrak{o}_{n+1}. \ar[l]_{(\gamma y)_{n+1}}}
           \end{equation*}
Dieses Diagramm induziert wiederum folgendes kommutierendes Diagramm:\par
\begin{equation*}
        \xymatrix@=3em{%
        P_{x_{n+1}} \ar[d]^{b^{*}} & & \Gamma(\mathcal{Y}_{n+1}, \pi_{n+1}^{*}P_{n+1}) \ar[d]^{a^{*}} \ar[ll]^{y_{n+1}^{*}}_{\sim} \ar[rr]_{(\gamma y_{n+1})^{*}}^{\sim} & & P_{x'_{n+1}} \ar[d]^{b^{*}} \\
        P_{x_{n}} & & \Gamma(\mathcal{Y}_{n}, \pi_{n}^{*}P_{n}) \ar[ll]^{y_{n}^{*}}_{\sim} \ar[rr]_{(\gamma y_{n})^{*}}^{\sim} & & P_{x_{n}}.}
           \end{equation*}
Die Abbildungen $b^{*}$ sind gerade die kanonischen Reduktionsabbildungen von den topologischen Räumen mit stetiger $G(\mathfrak{o}_{n+1})$"=Aktion $P_{x_{n+1}}$ bzw. $P_{x'_{n+1}}$ zu den topologischen Räumen mit stetiger $G(\mathfrak{o}_{n})$"=Aktion $P_{x_{n}}$ bzw. $P_{x'_{n}}$. Damit ist die Abbildung $\rho_{P,n}(\gamma)=(\gamma y)_{n}^{*}\circ(y_{n}^{*})^{-1}$ die Reduktion modulo $p^{n}$ der Abbildung $\rho_{P,n+1}(\gamma)=(\gamma y)^{*}_{n+1}\circ(y_{n+1}^{*})^{-1}$.\hfill$\Box$}\\\par

Man bemerke, dass man für einen fest gewählten Punkt $x\in U(\mathbb{C}_{p})$ einen stetigen Homomorphismus $\rho_{P}\colon\pi_{1}(U,x)\to\Aut_{\mathcal{P}(G(\mathfrak{o}))}(P_{x_{\mathfrak{o}}})$ erhält.\\\par
Bezeichnet man mit $\Rep_{\pi_{1}(U)}(G(\mathfrak{o}))$ die Kategorie der stetigen Funktoren von $\pi_{1}(U)$ nach $\mathcal{P}(G(\mathfrak{o}))$, so erhält man insgesamt einen Funktor \[\rho\colon\mathscr{B}_{\mathfrak{X}_{\mathfrak{o},D}}(G)\to \Rep_{\pi_{1}(U)}(G(\mathfrak{o})),\] indem man ein Objekt $P$ aus $\mathscr{B}_{\mathfrak{X}_{\mathfrak{o},D}}(G)$ auf $\rho_{P}$ und einen Morphismus \mbox{$f\colon P\to P'$} auf die Familie $(f_{x_{\mathfrak{o}}})_{x\in U(\mathbb{C}_{p})}$ abbildet. Dabei sei $f_{x_{\mathfrak{o}}}\colon P_{x_{\mathfrak{o}}}\to P'_{x_{\mathfrak{o}}}$ die Abbildung $f_{x_{\mathfrak{o}}}=(x_{\mathfrak{o}}^{*}f)(\mathfrak{o})$.\par
Dies ergibt einen wohldefinierten Funktor, da die Familie von Morphismen topologischer Räume mit stetiger $G(\mathfrak{o})$"=Aktion \[f_{x_{\mathfrak{o}}}=x_{\mathfrak{o}}^{*}f\colon P_{x_{\mathfrak{o}}}\to P'_{x_{\mathfrak{o}}},x\in U(\mathbb{C}_{p})=\ob\Pi_{1}(U)\] eine natürliche Transformation von $\rho_{P}$ nach $\rho_{P'}$ definiert, die mit $\rho_{f}$ bezeichnet werde:\par
Sei nämlich $\gamma\in\mor_{\Pi_{1}(U)}(x,x')$ ein étaler Weg. Dann existiert nach Korollar \ref{dom} für alle $n\geq 1$ ein Element $\pi\colon\mathcal{Y}\to\mathfrak{X}$ in $S_{\mathfrak{X},D}^{good}$, so dass sowohl $\pi_{n}^{*}P_{n}$ als auch $\pi_{n}^{*}P'_{n}$ triviale Torseure sind.\par
Es bezeichne nun $f_{n}$ die Reduktion von $f$ modulo $p^{n}$ und man setze \mbox{$f_{x_{n}}:=x_{n}^{*}(f)$}. Außerdem wähle man einen Punkt $y$ über $x$ und definiere $y':=\gamma y$. Dann zeigt das kommutierende Diagramm \begin{equation*}
        \xymatrix@=3em{%
        P_{x_{n}} \ar[rr]^{f_{x_{n}}} & & P'_{x_{n}}\\
        \Gamma(\mathcal{Y}_{n},\pi_{n}^{*}P_{n}) \ar[u]^{y_{n}^{*}}_{\sim} \ar[rr]^{\Gamma(\mathcal{Y}_{n},\pi^{*}_{n}f_{n})} \ar[d]_{(y'_{n})^{*}}^{\sim} & & \Gamma(\mathcal{Y}_{n},\pi_{n}^{*}P'_{n}) \ar[u]_{y_{n}^{*}}^{\sim}  \ar[d]^{(y'_{n})^{*}}_{\sim}\\
        P_{x'_{n}} \ar[rr]^{f_{x'_{n}}} & & P'_{x'_{n}},\\}
           \end{equation*} 
           dass $f_{x'_{n}}\circ\rho_{P,n}(\gamma)=\rho_{P',n}(\gamma)\circ f_{x_{n}}$ gilt. Mittels Übergang zum projektiven Limes erhält man $f_{x'_{\mathfrak{o}}}\circ\rho_{P}(\gamma)=\rho_{P'}(\gamma)\circ f_{x_{\mathfrak{o}}}$. Also ist $\rho_{f}=(f_{x_{\mathfrak{o}}})$ ein Morphismus von $\rho_{P}$ zu $\rho_{P'}$.\par
\newpage

\section{Ausdehnung auf Prinzipalbündel über $X_{\mathbb{C}_{p}}$}

Für die Definition des Paralleltransports interessiert man sich eigentlich eher für Torseure auf $X_{\mathbb{C}_{p}}$ als auf $\mathfrak{X}_{\mathfrak{o}}$ für ein bestimmtes Modell $\mathfrak{X}$, wie wir das bislang getan haben. Deshalb definieren wir wie in \cite{DW3}:\\\par

\begin{defn}
Für ein glattes und affines Gruppenschema $G$ von endlicher Präsentation über $\mathfrak{o}$ definiert man die Kategorie $\mathscr{B}_{X_{\mathbb{C}_{p}}, D}(G)$ als die volle Unterkategorie aller $G$"=Rechtstorseure (d.\,h. $G_{X_{\mathbb{C}_{p}}}$"=Torseure) $P$ auf $X_{\mathbb{C}_{p}}$ mit der folgenden Eigenschaft:\par
Es gibt ein Modell $\mathfrak{X}$ von $X$ über $\overline{\mathbb{Z}_{p}}$ und einen $G$"=Torseur $\widetilde{P}$ auf $\mathfrak{X}_{\mathfrak{o}}$, der zu der Kategorie $\mathscr{B}_{\mathfrak{X}_{\mathfrak{o}},D}(G)$ gehört, so dass die $G$"=Torseure $P$ und $j^{*}\widetilde{P}$ über $X_{\mathbb{C}_{p}}$ isomorph sind, wobei \mbox{$j\colon X_{\mathbb{C}_{p}}\to\mathfrak{X}_{\mathfrak{o}}$} die kanonische Immersion sei.\par
\end{defn}\par

Insgesamt erhält man damit aus der obigen Konstruktion zusammen mit den gleichen Argumenten wie im Vektorbündelfall einen Funktor \[\rho\colon\mathscr{B}_{X_{\mathbb{C}_{p}}, D}(G)\to \Rep_{\pi_{1}(U)}(G(\mathbb{C}_{p}))\] 
in die Kategorie der stetigen Funktoren von $\pi_{1}(U)$ in die Kategorie $\mathcal{P}(G(\mathbb{C}_{p}))$ der topologischen Räume mit einer einfach"=transitiven und stetigen $G(\mathbb{C}_{p})$"=Aktion von rechts:\par
Für ein Objekt $P$ von $\mathscr{B}_{X_{\mathbb{C}_{p}}, D}(G)$ definiere man den stetigen Funktor \[\rho(P)=\rho_{P}\colon\Pi_{1}(U)\to\mathcal{P}(G(\mathbb{C}_{p}))\] wie folgt:\par Für $x\in U(\mathbb{C}_{p})=\ob\Pi_{1}(U)$ setze man $\rho_{P}(x)=P_{x}=x^{*}P$. Zum anderen definiert man für zwei Punkte $x,x'\in U(\mathbb{C}_{p})$ die stetige Abbildung \[\rho_{P}=\rho_{P,x,x'}\colon\mor_{\Pi_{1}(U)}(x,x')\to\Hom_{\mathbb{C}_{p}}(P_{x},P_{x'})\] durch $\rho_{P}(\gamma)=\psi^{-1}_{x'}\circ(\rho_{\widetilde{P}}(\gamma)\otimes_{\mathfrak{o}}\mathbb{C}_{p})\circ\psi_{x}$. \par
Dabei sei $\mathfrak{X}$ ein gewähltes Modell von $X$ über $\overline{\mathbb{Z}_{p}}$ und $\widetilde{P}$ ein $G$"=Torseur in der Kategorie $\mathscr{B}_{\mathfrak{X}_{\mathfrak{o}},D}(G)$ zusammen mit einem Isomorphismus \mbox{$\psi\colon P\to j^{*}\widetilde{P}$} von $G$"=Torseuren über $X_{\mathbb{C}_{p}}$, wobei $j\colon X_{\mathbb{C}_{p}}\to\mathfrak{X}_{\mathfrak{o}}$ die kanonische Immersion sei. Außerdem ist $\psi_{x}$ die Faserabbildung \[\psi_{x}=x^{*}(\psi)\colon P_{x}\xrightarrow{\sim}(j^{*}\widetilde{P})_{x}=\widetilde{P}_{x_{\mathfrak{o}}}\otimes_{\mathfrak{o}}\mathbb{C}_{p}=\widetilde{P}_{x_{\mathfrak{o}}}\otimes_{\mathbb{Z}}\mathbb{Q}.\] Für einen Morphismus $f\colon P\to P'$ in $\mathscr{B}_{\mathfrak{X}_{\mathbb{C}_{p}},D}(G)$ ist der Morphismus \[\rho(f)=\rho_{f}\colon\rho_{P}\to\rho_{P'}\] gegeben durch die Familie von Abbildungen \mbox{$f_{x}=x^{*}(f)\colon P_{x}\to P_{x'}$} für alle $x\in U(\mathbb{C}_{p})$.\\\par
Dass die Definition des Funktors unabhängig von den getroffenen Wahlen, d.\,h. von der Wahl des Modells $\mathfrak{X}$ und des $G$"=Torseurs $\widetilde{P}\in\mathscr{B}_{\mathfrak{X}_{\mathfrak{o},D}}(G)$, ist, folgt dabei aus den analogen Resultaten zu \cite[Proposition 26]{DW1} und \cite[Proposition 27]{DW1}: \\\par
Es sei $\alpha\colon\mathfrak{X}\to\mathfrak{X'}$ ein Morphismus von Modellen über $\overline{\mathbb{Z}_{p}}$ der glatten und projektiven Kurven $X$ bzw. $X'$ über $\overline{\mathbb{Q}_{p}}$ und $D'$ ein Divisor auf $X'$. Man setze außerdem $U'=X'-D'$ und \mbox{$U=X-\alpha^{*}D'$}. Dann definiert die generische Faser von $\alpha$ einen Funktor \[A(\alpha)\colon\Rep_{\Pi_{1}(U')}(G(\mathfrak{o}))\to\Rep_{\Pi_{1}(U)}(G(\mathfrak{o}))\] wie folgt:\par 
Für ein Objekt $\Gamma\in\Rep_{\Pi_{1}(U')}(G(\mathfrak{o}))$ definiert man $A(\alpha)(\Gamma)$ als die Komposition von Funktoren \[A(\alpha)(\Gamma)\colon\Pi_{1}(U)\xrightarrow{\alpha_{*}}\Pi_{1}(U')\xrightarrow{\Gamma}\mathcal{P}(G(\mathfrak{o})).\] Für einen Morphismus $f\colon\Gamma_{1}\to\Gamma_{2}$ in der Kategorie $\Rep_{\Pi_{1}(U')}(G(\mathfrak{o}))$, der durch die Familie von Abbildungen $f_{x'}\colon\Gamma_{1x'}\to\Gamma_{2x'}$ für $x'\in U'(\mathbb{C}_{p})$ gegeben ist, definieren wir $A(\alpha)(f)$ als die Familie von Abbildungen \[A(\alpha)(\Gamma_{1})_{x}=\Gamma_{1,\alpha(x)}\xrightarrow{f_{\alpha(x)}}\Gamma_{2,\alpha(x)}=A(\alpha)(\Gamma_{2})_{x}.\]
Dies definiert den gewünschten Funktor $A(\alpha)$. Ist $\alpha'\colon\mathfrak{X'}\to\mathfrak{X''}$ ein weiterer Morphismus von Modellen über $\overline{\mathbb{Z}_{p}}$ der glatten und projektiven Kurven $X'$ bzw. $X''$ über $\overline{\mathbb{Q}_{p}}$, so gilt $A(\alpha'\circ\alpha)=A(\alpha)\circ A(\alpha')$.\par

Als analoges Resultat zu \cite[Proposition 26]{DW1} erhält man in der so gegebenen Situation:\\\par
 
\begin{sat}[vgl. \protect{\cite[Proposition 3]{DW3}}]\label{unabhaengig}
Es seien die obigen Notationen gegeben. Dann induziert das Pullback entlang $\alpha$ einen Funktor $\alpha^{*}\colon\mathscr{B}_{\mathfrak{X'}_{\mathfrak{o}},D'}(G)\to\mathscr{B}_{\mathfrak{X}_{\mathfrak{o}},\alpha^{*}D'}(G)$ und das folgende Diagramm von Kategorien und Funktoren kommutiert (bis auf kanonische Isomorphismen):\par 
\begin{equation*}
        \xymatrix@=3em{%
        \mathscr{B}_{\mathfrak{X'}_{\mathfrak{o}},D'}(G) \ar[rr]^{\rho} \ar[d]^{\alpha^{*}}& & \Rep_{\Pi_{1}(U')}(G(\mathfrak{o})) \ar[d]^{A(\alpha)} \\
  \mathscr{B}_{\mathfrak{X}_{\mathfrak{o}},\alpha^{*}D}(G) \ar[rr]^{\rho} & & \Rep_{\Pi_{1}(U)}(G(\mathfrak{o})). \\}
           \end{equation*} 
           Insbesondere gilt für alle $\widetilde{P}\in\mathscr{B}_{\mathfrak{X'}_{\mathfrak{o}},D'}$, dass $\rho_{\alpha^{*}\widetilde{P}}=\rho_{\widetilde{P}}\circ\alpha_{*}$ als Funktoren von $\Pi_{1}(U)$ nach $\mathcal{P}(G(\mathfrak{o}))$ gilt.
\end{sat}
\bewub{Der Beweis ist wortwörtlich derselbe wie der von \cite[Proposition 26]{DW1}, wenn man $\mathcal{E}$ durch $\widetilde{P}$, "`Vektorbündel"' durch "`$G$"=Torseur"', die Kategorien $\mathscr{B}_{\mathfrak{X'}_{\mathfrak{o}},D'}$ und $\mathscr{B}_{\mathfrak{X}_{\mathfrak{o}},\alpha^{*}D'}$ durch die Kategorien $\mathscr{B}_{\mathfrak{X'}_{\mathfrak{o}},D'}(G)$ bzw. $\mathscr{B}_{\mathfrak{X}_{\mathfrak{o}},\alpha^{*}D'}(G)$ und den Verweis auf \cite[Proposition 9]{DW1} durch den auf Bemerkung \ref{pullback3} ersetzt.\newline\text{}\hfill$\Box$}\\\par

Insbesondere erhält man also für jeden Morphismus $f\colon\mathfrak{X}\to\mathfrak{X'}$ von Modellen von $X$ über $\overline{\mathbb{Z}_{p}}$, der auf der generischen Faser die Identität ist, ein kommutatives Diagramm \begin{equation*}
        \xymatrix@=3em{%
        \mathscr{B}_{\mathfrak{X'}_{\mathfrak{o}},D}(G) \ar[rrr]^{f^{*}} \ar[drr]_{\rho} & & & \mathscr{B}_{\mathfrak{X}_{\mathfrak{o}},D}(G) \ar[dl]_{\rho} 
\\ & &  \Rep_{\Pi_{1}(U)}(G(\mathfrak{o})) & . \\}
           \end{equation*} 
Zusammen mit dem kanonischen Funktor \[\Rep_{\Pi_{1}(U)}(G(\mathfrak{o}))\otimes\mathbb{Q}\to\Rep_{\Pi_{1}(U)}(G(\mathbb{C}_{p}))\]
ergibt dieses das folgende kommutative Diagramm:\\\par
 \begin{equation*}
        \xymatrix@=3em{%
        \mathscr{B}_{\mathfrak{X'}_{\mathfrak{o}},D}(G)\otimes\mathbb{Q} \ar[rrr]^{f^{*}} \ar[drr]_{\rho} & & & \mathscr{B}_{\mathfrak{X}_{\mathfrak{o}},D}(G)\otimes\mathbb{Q} \ar[dl]_{\rho} 
\\ & &  \Rep_{\Pi_{1}(U)}(G(\mathbb{C}_{p})) & . \\}
           \end{equation*} 
          Zusammen mit dem folgenden Resultat in Analogie zu \cite[Proposition 27]{DW1} erhält man daraus, dass der obige Funktor \[\rho\colon\mathscr{B}_{\mathfrak{X}_{\mathbb{C}_{p}},D}(G)\to\Rep_{\Pi_{1}(U)}(G(\mathbb{C}_{p}))\] wohldefiniert ist:\\\par
          
\begin{sat}
Es sei $X$ eine glatte und projektive Kurve über $\overline{\mathbb{Q}_{p}}$ mit Modellen $\mathfrak{X}_{1}$ und $\mathfrak{X}_{2}$ über $\overline{\mathbb{Z}_{p}}$. Dann existiert ein drittes Modell $\mathfrak{X}_{3}$ von $X$ zusammen mit Morphismen \[\mathfrak{X}_{1}\xleftarrow{p_{1}}\mathfrak{X}_{3}\xrightarrow{p_{2}}\mathfrak{X}_{2},\] die auf den generischen Fasern nach ihrer Identifikation mit $X$ jeweils die Identität ergeben. Für jeden Divisor $D$ auf $X$ hat man das folgende kommutative Diagramm von volltreuen Funktoren
\begin{equation*}
        \xymatrix@=3em{%
        \mathscr{B}_{\mathfrak{X}_{1\mathfrak{o}},D}(G)\otimes\mathbb{Q} \ar[dr]^{p_{1}^{*}} \ar[drrr]_{j_{\mathfrak{X}_{1\mathfrak{o}}}^{*}} & & & \\
         & \mathscr{B}_{\mathfrak{X}_{3\mathfrak{o}},D}(G)\otimes\mathbb{Q} \ar[rr]^{j_{\mathfrak{X}_{3\mathfrak{o}}}^{*}} & & 
 \mathscr{B}_{X_{\mathbb{C}_{p}}, D}(G),\\
\mathscr{B}_{\mathfrak{X}_{2\mathfrak{o}},D}\otimes\mathbb{Q} \ar[ur]_{p_{2}^{*}} \ar[urrr]_{j^{*}_{\mathfrak{X}_{2\mathfrak{o}}}}
& & &.\\}
           \end{equation*} wobei $j_{\mathfrak{X}_{i\mathfrak{o}}}\colon X_{\mathbb{C}_{p}}\to\mathfrak{X}_{i\mathfrak{o}}$ für \mbox{$i\in\{1,2,3\}$} jeweils die kanonische Immersion bezeichne.
\end{sat}
\bewub{Der Beweis ist wortwörtlich derselbe wie der von \cite[Proposition 27]{DW1}, wenn man wie oben $\mathcal{E}$ durch $\widetilde{P}$ und die Kategorien $\mathscr{B}_{\mathfrak{X}_{\mathfrak{o}},D}$ und $\mathscr{B}_{X_{\mathbb{C}_{p}},D}$ durch die Kategorien $\mathscr{B}_{\mathfrak{X}_{\mathfrak{o}},D}(G)$ bzw. $\mathscr{B}_{X_{\mathbb{C}_{p}}, D}(G)$ ersetzt. \hfill$\Box$}

\section{Verkleben der Funktoren $\rho_{P}$}

Im folgenden soll nun erläutert werden, wie die Funktoren \[\rho_{P}\colon\mathscr{B}_{X_{\mathbb{C}_{p}}, D}(G)\to\mathcal{P}(G(\mathbb{C}_{p}))\] verklebt werden können, wenn $P$ ein $G$"=Torseur auf $X_{\mathbb{C}_{p}}$ ist, der zu den Kategorien $\mathscr{B}_{X_{\mathbb{C}_{p}}, D}(G)$ für verschiedene Divisoren $D$ gehört.\par
Als erstes benötigen wir dazu das folgende Seifert"=van Kampen"=Theorem, das wir aus \cite{DW1} zitieren:\\\par

\begin{sat}[\protect{\cite[Proposition 34]{DW1}}]\label{kleben}
Es seien $U_{1}$ und $U_{2}$ zwei offene Unterschemata einer Kurve $X$ sowie \[i_{1}\colon U_{1}\cap U_{2}\into U_{1},\,\, i_{2}\colon U_{1}\cap U_{2}\into U_{2},\] und \[j_{1}\colon U_{1}\into U_{1}\cup U_{2},\,\,j_{2}\colon U_{2}\into U_{1}\cup U_{2}\] die zugehörigen offenen Immersionen, die das kommutierende Diagramm étaler Fundamentalgruppoide
\begin{equation*}
        \xymatrix@=3em{%
        \Pi_{1}(U_{1}\cap U_{2}) \ar[r]^{i_{1*}} \ar[d]^{i_{2*}} & \Pi_{1}(U_{1}) \ar[d]_{j_{1*}} \\
        \Pi_{1}(U_{2}) \ar[r]^{j_{2*}} & \Pi_{1}(U_{1}\cup U_{2})}
           \end{equation*}
liefern.\par
Dann existiert für jede Hausdorffsch"=topologische Kategorie $\mathcal{C}$ und alle stetigen Funktoren $\rho_{1}\colon\Pi_{1}(U_{1})\longrightarrow\mathcal{C}$ und $\rho_{2}\colon\Pi_{1}(U_{2})\longrightarrow\mathcal{C}$, für die \mbox{$\rho_{1}\circ i_{1*}=\rho_{2}\circ i_{2*}$} gilt, ein eindeutig bestimmter stetiger Funktor \[\rho\colon\Pi_{1}(U_{1}\cup U_{2})\longrightarrow\mathcal{C},\] so dass die Relationen $\rho\circ j_{1*}=\rho_{1}$ und $\rho\circ j_{2*}=\rho_{2}$ gelten.  
\end{sat}   

Damit können wir zeigen:\\\par

\begin{kor}\label{seifert}
Es seien $D_{1}$ und $D_{2}$ zwei Divisoren auf $X$ und man setze $U_{1}=X-D_{1}$ und $U_{2}=X-D_{2}$. Ferner sei $P$ ein $G$"=Torseur auf $X$, der sowohl zu $\mathscr{B}_{X_{\mathbb{C}_{p}}, D_{1}}(G)$ als auch $\mathscr{B}_{X_{\mathbb{C}_{p}}, D_{2}}(G)$ gehöre, und es seien $\rho_{P}^{1}$ bzw. $\rho_{P}^{2}$ die dazu konstruierten stetigen Funktoren $\Pi_{1}(U_{i})\to\mathcal{P}(G(\mathbb{C}_{p}))$ $(i=1,2)$.\par
Dann existiert ein eindeutig bestimmter stetiger Funktor \[\rho_{P}\colon\Pi(U_{1}\cup U_{2})\to\mathcal{P}(G(\mathbb{C}_{p})),\] der auf $\Pi_{1}(U_{i})$ für $i=1,2$ die Funktoren $\rho_{P}^{i}$ induziert. 
\end{kor}
\bewub{Da $G(\mathfrak{o})$ versehen mit der prodiskreten Topologie ein Hausdorff"=Raum ist, ist für zwei Objekte $P_{1}$ und $P_{2}$ aus $\mathcal{P}(G(\mathfrak{o}))$ nach Definition der Topologie auf $\mor_{\mathcal{P}(G(\mathfrak{o}))}(P_{1},P_{2})$ diese Menge mit ihrer Topologie ein Hausdorff"=Raum, so dass die Aussage aus Satz \ref{kleben} folgt. \hfill$\Box$}\\\par
\newpage

\section{Funktorialitäten und Vergleich mit dem Vektorbündelfall}

In diesem Abschnitt stellen wir einige Resultate zusammen, die in \cite{DW3} ohne Beweis genannt werden und die zeigen, dass die konstruierten Funktoren $\mathscr{B}_{\mathfrak{X}_{\mathfrak{o}},D}(G)\to\Rep_{\Pi_{1}(U)}(G(\mathfrak{o}))$ und $\mathscr{B}_{X_{\mathbb{C}_{p}}, D}(G)\to\Rep_{\Pi_{1}(U)}(G(\mathbb{C}_{p}))$ jeweils mit vielen natürlichen Operationen verträglich sind und in kanonischer Weise die im Falle der Vektorbündel konstruierten Funktoren ausehnen.\par
Als erstes zeigen wir die Verträglichkeit mit Galoiskonjugation:\par
Für einen $\mathbb{Q}_{p}$"=Automorphismus $\sigma$ von $\overline{\mathbb{Q}_{p}}$ definiert die Galoiskonjugation einen Funktor $\sigma_{*}$ von der Kategorie der $G$"=Torseure auf $\mathfrak{X}_{\mathfrak{o}}$ in die der ${}^{\sigma}G$"=Torseure auf ${}^{\sigma}\mathfrak{X}_{\mathfrak{o}}$, der jedem $G$"=Torseur $E$ auf $\mathfrak{X}_{\mathfrak{o}}$ den ${}^{\sigma}G$"=Torseur ${}^{\sigma}E$ auf ${}^{\sigma}\mathfrak{X}_{\mathfrak{o}}$ zuordnet. Dabei ist ${}^{\sigma}\mathfrak{X}=\mathfrak{X}\otimes_{\overline{\mathbb{Z}_{p}},\sigma}\overline{\mathbb{Z}_{p}}$ und damit $({}^{\sigma}\mathfrak{X})_{\mathfrak{o}}=\mathfrak{X}_{\mathfrak{o}}\otimes_{\mathfrak{o},\sigma}\mathfrak{o}$.
Man beachte ferner, dass ${}^{\sigma}E_{n}$ durch das Element ${}^{\sigma}\pi\in S_{{}^{\sigma}\mathfrak{X}_{\mathfrak{o}},{}^{\sigma}D}$ trivialisiert wird, wenn $E_{n}$ durch das Element $\pi\in S_{\mathfrak{X}_{\mathfrak{o}},D}$ trivialisiert wird, so dass die Galoiskonjugation insgesamt einen Funktor 
\[\sigma_{*}\colon\mathscr{B}_{\mathfrak{X}_{\mathfrak{o}},D}(G)\to\mathscr{B}_{{}^{\sigma}\mathfrak{X}_{\mathfrak{o}},{}^{\sigma}D}({}^{\sigma}G)\] liefert. Es ist klar, dass $(\sigma\tau)_{*}=\sigma_{*}\tau_{*}$ und $\id_{*}=\id$ gilt.\par
Andererseits lässt sich ein Funktor \[\sigma_{*}\colon\mathcal{P}(G(\mathfrak{o}))\to\mathcal{P}({}^{\sigma}G(\mathfrak{o}))\] mit denselben Eigenschaften wie folgt definieren:\par
Indem man einen Schnitt $s\in G(\mathfrak{o})$ auf den Schnitt \[\sigma\circ s\circ\sigma^{-1}=\sigma\circ s\circ\spec\sigma\in{}^{\sigma}G(\mathfrak{o})\] abbildet, erhält man einen Isomorphismus topologischer Räume \[\sigma\colon G(\mathfrak{o})\xrightarrow{\sim}{}^{\sigma}G(\mathfrak{o}).\] Für einen $G(\mathfrak{o})$"=Raum $Z$ in der Kategorie $\mathcal{P}(G(\mathfrak{o}))$ definiert man nun $\sigma_{*}(Z)$ in der Kategorie $\mathcal{P}({}^{\sigma}G(\mathfrak{o}))$ als den topologischen Raum $Z$ zusammen mit der induzierten Aktion von ${}^{\sigma}G$, die durch $z\cdot g:=z({}^{\sigma^{-1}}g)$ für alle $z\in Z$ und alle $g\in{}^{\sigma}G(\mathfrak{o})$ gegeben ist. Mit $\sigma\colon Z\to{}^{\sigma}Z$ sei die Identitätsabbildung auf den unterliegenden topologischen Räumen bezeichnet. Setzen wir noch \mbox{$\sigma_{*}(\varphi)=\varphi$} auf der Ebene der Morphismen, so erhält man insgesamt den gewünschten wohldefinierten Funktor \[\sigma_{*}\colon\mathcal{P}(G(\mathfrak{o}))\to\mathcal{P}({}^{\sigma}G(\mathfrak{o})).\]
Außerdem betrachten wir die Wirkung der Galoiskonjugation auf die stetigen Funktoren von $\Pi_{1}(U)$ in die Kategorie $\mathcal{P}(G(\mathfrak{o}))$, die in der Kategorie $\Rep_{\Pi_{1}(U)}(G(\mathfrak{o}))$ liegen, wobei wir wie folgt einen Funktor \[C_{\sigma}\colon\Rep_{\Pi_{1}(U)}(G(\mathfrak{o}))\to\Rep_{\Pi_{1}({}^{\sigma}U)}({}^{\sigma}G(\mathfrak{o}))\] erhalten:\par
Für Objekte $\mathcal{Z}\colon\Pi_{1}(U)\to\mathcal{P}(G(\mathfrak{o}))$ setzt man $C_{\sigma}(\mathcal{Z}):=\sigma_{*}\circ\mathcal{Z}\circ\sigma_{*}^{-1}$, für Morphismen $f\colon\mathcal{Z}\to\mathcal{Z'}$ in $\Rep_{\Pi_{1}(U)}(G(\mathfrak{o}))$  definiert man $C_{\sigma}(f)$ als \linebreak \mbox{$C_{\sigma}(f):=(\sigma\circ f_{x}\circ\sigma^{-1})_{x\in U(\mathbb{C}_{p})}$}. Wie oben hat man kanonische Isomorphismen $C_{\tau\sigma}=C_{\tau}\circ C_{\sigma}$ und $C_{\id}=\id$.\par
Mit diesen Definitionen gilt nun das folgende Resultat in Analogie zu \cite[Proposition 25]{DW1}:\\\par

\begin{sat}\label{galoiskonj}
In der gegebenen Situation ist das Diagramm \begin{equation*}
        \xymatrix@=3em{%
        \mathscr{B}_{\mathfrak{X}_{\mathfrak{o}},D}(G) \ar[r]^{\rho^{\mathfrak{X}}} \ar[d]_{\sigma_{*}} & \Rep_{\Pi_{1}(U)}(G(\mathfrak{o})) \ar[d]^{C_{\sigma}}\\
        \mathscr{B}_{{}^{\sigma}\mathfrak{X}_{\mathfrak{o}},D}({}^{\sigma}G) \ar[r]^{\rho^{{}^{\sigma}\mathfrak{X}}} & \Rep_{\Pi_{1}({}^{\sigma}U)}({}^{\sigma}G(\mathfrak{o}))}
           \end{equation*} von Kategorien und Funktoren bis auf kanonische Isomorphismen von Funktoren kommutativ.\par
Insbesondere gilt $\rho_{{}^{\sigma}P}=\sigma_{*}\circ\rho_{P}\circ\sigma_{*}^{-1}$ als Funktoren von $\Pi_{1}({}^{\sigma}U)$ nach $\mathcal{P}({}^{\sigma}G(\mathfrak{o}))$ für jeden $G$"=Torseur $P\in\mathscr{B}_{\mathfrak{X}_{\mathfrak{o}},D}(G)$.
\end{sat}
\bewub{Es werde zunächst die Kommutativität des Diagramms auf Ebene der Objekte nachgeprüft:\par
Sei also $P\in\mathscr{B}_{\mathfrak{X}_{\mathfrak{o}},D}(G)$. Dann ist zu zeigen, dass \[\rho_{{}^{\sigma}P}=C_{\sigma}(\rho_{P})=\sigma_{*}\circ\rho_{P}\circ\sigma_{*}^{-1}\] als Funktoren von $\Pi_{1}({}^{\sigma}U)$ nach $\mathcal{P}({}^{\sigma}G(\mathfrak{o}))$ gilt.\par Für ein Objekt $x\in {}^{\sigma}U(\mathbb{C}_{p})$ folgt aus den obigen Definitionen und denen in §1.3, insbesondere aus der Tatsache, dass Galoiskonjugation mit Pullback vertauscht und also $({}^{\sigma}P)_{x_{\mathfrak{o}}}=\sigma_{*}(P_{(\sigma^{-1}(x))_{\mathfrak{o}}})$ gilt, dass \[({}^{\sigma}P)_{x_{\mathfrak{o}}}=(\sigma_{*}\circ\rho_{P}\circ\sigma_{*}^{-1})(x)\] ist und also \[\rho_{{}^{\sigma}P}(x)=(\sigma_{*}\circ\rho_{P}\circ\sigma_{*}^{-1})(x)\] gilt. \par Ist $\gamma$ ein étaler Weg in $\Pi_{1}({}^{\sigma}U)$ von $x$ nach $x'$, so folgt ebenfalls aus den genannten Definitionen, dass \[\rho_{{}^{\sigma}P}(\gamma)=(\sigma_{*}\circ\rho_{P}\circ\sigma_{*}^{-1})(\gamma)\] und also insgesamt \[\rho_{{}^{\sigma}P}=C_{\sigma}(\rho_{P})=\sigma_{*}\circ\rho_{P}\circ\sigma_{*}^{-1}\] gilt.\par
Es bleibt noch zu zeigen, dass \[(\rho^{{}^{\sigma}\mathfrak{X}}\circ\sigma_{*})(\varphi)=C_{\sigma}(\rho^{\mathfrak{X}}(\varphi))=(\sigma\circ\varphi_{x}\circ\sigma^{-1})_{x\in U(\mathbb{C}_{p})}\] für jeden Morphismus $\varphi\colon P\to P'$ in $\mathscr{B}_{\mathfrak{X}_{\mathfrak{o}},D}(G)$ gilt, was ebenfalls mit etwas Schreibaufwand direkt aus den obigen Definitionen der Funktoren und §1.3. folgt. \hfill$\Box$} \\\par

\begin{bem}\label{bemgaloisconj}
Ist $K$ eine algebraische Erweiterung von $\mathbb{Q}_{p}$ in $\overline{\mathbb{Q}_{p}}$ und sind $\mathfrak{X}, X, D$ und $G$ bereits über $\mathfrak{o}_{K}$ bzw. $K$ definiert, so können ${}^{\sigma}\mathfrak{X}, {}^{\sigma}X, {}^{\sigma}D$ und ${}^{\sigma}G$ über $\overline{\mathbb{Z}_{p}}$ bzw. $\overline{\mathbb{Q}_{p}}$ bzw. $\mathfrak{o}$ für alle \mbox{$\sigma\in\Gal(\overline{\mathbb{Q}_{p}}/K)$} mit $\mathfrak{X}, X, D$ bzw. $G$ identifiziert werden. In diesem Fall kommutiert dann der Funktor \[\rho\colon\mathscr{B}_{\mathfrak{X}_{\mathfrak{o}},D}(G)\to\Rep_{\Pi_{1}(U)}(G(\mathfrak{o}))\] mit den Aktionen von $\Gal(\overline{\mathbb{Q}_{p}}/K)$ auf $\mathfrak{X},X, D$ und $G$ vermittels $\sigma_{*}$ und $C_{\sigma}$.\\
\end{bem}\par

Als nächstes sei daran erinnert, dass nach Satz \ref{unabhaengig} die Konstruktion des Funktors \[\mathscr{B}_{\mathfrak{X}_{\mathfrak{o}},D}(G)\to\Rep_{\Pi_{1}(U)}(G(\mathfrak{o}))\] mit dem Wechsel der Kurve $X$ und ihres Modells $\mathfrak{X}$ verträglich ist, d.h. mit Morphismen $\alpha\colon\mathfrak{X}\to\mathfrak{X'}$ über $\overline{\mathbb{Z}_{p}}$ von Modellen der glatten und projektiven Kurven $X$ bzw. $X'$ über $\overline{\mathbb{Q}_{p}}$ vertauscht.\\\par
Gleiches gilt für den Wechsel des Gruppenschemas $G$:\\\par
Es sei $\varphi\colon G\to G'$ ein Morphismus glatter, affiner Gruppenschemata von endlicher Präsentation über $\mathfrak{o}$. Dann ist für jeden $G$"=Torseur $P$ auf einem $\mathfrak{o}$"=Schema $\xi$ die fppf"=Garbe $\varphi_{*}P=P\times^{G}G'$ auf $\xi$ ein Garbentorseur unter $G'$ und damit nach \cite[Theorem III.4.3a)]{M} ein $G'$"=Torseur, da $G'$ als affin gegeben ist. Setzt man für jeden Morphismus $\alpha\colon P\to P'$ von $G$"=Torseuren auf $\xi$ den induzierten Morphismus $\varphi_{*}P\to\varphi_{*}P'$ von $G'$"=Torseuren als $\varphi_{*}(\alpha)$, so erhält man einen wohldefinierten Funktor $\varphi_{*}$ von der Kategorie der $G$"=Torseure auf $\xi$ in die Kategorie der $G'$"=Torseure auf $\xi$.\par
Nach Konstruktion gilt \[f^{*}\varphi_{*}(P)=\varphi_{*}(f^{*}(P))\] für jeden Morphismus \mbox{$f\colon\widetilde{\xi}\to\xi$} von $\mathfrak{o}$"=Schemata. Daraus folgt, dass sich $\varphi_{*}$ zu einem Funktor \[\varphi_{*}\colon\mathscr{B}_{\mathfrak{X}_{\mathfrak{o}},D}(G)\to\mathscr{B}_{\mathfrak{X}_{\mathfrak{o}},D}(G')\] einschränkt.\par
Außerdem erhält man einen Funktor \[\varphi_{*}\colon\mathcal{P}(G(\mathfrak{o}))\to\mathcal{P}(G'(\mathfrak{o})),\] indem man einen $G(\mathfrak{o})$"=Raum $\mathcal{Z}$ auf den $G'(\mathfrak{o})$"=Raum \mbox{$\varphi_{*}(\mathcal{Z})=\mathcal{Z}\times^{G(\mathfrak{o})}G'(\mathfrak{o})$} und jeden Morphismus $\mathcal{Z}\to\mathcal{Z'}$ auf den induzierten Morphismus \mbox{$\varphi_{*}\mathcal{Z}\to\varphi_{*}\mathcal{Z'}$} abbildet.\par
Schließlich definiert man noch einen Funktor \[\varphi_{*}\colon\Rep_{\Pi_{1}(U)}(G(\mathfrak{o}))\to\Rep_{\Pi_{1}(U)}(G'(\mathfrak{o})),\] indem man ein Objekt $\mathcal{Z}$ auf $\varphi_{*}\circ\mathcal{Z}$ und jeden Morphismus $f=(f_{x})$ auf $(\varphi_{*}(f_{x}))$ abbildet.\par
Mit diesen Definitionen findet man das folgende Resultat:\\\par

\begin{sat}\label{gruppenwechsel}
In der gegebenen Situation ist das Diagramm \begin{equation*}
        \xymatrix@=3em{%
        \mathscr{B}_{\mathfrak{X}_{\mathfrak{o}},D}(G) \ar[r]^{\rho^{\mathfrak{X}}} \ar[d]_{\varphi_{*}} & \Rep_{\Pi_{1}(U)}(G(\mathfrak{o})) \ar[d]^{\varphi_{*}}\\
        \mathscr{B}_{\mathfrak{X}_{\mathfrak{o}},D}(G') \ar[r]^{\rho^{\mathfrak{X}}} & \Rep_{\Pi_{1}(U)}(G'(\mathfrak{o}))}
           \end{equation*} von Kategorien und Funktoren bis auf kanonische Isomorphismen von Funktoren kommutativ.\par
Insbesondere gilt $\rho_{\varphi_{*}P}=\varphi_{*}\circ\rho_{P}$ als Funktoren von $\Pi_{1}(U)$ nach $\mathcal{P}(G'(\mathfrak{o}))$ für jeden $G$"=Torseur $P\in\mathscr{B}_{\mathfrak{X}_{\mathfrak{o}},D}(G)$.
\end{sat}
\bewub{Es sei $P\in\mathscr{B}_{\mathfrak{X}_{\mathfrak{o}},D}(G)$. Dann gilt \[(\varphi_{*}\circ\rho)(P)(x)=(\varphi_{*}\circ\rho_{P})(x)=P_{x_{\mathfrak{o}}}\times^{G(\mathfrak{o})}G'(\mathfrak{o})=(x_{\mathfrak{o}}^{*}P(\mathfrak{o}))\times^{G(\mathfrak{o})}G'(\mathfrak{o})\] und \[(\rho\circ\varphi_{*})(P)(x)=\rho_{P\times^{G}G'}(x)=(P\times^{G}G')_{x_{\mathfrak{o}}}=(x_{\mathfrak{o}}^{*}P)(\mathfrak{o})\] für alle $x\in\Pi_{1}(U)$ jeweils als topologische Räume mit einer einfach"=transitiven und stetigen Aktion von $G'(\mathfrak{o})$ von rechts. Weil Pullback unter $x_{\mathfrak{o}}$ mit der Bildung des eingeschränkten Produktes vertauscht und für ein $\mathfrak{o}$"=Schema $A$ gilt, dass $(A\times^{G}G')(\mathfrak{o})$ und $A(\mathfrak{o})\times^{G(\mathfrak{o})}G'(\mathfrak{o})$ als topologische Räume mit einer einfach"=transitiven und stetigen Aktion von $G'(\mathfrak{o})$ von rechts kanonisch isomorph sind, folgt damit \[(\varphi_{*}\circ\rho)(P)(x)=(\rho\circ\varphi_{*})(P)(x)\] für alle $x\in \Pi_{1}(U)$.\par
Aus demselben Argument folgt für jeden étalen Weg $\gamma$ von $x$ nach $x'$ in $\Pi_{1}(U)$, dass \[(\varphi_{*}\circ\rho)(P)(\gamma)=\varphi_{*}(\rho_{P}(\gamma))=\rho_{P\times^{G}G'}(\gamma)=(\rho\circ\varphi_{*})(P)(\gamma)\] gilt. \par
Also gilt auf Ebene der Objekte \[(\rho\circ\varphi)(P)=(\varphi_{*}\circ\rho)(P)\] für alle \mbox{$P\in\mathscr{B}_{\mathfrak{X}_{\mathfrak{o}},D}(G)$}.\par
Es sei nun $\alpha\colon P\to P'$ ein Morphismus in $\mathscr{B}_{\mathfrak{X}_{\mathfrak{o}},D}(G)$. Dann ergibt sich wiederum aus der Tatsache, dass die Bildung des eingeschränkten Produkts $\_\times^{G}G'$ mit Pullback vertauscht, dass \[(\rho\circ\varphi_{*})(\alpha)=((\alpha\times^{G}G')_{x_{\mathfrak{o}}})_{x\in U(\mathbb{C}_{p})}=(\varphi_{*}(\alpha_{x_{\mathfrak{o}}}))_{x\in U(\mathbb{C}_{p})}=(\varphi_{*}\circ\rho)(\alpha)\] gilt, so dass insgesamt die Kommutativität des Diagramms folgt. \hfill $\Box$}\\\par

Mittels einer ähnlichen Konstruktion kann man zudem von $G$"=Torseuren zu Vektorbündeln wechseln:\par
Dazu sei $\Gamma$ ein freier $\mathfrak{o}$"=Modul von endlichem Rang und $\bf{\Gamma}$ der assoziierte affine Raum über $\Spec\,\mathfrak{o}$. Dann versteht man unter einer algebraischen Darstellung von $G$ auf $\Gamma$ einen Homomorphismus von $\mathfrak{o}$"=Gruppenschemata \[G\to\Gl(\bf{\Gamma}).\] Im folgenden schreiben wir kurz $\bf{\Gamma}$ für diese Darstellung und definieren $\Rep_{G}(\mathfrak{o})$ und $\Rep_{G(\mathfrak{o})}(\mathfrak{o})$ als die Kategorie der algebraischen bzw. stetigen Darstellungen von $G$ bzw. $G(\mathfrak{o})$ auf freien $\mathfrak{o}$"=Moduln von endlichem Rang. \par
Für jeden $G$"=Torseur $P$ auf einem $\mathfrak{o}$"=Schema $\xi$ ist dann die fppf"=Garbe \linebreak \mbox{$E=P\times^{G}\bf{\Gamma}$} ein Vektorbündel auf $\xi$. Diese Konstruktion ist funktoriell in $P$ und der Darstellung $\bf{\Gamma}$; für einen Morphismus $f\colon\widetilde{\xi}\to\xi$ gilt \[f^{*}(P\times^{G}{\bf{\Gamma}})=f^{*}P\times^{G}\bf{\Gamma}.\] Unter Verwendung der Notation aus \cite{DW1} folgt daraus, dass sich $\_\times^{G}\bf{\Gamma}$ zu einem Funktor \[\_\times^{G}{\bf{\Gamma}}\colon\normalfont\mathscr{B}_{\mathfrak{X}_{\mathfrak{o}},D}(G)\to\mathscr{B}_{\mathfrak{X}_{\mathfrak{o}},D}\] einschränkt, so dass man einen Bifunktor \[\_\times^{G}\_\colon\mathscr{B}_{\mathfrak{X}_{\mathfrak{o}},D}(G)\times\Rep_{G}(\mathfrak{o})\to\mathscr{B}_{\mathfrak{X}_{\mathfrak{o}},D}\] erhält. (Wie in \cite{DW1} ist dabei $\mathscr{B}_{\mathfrak{X}_{\mathfrak{o}},D}$ die Kategorie aller Vektorbündel $E$ auf $\mathfrak{X}_{\mathfrak{o}}$ mit der Eigenschaft, dass für alle $n\in\nat$ ein Element $\pi\in S_{\mathfrak{X},D}$ existiert, so dass $\pi_{n}^{*}E$ trivial ist.)\par
Andererseits gibt es einen Bifunktor \[\_\times^{G(\mathfrak{o})}\_\colon\Rep_{\Pi_{1}(U)}(G(\mathfrak{o}))\times\Rep_{G(\mathfrak{o})}(\mathfrak{o})\to\Rep_{\Pi_{1}(U)}(\mathfrak{o}),\] der wie folgt definiert ist:\par
Ein Paar $(\mathcal{Z},\Gamma)$ aus einem stetigen Funktor $\mathcal{Z}\colon\Pi_{1}(U)\to\mathcal{P}(G(\mathfrak{o}))$ und einer stetigen Darstellung von $\Gamma$ auf $G(\mathfrak{o})$ wird unter $\_\times^{G}\_$ auf den wie folgt definierten Funktor \[\mathcal{Z}\times^{G(\mathfrak{o})}\Gamma\colon\Pi_{1}(U)\to\Mod_{\mathfrak{o}}\] abgebildet:\par
Für $x\in U(\mathbb{C}_{p})=\Pi_{1}(U)$ setzt man $(\mathcal{Z}\times^{G(\mathfrak{o})}\Gamma)(x)=\mathcal{Z}(x)\times^{G(\mathfrak{o})}\Gamma$. Für einen étalen Weg $\gamma$ von $x$ nach $x'$ definiert man $(\mathcal{Z}\times^{G(\mathfrak{o})}\Gamma)(\gamma)$ als die $\mathfrak{o}$"=lineare Abbildung von $\mathcal{Z}(x)\times^{G(\mathfrak{o})}\Gamma$ nach $\mathcal{Z}(x')\times^{G(\mathfrak{o})}\Gamma$, die durch die stetige $G(\mathfrak{o})$"=lineare Abbildung $\mathcal{Z}(\gamma)\colon\mathcal{Z}(x)\to\mathcal{Z}(x')$ induziert wird.\par
Morphismen in $\Rep_{\Pi_{1}(U)}(G(\mathfrak{o}))$ bzw. $\Rep_{G(\mathfrak{o})}(\mathfrak{o})$ werden in offensichtlicher Weise abgebildet.\par
Wie in \cite{DW1} ist dabei $\Rep_{\Pi_{1}(U)}(\mathfrak{o})$ die $\mathfrak{o}$"=lineare Kategorie von stetigen Funktoren von $\Pi_{1}(U)$ in die Kategorie der freien $\mathfrak{o}$"=Modulen von endlichem Rang.\\\par

\begin{sat}\label{darstellung}
In der obigen Situation kommutiert das folgende Diagramm:\par
\begin{equation*}
        \xymatrix@=3em{%
        \mathscr{B}_{\mathfrak{X}_{\mathfrak{o}},D}(G)\times\Rep_{G}(\mathfrak{o}) \ar[rrr]^{\_\times^{G}\_}\ar[d]_{\rho\times\text{nat}} &&& \mathscr{B}_{\mathfrak{X}_{\mathfrak{o}},D} \ar[d]^{\rho}\\
        \Rep_{\Pi_{1}(U)}(G(\mathfrak{o}))\times\Rep_{G(\mathfrak{o})}(\mathfrak{o})   \ar[rrr]^{\_\times^{G(\mathfrak{o})}\_} & &&\Rep_{\Pi_{1}(U)}(\mathfrak{o}).}
           \end{equation*}
           Insbesondere gilt \[\rho_{P\times^{G}\bf{\Gamma}}=\rho_{P}\times^{G(\mathfrak{o})}\Gamma\]  als Funktoren von $\Pi_{1}(U)$ nach $\Mod_{\mathfrak{o}}$ für alle $P\in\mathscr{B}_{\mathfrak{X}_{\mathfrak{o}},D}(G)$ und $\Gamma\in\Rep_{G}(\mathfrak{o})$.
\end{sat}
\bewub{Auf Ebene der Objekte ist zu zeigen, dass \[\rho_{P\times^{G}\bf{\Gamma}}=\rho_{P}\times^{G(\mathfrak{o})}\Gamma\]  als Funktoren von $\Pi_{1}(U)$ nach $\Mod_{\mathfrak{o}}$ für alle $P\in\mathscr{B}_{\mathfrak{X}_{\mathfrak{o}},D}(G)$ und $\Gamma\in\Rep_{G}(\mathfrak{o})$ gilt.\par
Ist also $x\in U(\mathbb{C}_{p})$, so folgt aus der Relation $(P\times^{G}{\bf{\Gamma}})_{x_{\mathfrak{o}}}=P_{x_{\mathfrak{o}}}\times^{G(\mathfrak{o})}{\bf{\Gamma}}$, dass $(\rho_{P}\times^{G(\mathfrak{o})}{\bf{\Gamma}})(x)=\rho_{P\times^{G}{\bf{\Gamma}}}(x)$ gilt.\par
Ist ferner $\gamma$ ein étaler Weg von $x$ nach $x'$ in $\Pi_{1}(U)$, so folgt aus der genannten Relation, dass die Abbildung $(P\times^{G}  {\bf{\Gamma}})_{x_{\mathfrak{o}}}\to(P\times^{G}{\bf{\Gamma}})_{x'_{\mathfrak{o}}}$ gleich der von \mbox{$P_{x_{\mathfrak{o}}}\to P_{x'_{\mathfrak{o}}}$} induzierten $\mathfrak{o}$"=linearen Abbildung von $P_{x_{\mathfrak{o}}}\times^{G(\mathfrak{o})}{\bf{\Gamma}}$ nach $P_{x'_{\mathfrak{o}}}\times^{G(\mathfrak{o})}{\bf{\Gamma}}$ ist, woraus $(\rho_{P}\times^{G(\mathfrak{o})}{\bf{\Gamma}})(\gamma)=\rho_{P\times^{G}{\bf{\Gamma}}}(\gamma)$ folgt.\par
Mit den gleichen Überlegungen sieht man, dass die Funktoren \mbox{$\rho\circ\_\times^{G}\_$} und $\_\times^{G(\mathfrak{o})}\_\circ(\rho\times\text{nat})$ auch auf Morphismen in $\mathscr{B}_{\mathfrak{X}_{o},D}(G)$ bzw. in $\Rep_{G}(\mathfrak{o})$ übereinstimmen. \hfill$\Box$}\\\par

\begin{bem}
Zu jedem gegebenen $G$"=Torseur $P$ in $\mathscr{B}_{\mathfrak{X}_{\mathfrak{o}},D}(G)$ erhält man mittels des Bifunktors $\_\times^{G}\_$ ein Vektorbündel in $\mathscr{B}_{\mathfrak{X}_{\mathfrak{o}},D}$ für jede algebraische Darstellung $\bf{\Gamma}$ von $G$. Der vorangehende Satz zeigt somit, wie man den Paralleltransport entlang étaler Wege zwischen den Fasern des Vektorbündels über den Paralleltransport des Torseurs berechnen kann. \\
\end{bem}\par

Im folgenden erklären wir kurz, wie wir ein $\Gl_{n}$"=Prinzipalbündel in \linebreak $\mathscr{B}_{\mathfrak{X}_{\mathfrak{o}},D}(\Gl_{n})$ aus einem Vektorbündel vom Rang $n$ in $\mathscr{B}_{\mathfrak{X}_{\mathfrak{o}},D}$ erhalten. Ist $E$ ein Vektorbündel vom Rang $n$ auf einem Schema $\xi$, so ist die fppf"=Garbe $\Iso_{\xi}(E,\mathbb{A}^{n}_{\xi})$ der lokalen Isomorphismen von $E$ zu dem trivialen Vektorbündel $\mathbb{A}_{\xi}^{n}$ auf $\xi$ ein Rechts"=Garbentorseur unter $\Gl_{n,\xi}$. Da $\Gl_{n,\xi}$ affin über $\xi$ ist, ist dieser Garbentorseur durch einen $\Gl_{n}$"=Torseur ${\bf{P}}(E)$ darstellbar, das sogenannte Rahmenbündel von $E$. Ordnet man jedem Morphismus $E\to E'$ den induzierten Morphismus $\Iso_{\xi}(E,\mathbb{A}^{n}_{\xi})\to\Iso_{\xi}(E',\mathbb{A}^{n}_{\xi})$ zu, so erhält man einen wohldefinierten Funktor $\bf{P}$ von der Kategorie der Vektorbündel vom Rang $n$ auf $\xi$ in die Kategorie der $\Gl_{n}$"=Torseure auf $\xi$. Aufgrund der Konstruktion des Funktors gilt $f^{*}{\bf{P}}(E)={\bf{P}}(f^{*}E)$ für jeden Morphismus $f\colon\widetilde{\xi}\to\xi$ in funktorieller Weise. Daraus folgt, dass sich $\bf{P}$ zu einem Funktor von der vollen Unterkategorie $\mathscr{B}_{\mathfrak{X}_{\mathfrak{o}},D}^{(n)}$ der Vektorbündel vom Rang $n$ in $\mathscr{B}_{\mathfrak{X}_{\mathfrak{o}},D}$ in die Kategorie $\mathscr{B}_{\mathfrak{X}_{\mathfrak{o}},D}(\Gl_{n})$ einschränkt. \par
Auf der anderen Seite definieren wir einen Funktor $\bf{P}$ von der vollen Unterkategorie $\Rep_{\Pi_{1}(U)}^{n}(\mathfrak{o})$ von $\Rep_{\Pi_{1}(U)}(\mathfrak{o})$, die aus den Darstellungen vom Rang $n$ besteht, nach $\Rep_{\Pi_{1}(U)}(\Gl_{n}(\mathfrak{o}))$ wie folgt:\par
Für ein Objekt $\bf{\Gamma}$ aus $\Rep_{\Pi_{1}(U)}^{n}(\mathfrak{o})$, d.\,h. einen stetigen Funktor \[{\bf{\Gamma}}\colon\Pi_{1}(U)\to\Mod_{\mathfrak{o}}^{(n)},\] definiert man den stetigen Funktor \[{\bf{P}}({\bf{\Gamma}})\colon\Pi_{1}(U)\to\mathcal{P}(G(\mathfrak{o})),\] indem man ${\bf{P}}(\Gamma)(x)=\Iso_{\mathfrak{o}}({\bf}{\Gamma}(x),\mathfrak{o}^{n})$ für jeden Punkt \mbox{$x\in U(\mathbb{C}_{p})$} setzt, wobei $\Iso_{\mathfrak{o}}(\Gamma(x),\mathfrak{o}^{n})$ mit der in natürlicher Weise induzierten Topologie und der Aktion von $\Gl_{n}(\mathfrak{o})$ von rechts versehen sei. Für einen étalen Weg $\gamma$ von $x$ nach $x'$ setzt man \[{\bf{P}}({\bf{\Gamma}})(\gamma)=({\bf{\Gamma}}(\gamma)^{-1})^{*}\colon\Iso_{\mathfrak{o}}({\bf{\Gamma}}(x),\mathfrak{o}^{n})\to\Iso_{\mathfrak{o}}({\bf{\Gamma}}(x'),\mathfrak{o}^{n}),\varphi\mapsto\varphi\circ\Gamma(\gamma)^{-1},\] so dass man insgesamt einen wohldefinierten Funktor $\bf{P}(\bf{\Gamma})$ erhält.\par
Diese Definition des Funktors $\bf{P}$ auf der Ebene der Objekte setzt sich in natürlicher Weise auf Morphismen $\alpha\colon\bf{\Gamma}\to\bf{\Gamma'}$ in $\Rep_{\Pi_{1}(U)}^{(n)}(\mathfrak{o})$ fort, indem man die von $\alpha$ induzierte Abbildung von $\Iso_{\mathfrak{o}}({\bf{\Gamma}}(x),\mathfrak{o}^{n})$ nach $\Iso_{\mathfrak{o}}({\bf{\Gamma'}}(x),\mathfrak{o}^{n})$ für alle \mbox{$x\in U(\mathbb{C}_{p})$} als ${\bf{P}}(\alpha)(x)$ definiert, was eine natürliche Transformation $\bf{P}(\alpha)\colon\bf{P}(\bf{\Gamma})\to\bf{P}(\bf{\Gamma'})$ liefert.\\\par

\begin{sat}\label{pblvbl}
In der gegebenen Situation ist das Diagramm \begin{equation*}
        \xymatrix@=3em{%
        \mathscr{B}_{\mathfrak{X}_{\mathfrak{o}},D}^{(n)} \ar[r]^{\rho^{\mathfrak{X}}} \ar[d]_{\bf{P}} & \Rep_{\Pi_{1}(U)}^{(n)}(\mathfrak{o}) \ar[d]^{\bf{P}}\\
        \mathscr{B}_{\mathfrak{X}_{\mathfrak{o}},D}(\Gl_{n}) \ar[r]^{\rho^{\mathfrak{X}}} & \Rep_{\Pi_{1}(U)}(\Gl_{n}(\mathfrak{o}))}
           \end{equation*} von Kategorien und Funktoren bis auf kanonische Isomorphismen von Funktoren kommutativ.\par
Insbesondere gilt $\rho_{{\bf{P}}(E)}={\bf{P}}(\rho_{E})$ für jedes Vektorbündel $E$ vom Rang $n$ in $\mathscr{B}_{\mathfrak{X}_{\mathfrak{o}},D}$. Damit sind für alle $x\in U(\mathbb{C}_{p})$ die Darstellungen von $\Pi_{1}(U)$ auf ${\bf{P}}(E)_{x}=\Iso_{\mathfrak{o}}(E_{x_{\mathfrak{o}}},\mathfrak{o}^{n})$ vermittels $\rho_{{\bf{P}}(E)}$ und vermittels $\gamma\mapsto(\rho_{E}(\gamma)^{-1})^{*}$ gleich.
\end{sat}
\bewub{Es sei $E\in\mathscr{B}_{\mathfrak{X}_{\mathfrak{o}},D}^{(n)}$. Dann gilt \[({\bf{P}}\circ\rho)(E)(x)={\bf{P}}(\rho_{E})(x)=\Iso_{\mathfrak{o}}(\rho_{E}(x),\mathfrak{o}^{n})=\Iso_{\mathfrak{o}}(E_{x_{\mathfrak{o}}},\mathfrak{o}^{n})\] sowie \[(\rho\circ{\bf{P}})(E)(x)=\rho_{{\bf{P}}(E)}(x)=({\bf{P}}(E))_{x_{\mathfrak{o}}}=(\Iso_{\mathfrak{X}_{\mathfrak{o}}}(E,\mathbb{A}_{\mathfrak{X}_{\mathfrak{o}}}^{n}))_{x_{\mathfrak{o}}}\] jeweils für alle $x\in U(\mathbb{C}_{p})$. Da per Definition die Bildung von $\Iso(\_,\_)$ mit der Bildung der globalen Schnitte vertauscht und wie oben gesehen der Funktor $\bf{P}$ verträglich mit Pullback ist, stimmen damit $({\bf{P}}\circ\rho)(E)(x)$ und $(\rho\circ{\bf{P}})(E)(x)$ für alle $x\in U(\mathbb{C}_{p})$ miteinander überein.\par
Ist nun $\gamma$ ein étaler Weg von $x$ nach $x'$ und \mbox{$\varphi\in\Iso_{\mathfrak{o}}(E_{x_{\mathfrak{o}}},\mathfrak{o}^{n})$}, so gilt einerseits \[{\bf{P}}(\rho_{E})(\gamma)(\varphi)=\varphi\circ(\rho_{E}(\gamma)^{-1})=\varphi\circ\projlim_{m}\rho_{E,m}(\gamma)^{-1}=\projlim_{m}(\varphi\circ(\rho_{E,m}(\gamma))^{-1})\] als Elemente von $\Iso_{\mathfrak{o}}(E_{x'_{\mathfrak{o}}},\mathfrak{o}_{n})$
und andererseits \[\rho_{{\bf{P}}(E)}(\gamma)(\varphi)=(\projlim_{m}\rho_{{\bf{P}}(E),m}(\gamma))(\varphi)=\projlim_{m}(\rho_{{\bf{P}}(E),m}(\gamma))(\varphi),\] als Elemente von $(\Iso_{\mathfrak{X}_{\mathfrak{o}}}(E,\mathbb{A}^{n}_{\mathfrak{X}_{\mathfrak{o}}}))_{x'_{\mathfrak{o}}}$.\par
Identifiziert man daher wie oben $\Iso_{\mathfrak{o}}(E_{x_{\mathfrak{o}}},\mathfrak{o}^{n})$ mit $x_{\mathfrak{o}}^{*}\Iso_{\mathfrak{o}}(E,\mathfrak{o}^{n})$ sowie $\Iso_{\mathfrak{o}_{m}}(E_{x_{m}},\mathfrak{o}_{m}^{n})$ mit $x_{m}^{*}\Iso_{\mathfrak{o}}(E,\mathfrak{o}^{n})$, so genügt es daher,\[(\rho_{{\bf{P}}(E),m}(\gamma))(\varphi)= \varphi\circ\rho_{E,m}(\gamma)\] zu zeigen. Dies folgt aber aus der konkreten Definition von $\rho_{E,m}(\gamma)$ als \[\rho_{E,m}(\gamma)=(\gamma y)_{m}^{*}\circ(y_{m}^{*})^{-1}\colon E_{x_{m}}\to E_{x'_{m}}\] und der von $\rho_{{\bf{P}}(E),m}(\gamma)$ als \[\rho_{{\bf{P}}(E),m}(\gamma)=(\gamma y)_{m}^{*}\circ(y_{m}^{*})^{-1}\colon(\Iso_{\mathfrak{X}_{\mathfrak{o}_{m}}}(E,\mathbb{A}_{\mathfrak{X}_{\mathfrak{o}_{m}}}))_{x_{m}}\to(\Iso_{\mathfrak{X}_{\mathfrak{o}_{m}}}(E,\mathbb{A}_{\mathfrak{X}_{\mathfrak{o}_{m}}}))_{x'_{m}},\] wobei aus letzterer mit der obigen Identifikation von $\Iso_{\mathfrak{o}_{m}}(E_{x_{m}},\mathfrak{o}_{m}^{n})$ mit \linebreak $x_{m}^{*}\Iso_{\mathfrak{o}}(E,\mathfrak{o}^{n})$ folgt, dass $\rho_{{\bf{P}}(E),m}(\gamma)(\varphi)=\varphi\circ(\gamma y)_{n}^{*}\circ(y_{n}^{*})^{-1}$ als Abbildung von $\Iso_{\mathfrak{o}_{m}}(E_{x_{m}},\mathfrak{o}_{m}^{n})$ nach $\Iso_{\mathfrak{o}_{m}}(E_{x'_{m}},\mathfrak{o}_{m}^{n})$  gilt, woraus die Behauptung folgt.\par 
Ist $\alpha\colon E\to E'$ ein Morphismus in $\mathscr{B}_{\mathfrak{X}_{\mathfrak{o}},D}^{(n)}$, so zeigt die obige Rechnung, dass $\rho_{{\bf{P}}(E)}(\alpha)={\bf{P}}(\rho_{E}(\alpha))$ gilt. \hfill$\Box$}\\\par

Aufgrund der Definition der Kategorie $\mathscr{B}_{X_{\mathbb{C}_{p}}, D}(G)$ über die Kategorie $\mathscr{B}_{\mathfrak{X}_{\mathfrak{o}},D}(G)$ und aufgrund dessen, dass der Funktor \[\mathscr{B}_{X_{\mathbb{C}_{p}}, D}(G)\to\Rep_{\Pi_{1}(U)}(G(\mathbb{C}_{p}))\] in gleicher Weise auf Basis des Funktors \[\mathscr{B}_{X_{\mathfrak{o}}, D}(G)\to\Rep_{\Pi_{1}(U)}(G(\mathfrak{o}))\] in §2.2 konstruiert wurde, übertragen sich alle Eigenschaften der Kategorie $\mathscr{B}_{\mathfrak{X}_{\mathfrak{o}},D}(G)$ und der betrachteten Funktoren in der Situation über $\mathfrak{o}$ auf die Situation über $\mathbb{C}_{p}$:\\\par
Insbesondere induziert der für jeden $\mathbb{Q}_{p}$"=Automorphismus $\sigma$ von $\overline{\mathbb{Q}_{p}}$ durch die Galoiskonjugation gegebene Funktor \[\sigma_{*}\colon\mathscr{B}_{\mathfrak{X}_{\mathfrak{o}},D}(G)\to \mathscr{B}_{{}^{\sigma}\mathfrak{X}_{\mathfrak{o}},D}({}^{\sigma}G)\] einen wohldefinierten Funktor \[\sigma_{*}\colon\mathscr{B}_{X_{\mathbb{C}_{p}}, D}(G)\to \mathscr{B}_{{}^{\sigma}X_{\mathbb{C}_{p}}, D}({}^{\sigma}G).\] Ebenso induziert der Funktor \[C_{\sigma}\colon\Rep_{\Pi_{1}(U)}(G(\mathfrak{o}))\to\Rep_{\Pi_{1}({}^{\sigma}U)}({}^{\sigma}G(\mathfrak{o}))\] einen Funktor \[C_{\sigma}\colon\Rep_{\Pi_{1}(U)}(G(\mathbb{C}_{p}))\to\Rep_{\Pi_{1}({}^{\sigma}U)}({}^{\sigma}G(\mathbb{C}_{p})).\] Für die so induzierten Funktoren $\sigma_{*}$ und $C_{\sigma}$ folgt direkt aus Satz \ref{galoiskonj} und Bemerkung \ref{bemgaloisconj}:\\\par 

\begin{sat}\label{galoiskonj1}
In der gegebenen Situation ist das Diagramm \begin{equation*}
        \xymatrix@=3em{%
        \mathscr{B}_{X_{\mathbb{C}_{p}}, D}(G) \ar[r]^{\rho} \ar[d]_{\sigma_{*}} & \Rep_{\Pi_{1}(U)}(G(\mathbb{C}_{p})) \ar[d]^{C_{\sigma}}\\
        \mathscr{B}_{{}^{\sigma}X_{\mathbb{C}_{p}}, {}^{\sigma}D}({}^{\sigma}G) \ar[r]^{\rho} & \Rep_{\Pi_{1}({}^{\sigma}U)}({}^{\sigma}G(\mathbb{C}_{p}))}
           \end{equation*} von Kategorien und Funktoren bis auf kanonische Isomorphismen von Funktoren kommutativ.\par
Insbesondere gilt $\rho_{{}^{\sigma}P}=\sigma_{*}\circ\rho_{P}\circ\sigma_{*}^{-1}$ als Funktoren von $\Pi_{1}({}^{\sigma}U)$ nach $\mathcal{P}({}^{\sigma}G(\mathbb{C}_{p}))$ für jeden $G$"=Torseur $P\in\mathscr{B}_{X_{\mathbb{C}_{p}}, D}(G)$.\par
Ist zusätzlich $X=X_{K}\otimes_{K}\overline{\mathbb{Q}_{p}}$ und $D=D_{K}\otimes_{K}\overline{\mathbb{Q}_{p}}$ sowie $G=G_{K}\otimes_{\mathfrak{o}_{K}}\mathfrak{o}$ für einen gewissen Körper $\mathbb{Q}_{p}\subset K\subset\overline{\mathbb{Q}_{p}}$, so lassen sich ${}^{\sigma}X$, ${}^{\sigma}D$ und ${}^{\sigma}G$ für jedes $\sigma\in\Gal(\overline{\mathbb{Q}_{p}},K)$ kanonisch mit $X$,$D$ und $G$ über $\overline{\mathbb{Q}_{p}}$ bzw. $\mathfrak{o}$ identifizieren und der Funktor \[\rho\colon\mathscr{B}_{X_{\mathbb{C}_{p}}, D}(G)\to\Rep_{\Pi_{1}(U)}(G(\mathbb{C}_{p}))\] kommutiert mit den Aktionen von $\Gal(\overline{\mathbb{Q}_{p}}/K)$ auf diesen Kategorien via $\sigma_{*}$ bzw. $C_{\sigma}$ für $\sigma\in\Gal(\overline{\mathbb{Q}_{p}},K)$. \\
\end{sat}\par

Aus Satz \ref{unabhaengig} und \cite[Lemma 8]{DW1} folgt zudem, dass die Konstruktion des Funktors \[\rho\colon\mathscr{B}_{X_{\mathbb{C}_{p}}, D}(G)\to\Rep_{\Pi_{1}(U)}(G(\mathbb{C}_{p}))\] mit Pullback entlang von Morphismen $f\colon X\to X'$ glatter und projektiver Kurven über $\overline{\mathbb{Q}_{p}}$ vertauscht:\\\par

\begin{sat}\label{unabh}
Es seien $X$ und $X'$ zwei glatte und projektive Kurven über $\overline{\mathbb{Q}_{p}}$ und $f\colon X\to X'$ ein Morphismus zwischen ihnen. Außerdem sei $D'$ ein Divisor auf $X'$ und $G$ ein glattes und affines Gruppenschema von endlicher Präsentation über $\mathfrak{o}$. Dann induziert das Pullback von $G$"=Torseuren entlang $f$ einen Funktor $f^{*}\colon\mathscr{B}_{X'_{\mathbb{C}_{p}}, D'}(G)\to\mathscr{B}_{X_{\mathbb{C}_{p}}, f^{*}D'}(G)$ und das folgende Diagramm von Kategorien und Funktoren kommutiert (bis auf kanonische Isomorphismen):\par 
\begin{equation*}
        \xymatrix@=3em{%
        \mathscr{B}_{X_{\mathbb{C}_{p}}, D'}(G) \ar[rr]^{\rho} \ar[d]^{f^{*}}& & \Rep_{\Pi_{1}(U')}(G(\mathbb{C}_{p})) \ar[d]^{A(f)} \\
  \mathscr{B}_{X_{\mathbb{C}_{p}}, f^{*}D}(G) \ar[rr]^{\rho} & & \Rep_{\Pi_{1}(U)}(G(\mathbb{C}_{p})). \\}\end{equation*} 
  Dabei sei $A(f)\colon\Rep_{\Pi_{1}(U')}(G(\mathbb{C}_{p}))\to\Rep_{\Pi_{1}(U)}(G(\mathbb{C}_{p}))$ der durch den Funktor $A(f)\colon\Rep_{\Pi_{1}(U')}(G(\mathfrak{o}))\to\Rep_{\Pi_{1}(U)}(G(\mathfrak{o}))$ induzierte Funktor.\par
           Insbesondere gilt für alle $\widetilde{P}\in\mathscr{B}_{X'_{\mathbb{C}_{p}}, D'}(G)$, dass $\rho_{f^{*}\widetilde{P}}=\rho_{\widetilde{P}}\circ f_{*}$ als Funktoren von $\Pi_{1}(U)$ nach $\mathcal{P}(G(\mathbb{C}_{p}))$ gilt.\\
\end{sat}\par

Für einen Morphismus $\varphi\colon G\to G'$ affiner, glatter Gruppenschemata von endlicher Präsentation über $\mathfrak{o}$ induzieren ferner die oben definierten Funktoren \[\varphi_{*}\colon\mathscr{B}_{\mathfrak{X}_{\mathfrak{o}},D}(G)\to\mathscr{B}_{\mathfrak{X}_{\mathfrak{o}},D}(G')\] und \[\varphi_{*}\colon\Rep_{\Pi_{1}(U)}(G(\mathfrak{o}))\to\Rep_{\Pi_{1}(U)}(G'(\mathfrak{o}))\] Funktoren \[\varphi_{*}\colon\mathscr{B}_{X_{\mathbb{C}_{p}}, D}(G)\to\mathscr{B}_{X_{\mathbb{C}_{p}}, D}(G')\] und \[\varphi_{*}\colon\Rep_{\Pi_{1}(U)}(G(\mathbb{C}_{p}))\to\Rep_{\Pi_{1}(U)}(G'(\mathbb{C}_{p})).\]\\\par
Für die so definierten Funktoren $\varphi_{*}$ folgt dann direkt aus Satz \ref{gruppenwechsel}:\\\par

\begin{sat}\label{gruppen1}
In der gegebenen Situation ist das Diagramm \begin{equation*}
        \xymatrix@=3em{%
        \mathscr{B}_{X_{\mathbb{C}_{p}}, D}(G) \ar[r]^{\rho} \ar[d]_{\varphi_{*}} & \Rep_{\Pi_{1}(U)}(G(\mathbb{C}_{p})) \ar[d]^{\varphi_{*}}\\
        \mathscr{B}_{X_{\mathbb{C}_{p}}, D}(G') \ar[r]^{\rho} & \Rep_{\Pi_{1}(U)}(G'(\mathbb{C}_{p}))}
           \end{equation*} von Kategorien und Funktoren bis auf kanonische Isomorphismen von Funktoren kommutativ.\par
Insbesondere gilt $\rho_{\varphi_{*}P}=\varphi_{*}\circ\rho_{P}$ als Funktoren von $\Pi_{1}(U)$ nach $\mathcal{P}(G'(\mathbb{C}_{p}))$ für jeden $G$"=Torseur $P\in\mathscr{B}_{X_{\mathbb{C}_{p}}, D}(G)$.\\
\end{sat}\par

Als nächstes halten wir fest, dass für jede algebraische Darstellung ${\bf{\Gamma}}$ des Gruppenschemas $G$ der oben definierte Funktor \[\_\times^{G}{\bf{\Gamma}}\colon\mathscr{B}_{\mathfrak{X}_{\mathfrak{o}},D}(G)\to\mathscr{B}_{\mathfrak{X}_{\mathfrak{o}},D}\] einen Funktor \[\_\times^{G}{\bf{\Gamma}}\colon\mathscr{B}_{X_{\mathbb{C}_{p}}, D}(G)\to\mathscr{B}_{X_{\mathbb{C}_{p}}, D}\] induziert. Ebenso erhält man aus den Bifunktoren  \[\_\times^{G}\_\colon\mathscr{B}_{\mathfrak{X}_{\mathfrak{o}},D}(G)\times\Rep_{G}(\mathfrak{o})\to\mathscr{B}_{\mathfrak{X}_{\mathfrak{o}},D}\] und \[\_\times^{G(\mathfrak{o})}\_\colon\Rep_{\Pi_{1}(U)}(G(\mathfrak{o}))\times\Rep_{G(\mathfrak{o})}(\mathfrak{o})\to\Rep_{\Pi_{1}(U)}(\mathfrak{o})\] Bifunktoren \[\_\times^{G}\_\colon\mathscr{B}_{X_{\mathbb{C}_{p}}, D}(G)\times\Rep_{G}(\mathfrak{o})\to\mathscr{B}_{X_{\mathbb{C}_{p}}, D}\] und \[\_\times^{G(\mathbb{C}_{p})}\_\colon\Rep_{\Pi_{1}(U)}(G(\mathbb{C}_{p}))\times\Rep_{G(\mathbb{C}_{p})}(\mathbb{C}_{p})\to\Rep_{\Pi_{1}(U)}(\mathbb{C}_{p}),\] wobei $\Rep_{G(\mathbb{C}_{p})}(\mathbb{C}_{p})$ die Kategorie der stetigen Darstellungen von $G(\mathbb{C}_{p})$ auf endlichdimensionalen $\mathbb{C}_{p}$"=Vektorräumen sei.
Mit diesen Funktoren ergibt sich dann aus Satz \ref{darstellung}:\\\par

\begin{sat}\label{darstellung1}
In der obigen Situation kommutiert das folgende Diagramm:\par
\begin{equation*}
        \xymatrix@=3em{%
        \mathscr{B}_{X_{\mathbb{C}_{p}},D}(G)\times\Rep_{G}(\mathfrak{o}) \ar[rrr]^{\_\times^{G}\_}\ar[d]_{\rho\times\text{nat}} &&& \mathscr{B}_{X_{\mathbb{C}},D} \ar[d]^{\rho}\\
        \Rep_{\Pi_{1}(U)}(G(\mathbb{C}_{p}))\times\Rep_{G(\mathbb{C}_{p})}(\mathbb{C}_{p})   \ar[rrr]^{\_\times^{G(\mathbb{C}_{p})}\_} & &&\Rep_{\Pi_{1}(U)}(\mathbb{C}_{p}).}
           \end{equation*}
           Insbesondere gilt \[\rho_{P\times^{G}\bf{\Gamma}}=\rho_{P}\times^{G(\mathbb{C}_{p})}(\Gamma\otimes\mathbb{C}_{p})\] für alle $P\in\mathscr{B}_{X_{\mathbb{C}_{p}},D}(G)$ und $\Gamma\in\Rep_{G}(\mathfrak{o})$ als Funktoren von $\Pi_{1}(U)$ nach $\mathcal{P}(G(\mathbb{C}_{p}))$.\\
\end{sat}\par

Schließlich induziert auch der Funktor \[{\bf{P}}\colon\mathscr{B}_{\mathfrak{X}_{\mathfrak{o}},D}^{(n)}\to\mathscr{B}_{\mathfrak{X}_{\mathfrak{o}},D}(\Gl_{n})\] einen Funktor \[{\bf{P}}\colon\mathscr{B}_{X_{\mathbb{C}_{p}},D}^{(n)}\to\mathscr{B}_{X_{\mathbb{C}_{p}},D}(\Gl_{n}).\] Ebenso induziert der Funktor \[{\bf{P}}\colon\Rep_{\Pi_{1}(U)}^{(n)}(\mathfrak{o})\to\Rep_{\Pi_{1}(U)}(\Gl_{n}(\mathfrak{o}))\] einen Funktor \[{\bf{P}}\colon\Rep_{\Pi_{1}(U)}^{(n)}(\mathbb{C}_{p})\to\Rep_{\Pi_{1}(U)}(\Gl_{n}(\mathbb{C}_{p})).\]
Mit den so definierten Funktoren $\bf{P}$ folgt nun aus Satz \ref{pblvbl}:\\\par

\begin{sat}\label{vb1}
In der gegebenen Situation ist das Diagramm \begin{equation*}
        \xymatrix@=3em{%
        \mathscr{B}_{X_{\mathbb{C}_{P}},D}^{(n)} \ar[r]^{\rho} \ar[d]_{\bf{P}} & \Rep_{\Pi_{1}(U)}^{(n)}(\mathbb{C}_{p}) \ar[d]^{\bf{P}}\\
        \mathscr{B}_{X_{\mathbb{C}_{p}},D}(\Gl_{n}) \ar[r]^{\rho} & \Rep_{\Pi_{1}(U)}(\Gl_{n}(\mathbb{C}_{p}))}
           \end{equation*} von Kategorien und Funktoren bis auf kanonische Isomorphismen von Funktoren kommutativ.\par
Insbesondere gilt $\rho_{{\bf{P}}(E)}={\bf{P}}(\rho_{E})$ für jedes Vektorbündel $E$ vom Rang $n$ in $\mathscr{B}_{X_{\mathbb{C}_{p}},D}$. Damit sind für alle $x\in U(\mathbb{C}_{p})$ die Darstellungen von $\Pi_{1}(U)$ auf $({\bf{P}}(E))_{x}=\Iso_{\mathbb{C}_{p}}(E_{x},\mathbb{C}_{p}^{n})$ vermittels $\rho_{{\bf{P}}(E)}$ und vermittels $\gamma\mapsto(\rho_{E}(\gamma)^{-1})^{*}$ gleich.\\\par
\end{sat}

\chapter{Charakterisierung der Kategorie $\mathscr{B}_{\mathfrak{X}_{\mathfrak{o}},D}(G)$}

Anhand der im vorherigen Kapitel gegebenen Definition der Kategorie $\mathscr{B}_{\mathfrak{X}_{\mathfrak{o}},D}(G)$ lässt sich nur sehr schwer nachprüfen, ob ein gegebener $G$"=Torseur $P$ auf $\mathfrak{X}_{\mathfrak{o}}$ zu dieser Kategorie gehört. Im folgenden soll daher eine leichter nachprüfbare Charakterisierung der Objekte von $\mathscr{B}_{\mathfrak{X}_{\mathfrak{o}},D}(G)$ entwickelt werden, die analoge Resultate zu den Theoremen 16, 17 und 18 in \cite{DW1} liefert.\par

\section{Reduktion auf die spezielle Faser}

Als erstes zeigen wir:\\\par

\begin{thm}\label{thm16}
Es sei $\mathfrak{X}$ ein Modell über $\overline{\mathbb{Z}}_{p}$ der glatten und projektiven Kurve $X$ über $\overline{\mathbb{Q}}_{p}$. Es sei $k=\overline{\mathbb{F}_{p}}$ der Restklassenkörper von $\overline{\mathbb{Z}}_{p}$.\par
Dann liegt ein $G$"=Torseur $P$ auf $\mathfrak{X}_{\mathfrak{o}}$ genau dann in $\mathscr{B}_{\mathfrak{X}_{\mathfrak{o}},D}(G)$, wenn es ein Objekt $\pi\colon\mathcal{Y}\to\mathfrak{X}$ aus der Kategorie $S_{\mathfrak{X},D}$ gibt, so dass $\pi_{k}^{*}P_{k}$ trivial auf \mbox{$\mathcal{Y}_{k}=\mathcal{Y}\otimes_ {\overline{\mathbb{Z}}_{p}}k$} ist.
\end{thm}
\bewub{\begin{itemize}{\item[(i)]Dass die Bedingung notwendig ist, ist klar.\par
Es sei also $P$ ein $G$"=Torseur auf $\mathfrak{X}_{\mathfrak{o}}$, für den es ein Objekt $\pi\colon\mathcal{Y}\to\mathfrak{X}$ aus der Kategorie $S_{\mathfrak{X},D}$ gibt, so dass $\pi_{k}^{*}P_{k}$ trivial auf \mbox{$\mathcal{Y}_{k}=\mathcal{Y}\otimes_ {\overline{\mathbb{Z}}_{p}}k$} ist. Nach Korollar \ref{doms} darf angenommen werden, dass $\pi\in S_{\mathfrak{X},D}^{ss}$ ist.\par
\item[(ii)]
Als nächster Beweisschritt soll nun gezeigt werden, dass die Familie \[(\mathfrak{X}, D, G, P_{1},\pi\colon\mathcal{Y}\to\mathfrak{X})\] zu einer Familie \[(\mathfrak{X}^{0},D^{0},G^{0}, \overline{P}, \pi_{0}\colon\mathcal{Y}^{0}\to\mathfrak{X}^{0})\] über dem Ganzheitsring $\mathfrak{0}_{K}$ einer endlichen Erweiterung $K$ von $\mathbb{Q}_{p}$ absteigt.\par
Dabei ist $\mathfrak{X}^{0}$ ein Modell über $\mathfrak{o}_{K}$ der glatten und projektiven Kurve \mbox{$X^{0}=\mathfrak{X}^{0}\otimes_{\mathfrak{o}_{K}}K$} über $K$, $G^{0}$ ein Gruppenschema, das affin, glatt und von endlicher Präsentation über einer normalen, endlich erzeugten $\mathfrak{o}_{K}$"=Algebra $A$ ist, und $\overline{P}$ ein $G^{0}_{\mathfrak{X}^{0}_{A}\otimes_{A}\mathfrak{o}_{K}/p\mathfrak{o}_{K}}$"=Torseur auf $\mathfrak{X}^{0}_{A}\otimes_{A}\mathfrak{o}_{K}/p\mathfrak{o}_{K}$, dessen Einschränkung auf die spezielle Faser \mbox{$\mathfrak{X}^{0}\otimes_{\mathfrak{o}_{K}}\mathfrak{o}_{K}/\mathfrak{p}$} nach Pullback entlang $\pi_{0}\otimes\mathfrak{o}_{K}/\mathfrak{p}$ trivial ist. Außerdem ist $\pi_{0}$ ein Element von $S_{\mathfrak{X}_{0},D_{0}}^{ss}$.\par
Dies sieht man wie folgt:\par
Unter den gegebenen Voraussetzungen existieren nach Korollar \ref{kss} eine endliche Erweiterung $K^{1}$ von $\mathbb{Q}_{p}$ und eine glatte und projektive Kurve $X^{1}$ über $K^{1}$ mit Modell $\mathfrak{X}^{1}$ über $\mathfrak{o}_{K^{1}}$ sowie ein Divisor $D^{1}$ von $X^{1}$, so dass zum einen \mbox{$X=X^{1}\otimes_{K^{1}}\overline{K^{1}}$}, \mbox{$D=D^{1}\otimes_{K^{1}}\overline{K^{1}}$} und \mbox{$\mathfrak{X}=\mathfrak{X}^{1}\otimes_{\mathfrak{o}_K^{1}}\overline{\mathbb{Z}_{p}}$} sind, zum anderen ein Element \mbox{$\pi^{1}\colon\mathcal{Y}^{1}\rightarrow\mathfrak{X}^{1}$} aus $S_{\mathfrak{X}^{1},D^{1}}^{ss}$ mit der Eigenschaft existiert, dass \mbox{$\widetilde{\pi}^{1}=\pi^{1}\otimes_{\mathfrak{o}_{K^{1}}}\overline{\mathbb{Z}_{p}}$} den gegebenen Morphismus $\pi$ dominiert.\par
Da nach Voraussetzung $G$ von endlicher Präsentation über $\mathfrak{o}$ ist, existiert nach \cite[Proposition 8.9.1]{EGA IV} eine endliche Erweiterung $K^{2}$ von $\mathbb{Q}_{p}$, so dass $G$ zu einem Schema $G^{2}$ von endlichem Typ über einer normalen endlich erzeugten $\mathfrak{o}_{K^{2}}$"=Algebra $A^{2}$ absteigt. Da $\mathfrak{o}_{K^{2}}$ und damit auch $A^{2}$ noethersch ist, ist nach \cite[Définition 1.6.1]{EGA IV} $G^{2}$ von endlicher Präsentation über $A^{2}$. Nach \cite[Scholie 8.8.3]{EGA IV} ist $G^{2}$ zudem ein Gruppenschema.\par
Nach \cite[Théorème 8.10.5 (viii)]{EGA IV} und \cite[Pro\-po\-si\-tion 17.7.8]{EGA IV} existiert dann eine endliche Erweiterung $K^{3}$ von $\mathbb{Q}_{p}$ mit $K^{2}\subset K^{3}$, so dass das Gruppenschema $G^{3}:=G^{2}\otimes_{A^{2}}A^{3}$ affin, glatt und von endlicher Präsentation über $A^{3}$ ist, wobei $A^{3}:=A^{2}\otimes_{\mathfrak{o}_{K^{2}}}\mathfrak{o}_{K^{3}}$ ist.\par
Weil die endlichen Körpererweiterungen von $\mathbb{Q}_{p}$ ein induktives System bilden, gibt es nun eine endliche Erweiterung $K^{4}$ von $\mathbb{Q}_{p}$, die $K^{1}$ und $K^{3}$ als Zwischenerweiterungen enthält. Dann ist $X^{4}:=X^{1}\otimes_{K^{1}}K^{4}$ eine glatte und projektive Kurve über $K^{4}$, $\mathfrak{X}^{4}:=\mathfrak{X}^{1}\otimes_{\mathfrak{o}_{K^{1}}}\mathfrak{o}_{K^{4}}$ ein Modell von $X^{4}$ über $\mathfrak{o}_{K^{4}}$, $D^{4}:=D^{1}\otimes_{K^{1}}K^{4}$ ein Divisor von $K^{4}$, $G^{4}:=G^{3}\otimes_{A^{3}}A^{4}$ ein affines, glattes Gruppenschema von endlicher Präsentation über der normalen endlich erzeugten $\mathfrak{o}_{K^{4}}$"=Algebra $A^{4}:=A^{3}\otimes_{\mathfrak{o}_{K^{3}}}\mathfrak{o}_{K^{4}}$ und der Morphismus $\widetilde{\pi}^{4}$ dominiert den Morphismus $\pi$.\par  
Da $P_{1}$ insbesondere von endlicher Präsentation über $\mathfrak{X}_{\mathfrak{o}}\otimes_{\mathfrak{o}}\mathfrak{o}/p\mathfrak{o}$ ist, existiert nach \cite[Proposition 8.9.1]{EGA IV} eine endliche Körpererweiterung $K^{5}$ von $\mathbb{Q}_{p}$, die $K^{4}$ enthält, so dass $P_{1}$ zu einem Schema $P^{5}$ von endlichem Typ über $\mathfrak{X}^{5}_{A^{5}}\otimes_{A^{5}}\mathfrak{o}_{K^{5}}/p\mathfrak{o}_{K^{5}}$ absteigt. Da $\mathfrak{X}^{5}_{A^{5}}\otimes_{A^{5}}\mathfrak{o}_{K^{5}}/p\mathfrak{o}_{K^{5}}$ lokal noethersch ist, ist nach \cite[Définition 1.6.1]{EGA IV} $P^{5}$ auch von endlicher Präsentation über $\mathfrak{X}^{5}_{A^{5}}\otimes_{A^{5}}\mathfrak{o}_{K^{5}}/p\mathfrak{o}_{K^{5}}$.\par
Nach \cite[Théorème 8.10.5 (vi)]{EGA IV} und \cite[Proposition 11.2.6]{EGA IV} existiert dann eine endliche Erweiterung $K^{6}$ von $\mathbb{Q}_{p}$, die $K^{5}$ enthält, so dass $P^{6}:=P^{5}\times_{\mathfrak{X}^{5}_{A^{5}}\otimes_{A^{5}}\mathfrak{o}_{K^{5}}/p\mathfrak{o}_{K^{5}}}\Spec(\mathfrak{o}_{K^{6}}/p\mathfrak{o}_{K^{6}})$ treuflach und von endlicher Präsentation über $\mathfrak{X}^{6}_{A^{6}}\otimes_{A^{6}}\mathfrak{o}_{K^{6}}/p\mathfrak{o}_{K^{6}}$ ist. Es seien wieder $X^{6}, D^{6}, \mathfrak{X}^{6}$ etc. die entsprechenden Basiswechsel. Nach \cite[Proposition 8.9.1]{EGA IV} steigt der Morphismus \[f\colon P_{1}\times_{\mathfrak{X}_{\mathfrak{o}}\otimes_{\mathfrak{o}\textbf{}}\mathfrak{o}/p\mathfrak{o}}G_{\mathfrak{X}_{\mathfrak{o}}\otimes_{\mathfrak{o}}\mathfrak{o}/p\mathfrak{o}}\to P_{1},\] mittels dessen $G_{\mathfrak{X}_{\mathfrak{o}}\otimes_{\mathfrak{o}}\mathfrak{o}/p\mathfrak{o}}$ auf $P_{1}$ operiert, zu einem Morphismus \[f^{6}\colon P^{6}\times_{\mathfrak{X}^{6}_{A^{6}}\otimes_{A^{6}}\mathfrak{o}_{K^{6}}/p\mathfrak{o}_{K^{6}}}G^{6}_{\mathfrak{X}^{6}_{A^{6}}\otimes_{A^{6}}\mathfrak{o}_{K^{6}}/p\mathfrak{o}_{K^{6}}}\to P^{6}\] ab. Man beachte dabei, dass noetherscher Descent, wie bereits oben benutzt, mit der Bildung von Faserprodukten vertauscht.\par
Ebenso existiert nach \cite[Théorème 8.10.5 (i)]{EGA IV} eine endliche Erweiterung $K^{7}$ von $\mathbb{Q}_{p}$ mit $K^{6}$ als Zwischenkörper, so dass der Isomorphismus \[(\id_{P_{1}},f)\colon P_{1}\times_{\mathfrak{X}_{\mathfrak{o}}\otimes_{\mathfrak{o}}\mathfrak{o}/p\mathfrak{o}}G_{\mathfrak{X}_{\mathfrak{o}}\otimes_{\mathfrak{o}}\mathfrak{o}/p\mathfrak{o}}\longrightarrow P_{1}\times_{\mathfrak{X}_{o}\otimes_{\mathfrak{o}}\mathfrak{o}/p\mathfrak{o}} P_{1}\] zu \[(\id_{P^{7}},f^{7})\colon P^{7}\times_{\mathfrak{X}^{7}_{A^{7}}\otimes_{A^{7}}\mathfrak{o}_{K^{7}}/p\mathfrak{o}_{K^{7}}}G_{\mathfrak{X}^{7}_{A^{7}}\otimes_{A^{7}}\mathfrak{o}_{K^{7}}/p\mathfrak{o}_{K^{7}}}\longrightarrow P^{7}\times_{\mathfrak{X}^{7}_{A^{7}}\otimes_{A^{7}}\mathfrak{o}_{K^{7}}/p\mathfrak{o}_{K^{7}}} P^{7}\] absteigt und dies ein Isomorphismus ist.\par
Somit ist $P^{7}$ ein $G^{7}_{\mathfrak{X}^{7}_{A^{7}}\otimes_{A^{7}}\mathfrak{o}_{K^{7}}/p\mathfrak{o}_{K^{7}}}$"=Torseur auf $\mathfrak{X}^{7}_{A^{7}}\otimes_{A^{7}}\mathfrak{o}_{K^{7}}/p\mathfrak{o}_{K^{7}}$.\par
Nach \cite[Théorème 8.10.5 (i)]{EGA IV} existiert schließlich eine endliche Erweiterung $K^{8}$ von $\mathbb{Q}_{p}$ mit $K^{7}\subset K^{8}$, so dass die Einschränkung von $P^{8}$ auf die spezielle Faser $\mathfrak{X}^{8}_{A^{8}}\otimes_{A^{8}}\mathfrak{o}_{K^{8}}/\mathfrak{p}$ nach Pullback entlang $\pi^{8}\otimes\mathfrak{o}_{K^{8}}/\mathfrak{p}$ trivial ist:\par
Denn dass nach Voraussetzung die Einschränkung von $P_{1}$ auf die spezielle Faser $\mathfrak{X}_{\mathfrak{o}}\otimes\mathfrak{o}/\mathfrak{p}$ nach Pullback entlang $(\pi\otimes\id_{\mathfrak{o}})\otimes\id_{\mathfrak{o}/\mathfrak{p}}$ trivial wird, bedeutet gerade, dass es einen Isomorphismus \[\psi\colon((\pi\otimes\id_{\mathfrak{o}})\otimes\id_{\mathfrak{o}/\mathfrak{p}})^{*}P_{1}\mathop{|}\nolimits_{\mathfrak{X}_{\mathfrak{o}}\otimes_{\mathfrak{o}}\mathfrak{o}/\mathfrak{p}}\stackrel{\sim}{\longrightarrow}G_{\mathcal{Y}_{\mathfrak{o}}\otimes_{\mathfrak{o}}\mathfrak{o}/\mathfrak{p}}\] von $G_{\mathcal{Y}_{\mathfrak{o}}\otimes_{\mathfrak{o}}\mathfrak{o}/\mathfrak{p}}$"=Torseuren gibt. Nach \cite[Proposition 8.9.1]{EGA IV} und \linebreak\cite[Théorème 8.10.5 (i)]{EGA IV} existiert daher eine endliche Erweiterung $K^{8}$ von $\overline{\mathbb{Q}}_{p}$, so dass $\psi$ zu einem Isomorphismus \[\psi^{8}\colon(\pi^{8}\otimes\id_{\mathfrak{o}_{K^{8}}/\mathfrak{p}})^{*}P^{8}\mathop{|}\nolimits_{\mathfrak{X}^{8}_{A^{8}}\otimes_{A^{8}}\mathfrak{o}_{K^{8}}/\mathfrak{p}}\stackrel{\sim}{\longrightarrow} G_{\mathcal{Y}^{8}_{A^{8}}\otimes_{A^{8}}\mathfrak{o}_{K^{8}}/\mathfrak{p}}\] von $\mathcal{Y}^{8}_{A^{8}}\otimes_{A^{8}}\mathfrak{o}_{K^{8}}/\mathfrak{p}$"=Schemata absteigt. Dieser ist auch ein Isomorphismus von $G_{\mathcal{Y}^{8}_{A^{8}}\otimes_{A^{8}}\mathfrak{o}_{K^{8}}/\mathfrak{p}}$"=Torseuren, da die Relation \[m_{G_{\mathcal{Y}_{\mathfrak{o}}\otimes_{\mathfrak{o}}\mathfrak{o}/\mathfrak{p}}}\circ(\psi\otimes\id_{G_{\mathcal{Y}_{\mathfrak{o}}\otimes_{\mathfrak{o}}\mathfrak{o}/\mathfrak{p}}})=\psi\circ f'\] unter noetherschem Descent erhalten bleibt, wobei $m_{G_{\mathcal{Y}_{\mathfrak{o}}\otimes_{\mathfrak{o}}\mathfrak{o}/\mathfrak{p}}}$ die Gruppenmultiplikation von $G_{\mathcal{Y}_{\mathfrak{o}}\otimes_{\mathfrak{o}}\mathfrak{o}/\mathfrak{p}}$ und $f'$ die Aktion von $G_{\mathcal{Y}_{\mathfrak{o}}\otimes_{\mathfrak{o}}\mathfrak{o}/\mathfrak{p}}$ auf \linebreak \mbox{$((\pi\otimes\id_{\mathfrak{o}})\otimes\id_{\mathfrak{o}/\mathfrak{p}})^{*}P_{1}\mathop{|}_{\mathfrak{X}_{\mathfrak{o}}\otimes_{\mathfrak{o}}\mathfrak{o}/\mathfrak{p}}$} bezeichnet. Man benutzt dabei wieder, dass noetherscher Descent verträglich mit Faserprodukten ist.\par
Setzt man also $K:=K^{8}$ und $\overline{P}:=P^{8}$, $\mathfrak{X}^{0}:=\mathfrak{X}^{8}$, $D^{0}:=D^{8}$, $G^{0}:=G^{8}$ und $\pi^{8}:=\pi_{0}$, so erfüllt die Familie \[(\mathfrak{X}^{0},D^{0},G^{0}, \overline{P}, \pi_{0}\colon\mathcal{Y}^{0}\to\mathfrak{X}^{0})\] alle geforderten Eigenschaften.\par
\item[(iii)] Es sei $e$ der Verzweigungsgrad von $K$ über $\mathbb{Q}_{p}$ und man setze \[\mathfrak{o}_{\nu/e}=\mathfrak{o}/\mathfrak{p}^{\nu}\mathfrak{o}=\overline{\mathbb{Z}}_{p}/\mathfrak{p}^{\nu}\overline{\mathbb{Z}}_{p}.\] Man bemerke, dass diese Notation mit der bisherigen Notation \mbox{$\mathfrak{o}_{n}=\mathfrak{o}/p^{n}\mathfrak{o}$} verträglich ist.\par
Es sei $\pi_{\nu/e}, \widetilde{P}_{\nu/e}$ etc. jeweils der Basiswechsel mit $\mathfrak{o}_{\nu/e}$, wobei $\widetilde{P}$ die generische Faser von $P$ sei. Weil $\pi_{1/e}$ auch der Basiswechsel von \mbox{$\pi_{0}\otimes_{\mathfrak{o}_{K}}\mathfrak{o}_{K}/\mathfrak{p}$} mit $\mathfrak{o}_{1/e}$ ist, folgt, dass $\pi_{1/e}^{*}P_{1/e}=(\pi_{0}\otimes\mathfrak{o}_{K}/\mathfrak{p}\otimes\mathfrak{o}_{1/e})^{*}(\overline{P}\otimes\mathfrak{o}/\mathfrak{p})$ trivial auf $\mathcal{Y}_{1/e}$ ist.\par
Mittels Induktion genügt es daher, die folgende Behauptung (*) zu zeigen:\par
Ist $\nu\geq 2$ und $\pi\colon\mathcal{Y}\to\mathfrak{X}$ ein Element von $S_{\mathfrak{X},D}^{ss}$, so dass $\pi^{*}_{(\nu-1)/e}P_{(\nu-1)/e}$ trivial ist, so existiert ein Objekt $\mu\colon\mathcal{Z}\to\mathfrak{X}$ in $S_{\mathfrak{X},D}^{ss}$, so dass $\mu^{*}_{\nu/e}P_{\nu/e}$ ein trivialer $G$"=Torseur auf $\mathcal{Z}_{\nu/e}$ ist. 
\item[iv)] Wir beweisen daher die Behauptung (*):\par Es sei also \mbox{$\nu\geq 2$} und \mbox{$\pi\colon\mathcal{Y}\longrightarrow\mathfrak{X}$} ein Element von $S_{\mathfrak{X},D}^{ss}$, so dass \linebreak\mbox{$\pi^{*}_{(\nu-1)/e}P_{(\nu-1)/e}$} trivial ist.\par
Man betrachte die abgeschlossene Immersion $i\colon\mathcal{Y}_{(\nu-1)/e}\lhook\longrightarrow\mathcal{Y}_{\nu/e}$.\par
Ist $\mathcal{Y}_{\nu/e}$ affin, so ist nach \cite[Corollaire VII.1.3.2]{G} die kanonische Abbildung $H^{1}(\mathcal{Y}_{\nu/e},G_{\mathcal{Y}_{\nu/e}})\to H^{1}(\mathcal{Y}_{(\nu-1)/e},G_{\mathcal{Y}_{(\nu-1)/e}})$ sowohl für die fppf"=Topologie als auch für die étale Topologie eine Bijektion und damit auch $\pi^{*}_{\nu/e}P_{\nu/e}$ trivial. Also kann man in diesem Fall $\mu=\pi$ wählen.\par
Im folgenden sei daher angenommen, dass $\mathcal{Y}_{\nu/e}$ nicht affin ist. Es sei $\mathcal{I}$ die Idealgarbe, die die abgeschlossene Immersion $i$ definiert und für die $\mathcal{I}^{2}=0$ gilt.
Dann hat man nach \cite[Lemme VII.1.3.5]{G} folgende exakte Sequenz von Garben auf $(\mathcal{Y}_{\nu/e})_{\acute{e}t}$: \begin{equation*}
        \xymatrix@=1em{%
         1 \ar[r] & i_{*}(\mathscr{L}_{(\nu-1)/e}\otimes_{\mathcal{O}_{\mathcal{Y}_{(\nu-1)/e}}}\mathcal{I}) \ar[rrrr]^-{h} & & & & G \ar[rrr]^-{adj}& & & i_{*}G_{\mathcal{Y}_{(\nu-1)/e}} \ar[r] & 1,}
           \end{equation*}
 wobei $\mathscr{L}_{(\nu-1)/e}:=Lie(G_{\mathcal{Y}_{(\nu-1)/e}}/\mathcal{Y}_{(\nu-1)/e})$ sei.        
Man beachte dabei, dass $i$ ein universeller Homöomorphismus ist, so dass nach \cite[Remark II.3.17]{M} $\widetilde{(\mathcal{Y}_{(\nu-1)/e})}_{\acute{e}t}\cong\widetilde{(\mathcal{Y}_{\nu/e})}_{\acute{e}t}$ gilt und man $\mathcal{I}$ kanonisch mit einem $\mathcal{O}_{\mathcal{Y}_{(\nu-1)/e}}$"=Modul identifizieren kann. Desweiteren bezeichnet $\adj$ die Komposition des Adjunktionsmorphismus $G_{\mathcal{Y}_{\nu/e}}\longrightarrow i_{*}i^{*}G_{\mathcal{Y}_{\nu/e}}$ mit dem Morphismus $i_{*}i^{*}G_{\mathcal{Y}_{\nu/e}}\longrightarrow i_{*}G_{\mathcal{Y}_{(\nu-1)/e}}$, der von dem kanonischen Morphismus \linebreak\mbox{$i^{*}G_{\mathcal{Y}_{\nu/e}}\longrightarrow G_{\mathcal{Y}_{(\nu-1)/e}}$} induziert wird.
Der Morphismus $h$ schließlich wird durch die Bijektion \[i_{*}(\mathscr{L}_{(\nu-1)/e}\otimes_{\mathcal{O}_{\mathcal{Y}_{(\nu-1)/e}}}\mathcal{I})(S)\to\ker(G_{\mathcal{Y}_{\nu/e}}\to i_{*}G_{\mathcal{Y}_{(\nu-1)/e}})(S),\alpha\to e\alpha\] induziert, wobei $e$ der Einsschnitt von $G_{\mathcal{Y}_{(\nu-1)/e}}$ ist und $S\to\mathcal{Y}_{\nu/e}$ étale ist.\par        
Aus dieser exakten Sequenz erhält man die exakte Kohomologiesequenz \[H^{1}_{\acute{e}t}(\mathcal{Y}_{\nu/e},i_{*}(\mathscr{L}_{(\nu-1)/e}\otimes_{\mathcal{O}_{\mathcal{Y}_{(\nu-1)/e}}}\mathcal{I}))\xrightarrow{h} H^{1}_{\acute{e}t}(\mathcal{Y}_{\nu/e},G_{\mathcal{Y}_{\nu/e}})\to H^{1}_{\acute{e}t}(\mathcal{Y}_{\nu/e},i_{*}G_{\mathcal{Y}_{(\nu-1)/e}}).\]
Nach \cite[Corollaire VIII.5.6]{SGA 4} ist nun $R^{1}i_{*}G_{\mathcal{Y}_{(\nu-1)/e}}=0$, da $i$ als abgeschlossene Immersion insbesondere ganz ist.\par
Nach \cite[Proposition V.3.1.3]{G} ist damit $H^{1}_{\acute{e}t}(\mathcal{Y}_{\nu/e},i_{*}G_{\mathcal{Y}_{(\nu-1)/e}})$ isomorph zu $H^{1}_{\acute{e}t}(\mathcal{Y}_{(\nu-1)/e},G_{\mathcal{Y}_{(\nu-1)/e}})$, so dass also die Sequenz  \[H^{1}_{\acute{e}t}(\mathcal{Y}_{\nu/e},i_{*}(\mathscr{L}_{(\nu-1)/e}\otimes_{\mathcal{O}_{\mathcal{Y}_{(\nu-1)/e}}}\mathcal{I}))\xrightarrow{h} H^{1}_{\acute{e}t}(\mathcal{Y}_{\nu/e},G_{\mathcal{Y}_{\nu/e}})\xrightarrow{i^{*}} H^{1}_{\acute{e}t}(\mathcal{Y}_{(\nu-1)/e},G_{\mathcal{Y}_{(\nu-1)/e}})\] exakt ist.\par
Wie oben gilt nun nach \cite[Corollaire VIII.5.6]{SGA 4} \[R^{1}i_{*}(\mathscr{L}_{(\nu-1)/e}\otimes_{\mathcal{O}_{\mathcal{Y}_{(\nu-1)/e}}}\mathcal{I})=0.\] Nach \cite[Proposition V.3.1.3]{G} hat man somit einen Isomorphismus \begin{equation*}
 \xymatrix@=3em{%
H^{1}_{\acute{e}t}(\mathcal{Y}_{(\nu-1)/e},\mathscr{L}_{(\nu-1)/e}\otimes_{\mathcal{O}_{\mathcal{Y}_{(\nu-1)/e}}}\mathcal{I})\ar[r]^{j}_{\sim} & H^{1}_{\acute{e}t}(\mathcal{Y}_{\nu/e},i_{*}(\mathscr{L}_{(\nu-1)/e}\otimes_{\mathcal{O}_{\mathcal{Y}_{(\nu-1)/e}}}\mathcal{I})),}\end{equation*} so dass also die Sequenz  \[H^{1}_{\acute{e}t}(\mathcal{Y}_{(\nu-1)/e},\mathscr{L}_{(\nu-1)/e}\otimes_{\mathcal{O}_{\mathcal{Y}_{(\nu-1)/e}}}\mathcal{I})\xrightarrow{h\circ j} H^{1}_{\acute{e}t}(\mathcal{Y}_{\nu/e},G_{\mathcal{Y}_{\nu/e}})\xrightarrow{i^{*}} H^{1}_{\acute{e}t}(\mathcal{Y}_{(\nu-1)/e},G_{\mathcal{Y}_{(\nu-1)/e}})\] exakt ist.\par
Ist nun $\Omega$ die Klasse von $\pi_{\nu/e}^{*}P_{\nu/e}$ in $H^{1}_{\acute{e}t}(\mathcal{Y}_{\nu/e},G_{\mathcal{Y}_{\nu/e}})$, so bildet $i^{*}$ die Klasse $\Omega$ auf die Klasse von $i^{*}\pi_{\nu/e}^{*}P_{\nu/e}=\pi^{*}_{(\nu-1)/e}P_{(\nu-1)/e}$, d.\,h. auf die triviale Klasse in $H^{1}_{\acute{e}t}(\mathcal{Y}_{(\nu-1)/e},G_{\mathcal{Y}_{(\nu-1)/e}})$, ab.\par Aufgrund der Exaktheit der Sequenz existiert also eine Klasse $A$ in $H^{1}_{\acute{e}t}(\mathcal{Y}_{(\nu-1)/e},\mathscr{L}_{(\nu-1)/e}\otimes_{\mathcal{Y}_{(\nu-1)/e}}\mathcal{I})$, so dass $(h\circ j)(A)=\Omega$ ist.\par
Die kanonische $\mathcal{O}_{\mathcal{Y}_{(\nu-1)/e}}$"=Modulstruktur von $\mathscr{L}_{(\nu-1)/e}$ und $\mathcal{I}$ induziert ferner die folgende exakte Sequenz von Garben auf $(\mathcal{Y}_{(\nu-1)/e})_{\acute{e}t}$: 
\[0\to\ker(\alpha)\to\mathcal{I}\xrightarrow{\alpha}\mathscr{L}_{(\nu-1)/e}\otimes_{\mathcal{O}_{\mathcal{Y}_{(\nu-1)/e}}}\mathcal{I}\to0.\] Dabei ist $\alpha$ der kanonische Morphismus. Dies liefert eine Surjektion \begin{equation*}\xymatrix@=2em{%
H^{1}_{\acute{e}t}(\mathcal{Y}_{(\nu-1)/e},\mathcal{I})\ar@{>>}[rrr]^-{\alpha} & & & H^{1}_{\acute{e}t}(\mathcal{Y}_{(\nu-1)/e},\mathscr{L}_{(\nu-1)/e}\otimes_{\mathcal{Y}_{(\nu-1)/e}}\mathcal{I}),}\end{equation*} da $\mathcal{Y}_{(\nu-1)/e}$ eindimensional ist und nach \cite[Proposition III.3.3]{M} \[H^{2}_{\acute{e}t}(\mathcal{Y}_{(\nu-1)/e},\mathscr{F})\cong H^{2}(\mathcal{Y}_{(\nu-1)/e},\mathscr{F}_{Zar})\] für jede étale Garbe $\mathscr{F}$ auf $\mathcal{Y}_{(\nu-1)/e}$ gilt, wobei $\mathscr{F}_{Zar}$ die Einschränkung auf die Zariski-Topologie ist.\par
Die $\mathcal{O}_{\mathcal{Y}_{(\nu-1)/e}}$"=Modulstruktur von $\mathcal{I}$ liefert außerdem eine kanonische exakte Sequenz \[0\longto\ker g\longto\mathcal{O}_{\mathcal{Y}_{(\nu-1)/e}}\stackrel{g}{\longrightarrow}\mathcal{I}\longto 0\] von Garben auf $(\mathcal{Y}_{(\nu-1)/e})_{\acute{e}t}$. Diese liefert mit derselben Begründung wie oben eine Surjektion \begin{equation*}\xymatrix@=3em{%
H^{1}_{\acute{e}t}(\mathcal{Y}_{(\nu-1)/e},\mathcal{O}_{\mathcal{Y}_{(\nu-1)/e}})\ar@{>>}[rrr]^{g} &&& H^{1}_{\acute{e}t}(\mathcal{Y}_{(\nu-1)/e},\mathcal{I}).}\end{equation*} 
\par
Ist nun das Geschlecht von $Y_{K}$, womit die generische Faser des Modells $\mathcal{Y}^{0}$ bezeichnet werde, gleich Null, so folgt wie in \cite{DW1}, dass \[H^{1}_{Zar}(\mathcal{Y}_{(\nu-1)/e},\mathcal{O}_{\mathcal{Y}_{(\nu-1)/e}})=0\] ist:\par Da man ohne Beschränkung der Allgemeinheit annehmen darf, dass $Y_{K}$ einen rationalen Punkt besitzt (andernfalls gehe zu einer geeigneten endlichen Erweiterung $K'$ von $K$ über), ist nämlich in diesem Fall nach \cite[Proposition 7.4.1]{L1} \mbox{$Y_{K}=\mathbb{P}^{1}_{K}$}. Nach \cite[Example 5.3.27]{L1} ist dann \mbox{$\chi(Y_{K},\mathcal{O}_{Y_{K}})=1$.} Da $\mathfrak{o}_{K}$ ein diskreter Bewertungsring ist, folgt damit aus \cite[Proposition 5.3.28]{L1}, dass auch $\chi(\mathcal{Y}_{\kappa},\mathcal{O}_{Y_{\kappa}})=1$ ist, wobei $\mathcal{Y}_{\kappa}$ die spezielle Faser von $\mathcal{Y}^{0}$ bezeichne.\par Da $(\lambda_{0})_{*}\mathcal{O}_{\mathcal{Y}_{0}}=\mathcal{O}_{\Spec\,\mathfrak{o}_{K}}$ gilt, wobei $\lambda_{0}$ den Strukturmorphismus von $\mathcal{Y}_{0}$ bezeichne, ist nun $H^{0}(\mathcal{Y}_{0},\mathcal{O}_{\mathcal{Y}_{0}})=\mathfrak{o}_{K}$, also nach \cite[Corollary 5.3.22]{L1} $H^{0}(\mathcal{Y}_{\kappa},\mathcal{O}_{\mathcal{Y}_{\kappa}})=\kappa$, so dass $H^{1}(\mathcal{Y}_{\kappa},\mathcal{O}_{\mathcal{Y}_{\kappa}})=0$ aus $\chi(\mathcal{Y}_{\kappa},\mathcal{O}_{Y_{\kappa}})=1$ folgt. Nach \cite[p. 53, Corollary 3]{Mf2} ist also $H^{1}(\mathcal{Y}_{\nu'/e},\mathcal{O}_{\mathcal{Y}_{\nu'/e}})=0$ für alle $\nu'\in\nat$.\par
Da also \[H^{1}_{Zar}(\mathcal{Y}_{(\nu-1)/e},\mathcal{O}_{\mathcal{Y}_{(\nu-1)/e}})=0\] ist, ist nach \cite[Proposition I.2.5.]{FK} auch \[H^{1}_{\acute{e}t}(\mathcal{Y}_{(\nu-1)/e},\mathcal{O}_{\mathcal{Y}_{(\nu-1)/e}})=0.\]
Nach dem bisher gezeigten ist dann aber $A=0$ und damit auch $\Omega=0$, so dass also $\pi_{\nu/e}^{*}P_{\nu/e}$ trivial ist. Im Fall, dass das Geschlecht von $Y_{K}$ gleich Null ist, kann man daher $\mu=\pi$ wählen.\par
Sei also das Geschlecht von $Y_{K}$ von Null verschieden. Dann gibt es nach dem Beweis von \cite[Theorem 11]{DW1} einen Morphismus \[\rho\colon(\mu\colon\mathcal{Z}\to\mathfrak{X})\to(\pi\colon\mathcal{Y}\to\mathfrak{X})\] in $S_{\mathfrak{X},D}$, so dass die induzierte Abbildung \[\rho^{*}\colon H^{1}(\mathcal{Y},\mathcal{O}_{\mathcal{Y}})\to H^{1}(\mathcal{Z},\mathcal{O}_{\mathcal{Z}})\] die Relation \[\rho^{*}(H^{1}(\mathcal{Y},\mathcal{O}_{\mathcal{Y}}))\subset p^{\nu/e}H^{1}(\mathcal{Z},\mathcal{O}_{\mathcal{Z}})\] erfüllt. Nach \cite[Proposition I.2.5]{FK} hat man dann dieselbe Aussage für den Fall, dass man die obigen Kohomologiegruppen für die Zariski-Topologie durch die entsprechenden Kohomologiegruppen für die étale Topologie ersetzt. Der 
Morphismus sei in diesem Fall mit $\rho_{\acute{e}t}^{*}$ bezeichnet. Wie im Beweis von \cite[Theorem 16]{DW1} folgt, dass die induzierte Abbildung \[(\rho_{\acute{e}t})_{(\nu-1)/e}^{*}\colon H^{1}_{\acute{e}t}(\mathcal{Y}_{(\nu-1)/e},\mathcal{O}_{\mathcal{Y}_{(\nu-1)/e}})\to H^{1}_{\acute{e}t}(\mathcal{Z}_{(\nu-1)/e},\mathcal{O}_{\mathcal{Z}_{(\nu-1)/e}})\] trivial ist. Aus der Kommutativität des Diagramms \begin{equation*}
 \xymatrix@=3em{%
H^{1}_{\acute{e}t}(\mathcal{Y}_{(\nu-1)/e},\mathcal{O}_{\mathcal{Y}_{(\nu-1)/e}})\ar[rrr]^{h\circ j\circ \alpha\circ g} \ar[d]^{(\rho^{*}_{\acute{e}t})_{(\nu-1)/e}=0} & && H^{1}_{\acute{e}t}(\mathcal{Y}_{\nu/e},G_{\mathcal{Y}_{\nu/e}}) \ar[d]^{(\rho^{*}_{\acute{e}t})_{\nu/e}}\\
H^{1}_{\acute{e}t}(\mathcal{Z}_{(\nu-1)/e},\mathcal{O}_{\mathcal{Z}_{(\nu-1)/e}})\ar[rrr]^{h\circ j\circ \alpha\circ g} & & & H^{1}_{\acute{e}t}(\mathcal{Z}_{\nu/e},G_{\mathcal{Y}_{\nu/e}})}\end{equation*} folgt dann, dass $(\rho_{\acute{e}t})_{\nu/e}^{*}\Omega$ trivial ist, da nach dem gezeigten $\Omega$ ein Urbild in $H^{1}_{\acute{e}t}(\mathcal{Y}_{(\nu-1)/e},\mathcal{O}_{\mathcal{Y}_{(\nu-1)/e}})$ unter $h\circ j\circ\alpha\circ g$ besitzt.\par
Also ist auch $\mu_{\nu/e}^{*}P_{\nu/e}=(\rho_{\acute{e}t})_{\nu/e}^{*}(\pi_{\nu/e}^{*}P_{\nu/e})$ ein trivialer Torseur auf $\mathcal{Z}_{\nu/e}$ für die étale Topologie und damit auch für die fppf"=Topologie, was zu zeigen war. \hfill$\Box$}\end{itemize}}\par

\section{Trivialisierbarkeit étaler $G$"=Torseure}

Für eine Charakterisierung der Kategorie $\mathscr{B}_{\mathfrak{X}_{\mathfrak{o}},D}(G)$ ohne Zuhilfenahme der Überdeckungskategorien $S_{\mathfrak{X},D}, S_{\mathfrak{X},D}^{good}$ bzw. $S^{ss}_{\mathfrak{X},D}$ benötigen wir die folgende Definion:\\\par

\begin{defn}
Es sei $G$ eine zusammenhängende algebraische Gruppe über einem endlichen Körper $\mathbb{F}_{q}$ und $X$ ein beliebiges $\mathbb{F}_{q}$"=Schema. Dann heißt ein étaler Torseur $P$ trivialisierbar, wenn es einen endlichen und surjektiven Morphismus \mbox{$f\colon Y\to X$} gibt, so dass $f^{*}P$ ein trivialer $G$"=Torseur auf $Y$ ist.
\end{defn}\par

Im folgenden sei nun stets zus\"atzlich vorausgesetzt, dass $G$ eine zusammenh\"angende reduktive Gruppe ist. Dazu sei an die Definition einer reduktiven Gruppe erinnert:\\\par
\begin{defn}[\protect{\cite[Exposé XIX, 1.2, 1.6, Définition 2.7]{SGA3}}]
\begin{enumerate}
\item[(i)] Es sei $G$ ein glattes, affines und zusammenhängendes Gruppenschema über einem algebraisch abgeschlossenen Körper $k$. Dann bezeichnet man die reduzierte Untergruppe $rad(G)$, die zu der neutralen Komponente des Schnitts der Boreluntergruppen von $G$ assoziiert ist, als das \textit{Radikal von $G$}. Der unipotente Teil $rad^{u}(G)$ von $rad(G)$ wird als das unipotente Radikal von $G$ bezeichnet. Dabei ist eine \textit{Boreluntergruppe} $B$ von $G$ eine maximale auflösbare Untergruppe und der \textit{unipotente} Teil einer Gruppe $H$ besteht aus den Elementen $h$, für die der Morphismus $h-1$ nilpotent ist (vgl. jeweils \cite{B}).
\item[(ii)]	Ein Gruppenschema $G$ über einem algebraisch abgeschlossenen Körper $k$ heißt \textit{reduktiv}, falls es affin, glatt und zusammenhängend ist und desweiteren $rad^{u}G=e$ gilt. 
\item[(iii)] Ein Gruppenschema $G$ über einer beliebigen Basis $S$ heißt \textit{reduktiv}, falls es affin und glatt über $S$ ist und seine geometrischen Fasern zusammenhängend und reduktiv sind.
\end{enumerate}
\end{defn}
\par

Dann liefert uns das folgende Resultat von Deligne einen Zusammenhang zwischen trivialisierbaren étalen $G$"=Torseuren und solchen, die isomorph zu ihrem Pullback unter dem absoluten Frobenius sind:\\\par

\begin{sat}[vgl. \protect{\cite[Beweis von Lemma 3.3.]{Las}}]\label{Laszlo}
Es sei $G$ eine glatte zusammenh\"angende algebraische Gruppe \"uber $\mathbb{F}_{q}$, wobei $q=p^{r}$ sei, $X$ ein Schema \"uber $\mathbb{F}_{q}$ und $F$ der absolute Frobenius, wobei wir in der Notation zur Vereinfachung derer nicht zwischen dem absoluten Frobenius von $X$ und dem von $G$ unterscheiden. Dann induziert die Einbettung \mbox{$G(\mathbb{F}_{q})= G^{F}\into G$} einen Funktor von der Kategorie der $G(\mathbb{F}_{q})$"=Torseure auf $X$ f\"ur die étale Topologie in die Kategorie der étalen $G$"=Torseure $P$ auf $X$, die als $G$"=Torseure isomorph zu $F^{*}P$ sind. Dieser Funktor ist eine Äquivalenz von Kategorien.  
\end{sat}
\bewub{Ist $Q$ ein $G(\mathbb{F}_{q})$"=Torseur auf $X$, so ist $P(Q)=Q\wedge^{G^{F}}G$ ein $G$-Torseur auf $X$. Nach Konstruktion ist dieser $G$"=Torseur $P(Q)$ als $G$"=Torseur isomorph zu seinem Pullback unter dem Frobenius $F$. Man bemerke dazu, dass $G^{F}$ glatt über $\mathbb{F}_{q}$ ist, so dass $Q$ darstellbar ist, die Bildung von eingeschränkten Produkten mit Basiswechsel vertauscht und $Q$ per Definition isomorph zu $F^{*}Q$ als $G^{F}$"=Torseur ist.\par
Es sei nun $P$ ein étaler $G$"=Torseur auf $X$, f\"ur den es einen Isomorphismus $\alpha\colon P\xrightarrow{\sim} F^{*}P$ von $G$"=Torseuren gibt. Dann lässt sich wie folgt ein kontravarianter Funktor $\mathscr{F}_{P,\alpha}$ von der Kategorie der étalen Schemata über $X$ in die Kategorie der abelschen Gruppen definieren:\par Für jedes étale $X$"=Schema $U$ setzt man \[\mathscr{F}_{P,\alpha}(U):=\{\sigma\in P(U); F^{*}\sigma=\alpha\circ\sigma\circ\iota_{U}\},\] wobei $\iota_{U}\colon F^{*}U\to U$ den kanonischen Morphismus bezeichne und man beachte, dass $F^{*}$ aufgrund seiner Funktorialität zu jedem Morphimus $\sigma\colon U\to P$ einen Morphismus $F^{*}\sigma\colon F^{*}U\to F^{*}P$ induziert. Dies ergibt einen Funktor $\mathscr{F}_{P,\alpha}$ von dem étalen Situs auf $X$ in die Kategorie der Mengen, der auch eine étale Garbe von Mengen auf $X$ ist, da die in der Definition von $\mathscr{F}_{P,\alpha}$ geforderte Bedingung an die lokalen Schnitte von $P$ verträglich mit den Garbenbedingungen ist. $\mathscr{F}_{P,\alpha}$ ist auch eine Garbe von abelschen Gruppen, da die Gruppenstruktur auf $P(U)$ mit Pullback unter dem Frobenius verträglich ist.\par 
Wir behaupten nun, dass $\mathscr{F}_{P,\alpha}$ ein $G(\mathbb{F}_{q})$"=Torseur ist:\par Als erstes zeigen wir dazu, dass $\mathscr{F}_{P,\alpha}$ ein formell-prinzipalhomogener Raum unter $G_{X}(\mathbb{F}_{q})$ ist, d.h. dass $G^{F}$ auf $\mathscr{F}_{P,\alpha}$ operiert und dann der kanonische Morphismus \[\mathscr{F}_{P,\alpha}(U)\otimes_{X(U)}G_{X}(\mathbb{F}_{q})(U)\xrightarrow{\sim}\mathscr{F}_{P,\alpha}(U)\otimes_{X(U)}\mathscr{F}_{P,\alpha}(U), (x,g)\mapsto(x,xg)\] für alle Schemata $U$, die étale über $X$ sind, ein Isomorphismus ist.\\
Zunächst zeigen wir dazu, dass sich die durch die Aktion von $G_{X}(U)$ auf $P(U)$ gegebene Aktion von $G_{X}(\mathbb{F}_{q})(U)$ auf $P(U)$ zu einer Aktion von $G_{X}(\mathbb{F}_{q})(U)$ auf $\mathscr{F}_{P,\alpha}(U)$ einschränkt, also \mbox{$\sigma\cdot g\in\mathscr{F}_{P,\alpha}(U)$} für alle \mbox{$\sigma\in\mathscr{F}_{P,\alpha}(U)$} und \mbox{$g\in G_{X}(\mathbb{F}_{q})$} gilt.\par 
Da Pullback unter dem Frobenius verträglich mit der Bildung von Faserprodukten ist, ist $F^{*}(\sigma\cdot g)=F^{*}(\sigma)\cdot g$. Nach Definition von $\mathscr{F}_{P,\alpha}$ ist ferner $F^{*}(\sigma)=\alpha\circ\sigma\circ\iota_{U}$. Weil außerdem $\alpha$ nach Voraussetzung ein Isomorphismus von $G$"=Torseuren ist und also $\alpha(yg)=\alpha(y)g$ für alle $y\in P$ und $g\in G_{X}$ gilt, folgt \[F^{*}(\sigma\cdot g)=F^{\sigma}\cdot g=(\alpha\circ\sigma\iota_{U})\cdot g=\alpha\circ(\sigma\cdot g)\circ\iota_{U},\] so dass $\sigma\cdot g\in\mathscr{F}_{P,\alpha}(U)$ ist.\par
Also schränkt sich der kanonische Morphismus \[\beta_{P}\colon P(U)\otimes_{X(U)}G_{X}(U)\to P(U)\otimes_{X(U)}P(U), (\sigma,g)\mapsto(\sigma,\sigma\cdot g)\] zu einem Morphismus  \[\beta_{\mathscr{F}_{P,\alpha}}\colon\mathscr{F}_{P,\alpha}(U)\otimes_{X(U)}G_{X}(\mathbb{F}_{q})(U)\to \mathscr{F}_{P,\alpha}(U)\otimes_{X(U)}\mathscr{F}_{P,\alpha}(U),(\sigma,g)\mapsto(\sigma,\sigma\cdot g)\] ein. Da ersterer ein Isomorphismus ist, ist damit auch $\beta_{\mathscr{F}_{P,\alpha}}$ ein Isomorphismus und also $\mathscr{F}_{P,\alpha}$ ein formell"=prinzipalhomogener Raum unter $G_{X}(\mathbb{F}_{q})$.\par  
Es bleibt noch zu zeigen, dass $\mathscr{F}_{P,\alpha}(U)\neq\emptyset$ für ein gewisses étales Schema $U$ über $X$ ist, so dass insgesamt $\mathscr{F}_{P,\alpha}$ ein étaler $G(\mathbb{F}_{q})$"=Torseur auf $X$ ist. \par
Dazu benutzen wir die Lang-Isogenie (siehe \cite{La}, \cite[3.4.]{Gi}) \[\lambda\colon G\to G, \lambda(\sigma)=\sigma^{-1}F(\sigma),\] die étale und surjektiv und eine Galoisüberlagerung von $G$ mit Gruppe \linebreak\mbox{$G(\mathbb{F}_{q})=G^{F}$} ist (vgl. \cite[Théorème 3.5.]{Gi}, man beachte, dass dort $G$ zusätzlich als reduktiv vorausgesetzt wird, aber dies nicht für die hier zitierten Aussagen benötigt wird), aber im allgemeinen kein Gruppenhomomorphismus.
Es sei nun $U$ ein étales Schema über $X$, bezüglich dessen der Torseur $P$ trivialisiert, d.\,h. es sei $P\times_{X}U$ isomorph zu dem trivialen $G_{X}\times_{X}U$"=Torseur $G_{X}\times_{X}U$. Dann ist $\alpha\otimes\id_{U}\in (F^{*}(G_{X}\times_{X}U))(G_{X}\times_{X}U)=G_{U}(G_{U})$ und die Surjektivität der Lang-Isogenie liefert, dass es ein $\sigma\in G_{U}(G_{U})$ gibt, so dass $\alpha\otimes\id_{U}=\sigma^{-1}F(\sigma)$ gilt.\par
Dass die beiden Abbildungen $Q\longmapsto(P(Q),P(Q)\xrightarrow{\sim}F^{*}(P(Q))$ und \linebreak \mbox{$(P,\alpha)\longmapsto\mathscr{F}(P)$} jeweils funktoriell sind, ist offensichtlich, so dass sie also Funktoren von der Kategorie der $G^{F}$"=Torseure auf $X$ in die Kategorie der Paare $(P,\alpha)$ aus einem $G$"=Torseur $P$ auf $X$ und einem Isomorphismus $\alpha\colon P\xrightarrow{\sim}F^{*}P$ von $G$"=Torseuren bzw. in umgekehrter Richtung von ebendieser Kategorie in die Kategorie der $G^{F}$"=Torseure auf $X$ definieren.\par
Dass die beiden Abbildungen $Q\longmapsto (P(Q),P(Q)\xrightarrow{\sim}F^{*}(P(Q)))$ und \linebreak $(P,\alpha)\longmapsto\mathscr{F}(P)$ invers zu einander sind, ergibt sich schließlich direkt aus ihrer oben dargestellten Konstruktion.\hfill$\Box$}\\\par

\begin{bem}
 Das originale Resultat von Deligne setzt nur voraus, dass $G$ eine zusammenhängende algebraische Gruppe ist, d.\,h. $G$ ist dort nicht notwendigerweise glatt. Der Beweis wird allerdings in \cite{Las} lediglich skizziert.\\
\end{bem}\par

\begin{kor}\label{pullback2}
Es sei $G$ ein zusammenhängendes und reduktives Gruppenschema von endlicher Präsentation über $\mathbb{F}_{q}$, $q=p^{r}$, und $X$ ein Schema über $\mathbb{F}_{q}$.\par
Dann ist jeder étale $G$"=Torseur $P$ auf $X$, der als étaler $G$"=Torseur $P$ isomorph zu seinem Pullback $F^{*}P$ unter dem absoluten Frobenius von $X$ ist, durch einen endlichen étalen Morphismus $f\colon Y\to X$ trivialisierbar. 
\end{kor}
\bewub{Es sei $P$ ein étaler $G$"=Torseur auf $X$, der als étaler $G$"=Torseur isomorph zu seinem Pullback $F^{*}P$ unter dem Frobenius ist. Dann existiert nach Satz \ref{Laszlo} ein $G^{F}$"=Torseur $Q$, so dass $P=Q\wedge^{G^{F}}G$ ist. Weil $G^{F}=G(\mathbb{F}_{q})$ glatt und endlich ist, ist $Q$ nach \cite[Proposition III.4.2]{M} endlich über $X$, so dass $f\colon Q\to X$ eine endliche étale Überlagerung von $X$ ist. Man beachte dabei, dass unter den gegebenen Voraussetzungen gilt, dass jeder $G$"=Torseur für die fppf"=Topologie auch ein Torseur für die étale Topologie ist und umgekehrt.\par
Es ist nun $f^{*}Q=Q\times_{X}Q$ trivial als étaler $G^{F}\times_{X}Q$"=Torseur auf $Q$, so dass auch $f^{*}P=P\times_{X}Q=(Q\times_{X}Q)\wedge^{G^{F}}G$ trivial als étaler $G$"=Torseur ist.\hfill$\Box$}\\\par

\begin{kor}\label{potenz}
Es sei $G$ ein zusammenhängendes und reduktives Gruppenschema von endlicher Präsentation über $\mathbb{F}_{q}$, $q=p^{r}$, und $X$ ein Schema über $\mathbb{F}_{q}$. Es sei ferner $P$ ein étaler $G$"=Torseur auf $X$ mit der Eigenschaft, dass es ein $s\in\nat$ gibt, so dass $(F^{t})^{*}P\cong(F^{s})^{*}P$ für alle $t\geq s$ als étale $G$"=Torseure gilt, wobei $F$ den absoluten Frobenius von $X$ bezeichne. Dann ist $P$ trivialisierbar.
\end{kor} 
\bewub{Es ist nach Voraussetzung $(F^{s})^{*}P\cong(F^{s+1})^{*}P=F^{*}((F^{s})^{*}P)$ als étale $G$"=Torseure auf $X$. Also ist nach Korollar \ref{pullback2} $(F^{s})^{*}P$ durch einen endlichen étalen Morphismus $f:Y\to X$ trivialisierbar. Damit wird aber auch $P$ durch den endlichen Morphismus $F^{s}\circ f$ trivialisiert.\hfill$\Box$}\\\par

\section{Semistabilität von Prinzipalbündeln}

Ähnlich wie im Fall der Vektorbündel werden uns diese Resultate eine Beschreibung der Kategorie $\mathscr{B}_{\mathfrak{X}_{\mathfrak{o},D}}(G)$ mittels des Begriffs der \emph{Semistabilität} ermöglichen. Dazu benötigen wir aber noch die folgende Konstruktion:\\\par
Ist $G$ eine zusammenhängende reduktive Gruppe über einem Körper $k$ beliebiger Charakteristik, $X$ eine glatte und projektive Kurve über $k$, $E$ ein $G$"=Torseur auf $X$ und $F$ ein quasi"=projektives Schema über $k$, auf dem $G$ von links operiert, so wird durch $g(e,f):=(e(g\otimes 1),g^{-1}f)$ für alle \mbox{$g\in G$}, $e\in E$ und $f\in F$ eine Aktion von $G$ auf $E\times_{\Spec\,k}F$ definiert. Der Quotient von $E\times_{\Spec\,k}F$ modulo dieser Aktion von $G$ ist durch ein eindeutig bestimmtes Schema $E(F)$ darstellbar, so dass der Morphismus \mbox{$E\times_{\Spec\,k}F\to E(F)$} dem Faserprodukt $E\times_{\Spec\,k}F$ die Struktur eines étalen $G$"=Torseurs über $E(F)$ gibt. Dies sieht man wie folgt (vgl. \cite[Proposition 4]{Se}):\par
Da $G$ reduktiv, also insbesondere glatt und affin, ist und folglich $E$ nach der Argumentation aus Satz \ref{etale} auch ein $G$"=Torseur auf $X$ für die étale Topologie ist, darf man aufgrund der Tatsache, dass die Existenz von $E(F)$ ein lokales Problem ist, annehmen, dass es eine étale Überlagerung $f\colon Y\to X$ gibt, so dass $f^{*}E$ trivial als étaler $G$"=Torseur ist. Nach Bemerkung \ref{galo} darf man ohne Beschränkung der Allgemeinheit voraussetzen, dass $f$ galoissch mit Galoisgruppe $H$ ist. Also ist $E\cong Y\times_{\Spec\,k}G/H$.\par
Mittels der Identifikation von \v{C}ech"=Kohomologie mit Galois"=Kohomologie aus \cite[Example III.2.6]{M} sieht man daraus, dass $E(F)=Y\times_{\Spec\,k}F/H$ gilt. Nach \cite[Exercise 3.3.23]{L1} ist aber $Y\times_{\Spec\,k}F/H$ durch ein Schema darstellbar.\par
Aus der Konstruktion ist offensichtlich, dass man für einen Morphismus $\alpha\colon F_{1}\to F_{2}$ von quasi"=projektiven $k$"=Schemata mit einer Aktion von $G$ von links, der verträglich mit der Aktion von $G$ auf $F_{1}$ bzw. $F_{2}$ ist, stets einen Morphismus $\alpha_{E}\colon E(F_{1})\to E(F_{2})$ erhält.\par
Außerdem halten wir fest, dass, wenn $\sigma\colon G\to H$ ein Morphismus zusammenhängender reduktiver algebraischer Gruppen über $k$ ist und man vermittels dieses $H$ als quasi"=projektive Varietät mit einer Aktion von $G$ von links auffasst, das wie oben definierte Schema $E_{\rho}(H)$ ein $H$"=Torseur ist.\par
Ist weiterhin $H$ eine normale abgeschlossene Untergruppe von $G$, so folgt aus dem Theorem von Chevalley, dass $G/H$ eine quasi"=projektive Varietät ist (vgl. z.B. \cite[Corollary 1.2]{C}). Versieht man nun $G/H$ mit der kanonischen Aktion von $G$ von links, so existiert also das oben definierte Schema $E(G/H)$.\\\par
Mittels dieser Konstruktion von Serre können wir nun Erweiterungen und Reduktionen der Strukturgruppe definieren:\\\par

    \begin{defn}
Ist $\rho\colon G\to H$ ein Morphismus von zusammenhängenden reduktiven Gruppenschemata von endlicher Präsentation über einem Körper $k$ und $X$ eine glatte und projektive Kurve über $k$, so induziert $\rho$ eine Abbildung \[H^{1}(X,G)\to H^{1}(X,H),\] in dem man jedem $G$"=Torseur $E$ den $H$"=Torseur $E_{\rho}(H)$ zuordnet. Man sagt, dass man $E_{\rho}(H)$ mittels Erweiterung der Strukturgruppe von $E$ zu $H$ erhält.\par
Gibt es umgekehrt zu einem gegebenen $H$"=Torseur $E'$ einen $G$"=Torseur $E$ und einen Isomorphismus von $H$"=Torseuren \mbox{$\Phi\colon E_{\rho}(H)\to E'$}, so bezeichnet man $(E,\Phi)$ als eine Reduktion der Strukturgruppe von $E'$ zu $G$.\\
\end{defn}\par

Für die weitere Beschreibung der Kategorie $\mathscr{B}_{\mathfrak{X}_{\mathfrak{o}}, D}(G)$ benötigen wir den Begriff der \emph{Semistabilität} von Prinzipalbündeln. Dazu sei zunächst an die Definition von \emph{Semistabilität} für Vektorbündel erinnert:\\\par

\begin{defn}[Mumford,\protect{\cite{Mf1}}]\label{mumford}
Es sei $\mathscr{F}$ ein Vektorbündel auf einer glatten und projektiven Kurve über einem Körper $k$. Dann bezeichnet man den Quotienten $\mu(\mathscr{F}):=\deg\mathscr{F}/\rk\mathscr{F}$ als den slope von $\mathscr{F}$ und es heißt $\mathscr{F}$ \emph{semistabil}, falls für jedes nichttriviale echte Unterbündel $\mathscr{G}$ von $\mathscr{F}$ die Ungleichung $\mu(\mathscr{G})\leq\mu(\mathscr{F})$ erfüllt ist.\\
\end{defn}\par

Man bemerke, dass jedes nichttriviale echte Unterbündel $\mathscr{G}$ von $\mathscr{F}$ Anlass gibt zu der Reduktion der Strukturgruppe des assoziierten $Gl_{\rk\mathscr{F}}$"=Prinzipalbündels zu der Untergruppe, die den Unterraum $k^{\rk\mathscr{G}}$ von $k^{\rk\mathscr{F}}$ festhält (vgl. den späteren Beweis von Lemma \ref{mumram}).\par
Die allgemeine Strukturtheorie algebraischer Gruppen zeigt, dass diese Untergruppe ein Spezialfall des allgemeinen Begriffs einer maximalen parabolischen Untergruppe einer gegebenen reduktiven algebraischen Gruppe ist (vgl. ebenfalls den Beweis von Lemma \ref{mumram}).\par
Dabei heißt nach \cite[Exposé XXVI, Définition 1.1]{SGA3} eine Untergruppe $P$ einer reduktiven algebraischen Gruppe $G$ über einem Basisschema $S$ \emph{parabolisch}, wenn zum einen $P$ glatt über $S$ ist, zum anderen $P_{\overline{s}}$ für alle $s\in S$ eine Boreluntergruppe von $G_{\overline{s}}$ enthält. \par
Es seien nun $X$ eine glatte und projektive Kurve über einem Körper $k$, $G$ eine zusammenhängende reduktive algebraische Gruppe über $k$ sowie $E$ ein $G$"=Torseur auf $X$. Man beachte dabei, dass $E$ unter den gegebenem Voraussetzungen ein Torseur sowohl für die fppf"=Topologie als auch für die étale Topologie ist. Es sei ferner $P$ eine  parabolische Untergruppe von $G$. Dann ist $G/P$ nach \cite[Exposé XXVI, Proposition 1.2]{SGA3} durch ein glattes und projektives $k$"=Schema darstellbar. Mittels \[G\times_{\Spec\,k}G/P\to G/P, (g,\overline{h})\mapsto \overline{gh}\] lässt sich eine Aktion von $G$ von links auf $G/P$ definieren. Also existiert das oben definierte Schema $E(G/P)$. Zur Vereinfachung der Notation sei im folgenden $E/P:=E(G/P)$.\par
Es sei ferner $T_{E/P}$ das relative Tangentialbündel bezüglich der Projektion \mbox{$E/P\to X$.} Dieses ist nach \cite[Exposé II, §3]{SGA3} als der kovariante Funktor \[M\mapsto \mathcal{H}om_{X}(\Spec(\mathcal{O}_{X}\oplus M), E/P)\] von der Kategorie der freien $\mathcal{O}_{X}$"=Moduln von endlichem Typ in die Kategorie der Funktoren über $E/P$ definiert, welcher nach \cite[Exposé II, Proposition 3.3]{SGA3} durch ein Vektorbündel auf $E/P$ darstellbar ist.\par
Ist $\sigma\colon X\to E/P$ ein Schnitt der kanonischen Projektion $E/P\to X$, so findet man (siehe z.B. \cite[S.2]{HN}), dass $\sigma^{*}T_{E/P}=E(\mathfrak{g}/\mathfrak{p})$ gilt, wobei $\mathfrak{g}=\Lie G$ und $\mathfrak{p}=\Lie P$ sei, $E$ als étaler $P$"=Torseur über $E/P$ aufgefasst werde und die Aktion von $P$ auf $\mathfrak{g}/\mathfrak{p}$ durch die Adjunktionsoperation gegeben ist.\par
Für jeden Schnitt \mbox{$\sigma\colon X\to E/P$} ist damit $\sigma^{*}T_{E/P}$ ein Vektorbündel auf $X$ vom Grad $\deg(\sigma^{*}T_{E/P})$.\par
Desweiteren gilt:\\\par

\begin{lem}\label{sorger}
Die Schnitte $\sigma\colon X\to E/P$ stehen in bijektiver Korrespondenz zu den étalen $P$"=Torseuren auf $X$, die man mittels Reduktionen der Strukturgruppe von $E$ zu $P$ erhält.
\end{lem} 
\bewub{Man fasse $E\to E/P$ als étalen $P$"=Torseur auf und betrachte für $\sigma\colon X\to E/P$ das Pullback"=Diagramm \begin{equation*}
        \xymatrix@=3em{%
         \sigma^{*}E \ar[r] \ar[d]^{P} & E \ar[d]^{P}\\
          X \ar[r]^{\sigma}& E/P. }
           \end{equation*} 
           Dieses definiert die gewünschte Reduktion der Strukturgruppe von $E$ zu $P$. \hfill$\Box$}\\

Daher ist es legitim, im folgenden Reduktionen der Strukturgruppe von $E$ zu $P$ mit $(P,\sigma)$ zu bezeichnen.\\\par

Dies motiviert die folgende Definition der Semistabilität für den Fall von Prinzipalbündeln:\\\par

\begin{defn}[\protect{\cite[Definition 1.1]{R}}]\label{semistabil}
Es sei $G$ eine zusammenhängende reduktive algebraische Gruppe über einem Körper $k$ und $X$ eine glatte und projektive Kurve über $k$. Dann heißt ein $G$"=Torseur $E$ auf $X$ semistabil, wenn für jede Reduktion $(P,\sigma)$ der Strukturgruppe zu einer maximalen parabolischen Untergruppe $P$, wobei $\sigma$ einen Schnitt $\sigma\colon X\to E/P$ bezeichnet, gilt, dass $\deg(\sigma^{*}T_{E/P})\geq 0$ ist. 
\end{defn}

Diese Definition der Semistabilität von Prinzipalbündeln nach Ramanathan ist verträglich mit der Definition von Semistabilität für Vektorbündeln von Mumford im folgenden Sinne:\\\par

\begin{lem}\label{mumram}
Ein $Gl(n)$"=Torseur $E$ auf einer glatten und projektiven Kurve $X$ über einem Körper $k$ ist genau dann semistabil, wenn für jedes Unterbündel $\mathscr{F}$ des assoziierten Vektorbündels $\mathscr{E}$ vom Rang $n$ gilt, dass \mbox{$\mu(\mathscr{F})\leq\mu(\mathscr{E})$} ist, d.\,h. wenn $\mathscr{E}$ ein semistabiles Vektorbündel im Sinne von Definition \ref{mumford} ist.
\end{lem}
\bewub{\begin{enumerate}
\item[(i)] Wir imitieren den Beweis der Tatsache für den Fall $\Char k=0$ aus \cite{HM}:\par Es sei $P(n,k)$ die Untergruppe von $Gl(n,k)$ der Matrizen der Form 
$ \left(\begin{array}{cc}
A & B\\
0 & C	\\
\end{array}\right)$ mit $A\in Gl(r,k)$ und $C\in Gl(n-r,k)$ bezüglich einer gewählten Basis von $k^{n}$. Dies ist eine parabolische Untergruppe, da $P(n,k)$ die Gruppe der oberen Dreiecksmatrizen enthält und diese wiederum eine Boreluntergruppe von $Gl(n,k)$ ist. Andererseits ist jede parabolische Untergruppe von $Gl(n,k)$ von dieser Gestalt.\par
Es sei $(P,\sigma)$ eine Reduktion der Strukturgruppe von $E$ zu $P$. Wie in Lemma \ref{sorger} liefert diese einen $P$"=Torseur $\sigma^{*}E$ auf $X$. Es sei $(U_{\alpha})_{\alpha\in I}$ eine Überdeckung von $X$, bezüglich derer $\sigma^{*}E$ trivialisiert. Dann definiert der $P$"=Torseur $\sigma^{*}E$ eine Klasse in $\check{H}^{1}_{\acute{e}t}(X,P)$ und der Kozykel zu der Überdeckung $(U_{\alpha})_{\alpha\in I}$ besitzt die Form $\left(\begin{array}{cc}
h_{\alpha\beta} & b_{\alpha\beta}\\
0 & q_{\alpha\beta}	\\
\end{array}\right)$.\par
Es sei $\mathscr{F}$ das Unterbündel vom Rang $r$ von $\mathscr{E}$, das zu $(h_{\alpha\beta})$ korrespondiert. Dann korrespondiert das Quotientenbündel $\mathscr{E}/\mathscr{F}$ zu $(q_{\alpha\beta})$ und wir behaupten, dass es einen natürlichen Isomorphismus
\[\sigma^{*}T_{E/P}\cong\mathscr{F}^{v}\otimes(\mathscr{E}/\mathscr{F})\] gibt.
Dies sieht man wie folgt:\par
Die Adjunktionsoperation von $P$ auf der Lie"=Algebra $\mathfrak{g}\mathfrak{l}(n,k)$ ist durch die Matrizenmultiplikation \[\left(\begin{array}{cc}
h & b\\
0 & q\\
\end{array}\right)\cdot\left(\begin{array}{cc}
A & B\\
0 & C\\
\end{array}\right)\cdot\left(\begin{array}{cc}
h^{-1} & b'\\
0 & q^{-1}	\\
\end{array}\right)=\left(\begin{array}{cc}
* & *\\
q\cdot C\cdot h^{-1} & *	\\
\end{array}\right)\] gegeben. Dabei ist $\left(\begin{array}{cc}
h^{-1} & b'\\
0 & q^{-1}	\\
\end{array}\right)$ die inverse Matrix zu $\left(\begin{array}{cc}
h & b\\
0 & q\\
\end{array}\right)$. Die Einträge $*$ gehören zu $\mathfrak{p}$, so dass der einzige Anteil, der für $\mathfrak{g}/\mathfrak{p}$ relevant ist, der Eintrag links unten ist. Also ist die Adjunktionsoperation von $P$ auf $\mathfrak{g}/\mathfrak{p}$ durch
\[(\left(\begin{array}{cc}
h & b \\
0 & q \\
\end{array}\right),C)\mapsto q\cdot C\cdot h^{-1}\] gegeben. Daher ist der Kozykel zu $\sigma^{*}E(\mathfrak{g}/\mathfrak{p})=\sigma^{*}T_{E/P}$ in $\check{H}^{1}_{\acute{e}t}(X,Gl(n,k))$ durch $((h^{-1}_{\alpha\beta})^{t}\otimes q_{\alpha\beta})$ bestimmt, der aber wiederum zu dem Vektorbündel $F^{v}\otimes(E/F)$ korrespondiert. (vgl. z.\,B. \cite[§0.5]{GH}, wobei man dort $\mathbb{C}$ stets durch einen beliebigen Körper beliebiger Charakteristik ersetzen darf.)
\item[(ii)] Mit dieser Vorarbeit können wir nun die Aussage des Lemmas beweisen:\par Es sei zunächst das zu $E$ assoziierte Vektorbündel $\mathscr{E}$ vom Rang $n$ semistabil im Sinne von Mumford und $(P,\sigma)$ eine Reduktion der Strukturgruppe von $E$ zu einer maximalen parabolischen Untergruppe $P$ von $Gl(n,k)$. Dann korrespondiert $P$ zu einer zweischrittigen Flagge \[0=V_{0}\subset V_{1}=(k^{n})^{P}\subset V_{2}=k^{n}.\]
Man setze nun $\mathscr{F}:=\sigma^{*}E(V_{1})$. Nach Konstruktion ist dies ein nichttriviales echtes Unterbündel von $\mathscr{E}$, so dass aufgrund der Semistabilität von $\mathscr{E}$ folgt, dass $\mu(\mathscr{F})\leq\mu(\mathscr{E})$ ist. Dies ist äquivalent dazu, dass $\mu(\mathscr{F})\leq\mu(\mathscr{E}/\mathscr{F})$ ist. Da aber nach den Rechenregeln für $\deg$ (siehe z.B. \cite[Exercise 2.4]{Lp}) \begin{align*} \deg(\mathscr{F}^{*}\otimes(\mathscr{E}/\mathscr{F}))& =-\deg\mathscr{F}\cdot\rg(\mathscr{E}/\mathscr{F})+\rg\mathscr{F}\cdot\deg(\mathscr{E}/\mathscr{F})\\ & -\deg\mathscr{F}\cdot(\rg\mathscr{E}-\rg\mathscr{F})+\rg\mathscr{F}\cdot(\deg\mathscr{E}-\deg\mathscr{F}) \\ & =(\mu(\mathscr{E}/\mathscr{F})-\mu(\mathscr{F}))\cdot\rg\mathscr{F}\cdot\rg(\mathscr{E}/\mathscr{F})\end{align*} gilt, ist damit $\deg(\mathscr{F}\otimes(\mathscr{E}/\mathscr{F}))\geq 0$. Nach $(i)$ ist damit \mbox{$\deg(\sigma^{*}T_{E/P})\geq 0$,} also $E$ semistabil als $Gl(n,k)$"=Torseur nach Definition \ref{semistabil}.\par
Es sei umgekehrt $E$ ein semistabiler $Gl(n,k)$"=Torseur. Dann ist jedes Untervektorbündel $\mathscr{F}$ von $\mathscr{E}$ von der Form $\mathscr{F}=\sigma^{*}(E(V_{1}))$ für eine Reduktion der Strukturgruppe $(P,\sigma)$ zu der parabolischen Untergruppe $P$, die zu der Flagge $0=V_{0}\subset V_{1}\subset V_{2}=k^{n}$ korrespondiert, wobei $\mathscr{F}$ den Unterraum $V_{1}$ festhält. Nach $(i)$ ist $\sigma^{*}T_{E/P}=\mathscr{F}^{v}\otimes(\mathscr{E}/\mathscr{F})$. Weil nach Voraussetzung $E$ semistabil ist, ist also \[\deg(\mathscr{F}^{v}\otimes(\mathscr{E}/\mathscr{F}))=\deg\sigma^{*}T_{E/P}\geq 0.\] Wie oben gilt aber \[\deg(\mathscr{F}^{*}\otimes(\mathscr{E}/\mathscr{F}))=(\mu(\mathscr{E})-\mu(\mathscr{F}))\cdot\rg\mathscr{F}\cdot\rg(\mathscr{E}/\mathscr{F}),\] so dass $\mu(\mathscr{F})\leq\mu(\mathscr{E}/\mathscr{F})$ folgt und damit $\mathscr{E}$ ein semistabiles Vektorbündel im Sinne von Mumford ist. \hfill$\Box$\\ \end{enumerate}}\par

Wegen $\sigma^{*}T_{E/P}=E(\mathfrak{g}/\mathfrak{p})$ gilt außerdem:\\\par

\begin{lem}[vgl. \protect{\cite[Remark 2.2]{R}}, \protect{\cite[Lemma 2.5]{BH}}]\label{adsemi}
Es sei $G$ eine zusammenhängende reduktive algebraische Gruppe über einem Körper $k$ beliebiger Charakteristik und $X$ eine glatte und projektive Kurve über $k$. Dann ist ein $G$"=Torseur $E$ auf $X$ semistabil, falls das assoziierte Vektorbündel $E(\mathfrak{g})$ semistabil ist.
\end{lem}
\bewub{
Es sei $(P,\sigma)$ eine Reduktion der Strukturgruppe von $E$ zu einer maximalen parabolischen Untergruppe $P$ von $G$. Dann liefert die exakte Sequenz \begin{equation*}
        \xymatrix@=1em{%
         0 \ar[r] & \mathfrak{p} \ar[r] & \mathfrak{g} \ar[r] & \mathfrak{g}/\mathfrak{p} \ar[r] & 0}
           \end{equation*} wegen der Funktorialität der Konstruktion des assoziierten Faserbündels (siehe auch \cite[3.2, Exemple b),c)]{Se}) eine exakte Sequenz 
           \begin{equation*}
        \xymatrix@=1em{%
         0 \ar[r] & E(\mathfrak{p}) \ar[r] & E(\mathfrak{g}) \ar[r] & E(\mathfrak{g}/\mathfrak{p}) \ar[r] & 0}
           \end{equation*} von Vektorbündeln auf $X$. Da $E(\mathfrak{g})$ semistabil als Vektorbündel ist, ist also $\mu(E(\mathfrak{g}/\mathfrak{p}))\geq\mu(E(\mathfrak{g}))$, d.\,h.  $\mu(\sigma^{*}T_{E/P})\geq\mu(E(\mathfrak{g}))$. Weil nach \cite[Note 4.2]{Be} $\deg(E(\mathfrak{g}))=0$ gilt, so dass $\mu(E(\mathfrak{g}))=0$ ist, folgt daraus $\deg(\sigma^{*}T_{E/P})\geq 0$ und damit die Semistabilität von $E$. \hfill$\Box$}\\\par

\begin{kor}\label{torsortrivial}
Es sei $G$ eine zusammenhängende reduktive algebraische Gruppe über einem Körper $k$ und $X$ eine glatte und projektive Kurve über $k$. Dann ist der triviale $G$"=Torseur $G_{X}$ auf $X$ semistabil.
\end{kor}
\bewub{In der gegebenen Situtation gilt \[G_{X}(\mathfrak{g})=(G\times_{\Spec\,k}X)\times_{\Spec\,k}\mathfrak{g}/(x,eg,g^{-1}f)\cong X\times_{\Spec\,k}(G\times_{\Spec\,k}\mathfrak{g}/(eg,g^{-1}f))\] mit $x\in X$, $e\in G$ und $f\in\mathfrak{g}$, so dass $G_{X}(\mathfrak{g})$ ein triviales Vektorbündel ist. Nach Lemma \ref{adsemi} ist also $G_{X}$ semistabil.\hfill$\Box$}\\\par

\begin{bem}
Es seien $G, X$ und $k$ wie bisher.
Dann zeigt der Beweis von Lemma \ref{adsemi} auch, dass ein $G$"=Torseur $E$ auf $X$ genau dann semistabil ist, wenn \[\deg(\sigma^{*}E(\mathfrak{p}))\leq 0\] für jede Reduktion $(P,\sigma)$ der Strukturgruppe zu einer maximalen parabolischen Untergruppe $P$ von $G$ gilt.
\end{bem}\par

Wie im Fall der Vektorbündel lässt sich außerdem auch bei $G$"=Torseuren deren Semistabilität nach Pullback entlang einer endlichen Überlagerung nachprüfen:\\\par

\begin{lem}[vgl. \protect{\cite[Lemma 6.8]{BH}}]\label{endlichsemist}
Es sei $X$ eine glatte und projektive Kurve über einem Körper $k$, $G$ eine zusammenhängende reduktive algebraische Gruppe über $k$, $E$ ein $G$"=Torseur auf $X$ und $f\colon Y\to X$ eine endliche Überlagerung. Dann ist $E$ semistabil, falls $f^{*}E$ es ist.
\end{lem} 
\bewub{Es sei also $f^{*}E$ semistabil. Angenommen, $E$ sei nicht semistabil. Nach Definition \ref{semistabil} existieren somit eine maximale parabolische Untergruppe $P$ von $G$ und ein Schnitt $\sigma\colon X\to E/P$, so dass $\deg(\sigma^{*}T_{E/P})<0$ ist.\par
Dann ist aber auch $P_{Y}$ eine maximale parabolische Untergruppe von $G_{Y}$, da alle in der Definition der Eigenschaft "`maximal parabolisch"' auftretenden Bedingungen stabil unter Basiswechsel sind, und \mbox{$\sigma\otimes\id_{Y}:Y\to E\times_{X}Y/P_{Y}$} ein Schnitt des Strukturmorphismus $E\times_{X}Y/P_{Y}\to Y$. Da die Bildung von $E/P$ nach \cite[Remarque 3.2]{Se} und die des Tangentialbündels (siehe \cite[Proposition II.3.4]{SGA3})) jeweils mit Basiswechsel vertauschen, gilt nach den Rechenregeln für den Grad eines Vektorbündels (vgl. \cite[Proposition 7.3.8]{L1}) \[\deg((\sigma\otimes\id_{Y})^{*}T_{f^{*}E/P_{Y}})=\deg(f^{*}(\sigma^{*}T_{E/P}))=[K(Y):K(X)]\deg(\sigma^{*}T_{E/P}).\]  Also ist auch \[\deg((\sigma\otimes\id_{Y})^{*}T_{E\times_{X}Y/P_{Y}})<0\] und damit $E\times_{X}Y$ ebenfalls nicht semistabil, was einen Widerspruch zu der Voraussetzung, dass $E\times_{X}Y$ semistabil ist, darstellt. Also ist $E$ semistabil.
\hfill $\Box$}\\\par

Ferner gilt:\\\par

\begin{lem}\label{noether}
Es seien $X_{0}$ eine glatte und projektive Kurve über einem Körper $k$, $\overline{k}$ der algebraische Abschluss von $k$, $X:=X_{0}\otimes_{k}\overline{k}$ und $f:X\to X_{0}$ die kanonische Projektion. Ferner sei $G$ eine zusammenhängende reduktive algebraische Gruppe über $k$ und $E$ ein $G$"=Torseur auf $X_{0}$. Dann ist $E$ semistabil, falls $f^{*}E$ es ist.
\end{lem} 
\bewub{Es sei $f^{*}E$ semistabil. Angenommen, $E$ sei nicht semistabil. Dann existieren nach Definition \ref{semistabil} eine maximale parabolische Untergruppe $P$ von $G_{X_{0}}$ und ein Schnitt $\sigma\colon X\to E/P$, so dass $\deg(\sigma^{*}T_{E/P})<0$ ist.\par
Weil $f^{*}E=E\otimes_{k}\overline{k}$ ist und, wie schon im Beweis des vorherigen Resultats verwendet, alle Konstruktionen mit Basiswechsel verträglich sind, folgt mit \cite[Proposition 7.3.7]{L1} für die maximale parabolische Untergruppe \mbox{$P_{X}=P\otimes_{k}\overline{k}$} und den Schnitt $\sigma\otimes_{k}\overline{k}\colon X_{0}\to E\otimes_{k}\overline{k}/P_{X}$:
\[\deg((\sigma\otimes_{k}\overline{k})^{*}T_{f^{*}E/P_{X}})=\deg((\sigma^{*}T_{E/P})\otimes_{k}\overline{k})=\deg(\sigma^{*}T_{E/P})<0.\]
Damit ist auch $f^{*}E$ nicht semistabil, was einen Widerspruch zu der Voraussetzung darstellt, dass $f^{*}E$ semistabil ist. Also ist $E$ semistabil. \hfill$\Box$}\\\par

\begin{bem}
Es sei angemerkt, dass die Umkehrung in den beiden vorangegangenen Resultaten im allgemeinen nicht richtig ist: Ist $E$ semistabil, so ist im allgemeinen $f^{*}E$ nicht notwendigerweise ebenfalls semistabil.\par
Für Vektorbündel weiß man nach einem Resultat von Gieseker (siehe \cite[Lemma 1.1]{Gie}), dass die Umkehrung gilt, wenn der Morphismus $f$ endlich und separabel ist.\par
Für den Fall eines $G$"=Torseurs unter einer zusammenhängenden reduktiven algebraischen Gruppe $G$ von endlicher Präsentation über $k$ gilt nach \cite[Lemma 6.8]{BH} die Umkehrung, wenn ebenfalls $f$ endlich und separabel ist, aber dies nur für den Fall, dass $G$ eine Darstellung "`of low height"' besitzt.
\end{bem}

\section{Der Grad eines Torseurs} 

Wir benötigen nun noch den Begriff des Grades eines Torseurs:\\\par

\begin{defn}[vgl. \protect{\cite[Definition 3.2]{HN}}]
Es sei $X$ eine glatte und projektive Kurve vom Geschlecht $g$ über einem Körper $k$, $G$ eine zusammenhängende reduktive algebraische Gruppe über $k$ und $E$ ein $G$"=Torseur auf $X$. Dann heißt der Homomorphismus \[d_{E}\colon\Hom(G,k^{*})\to\mathbb{Z}, \chi\longmapsto\deg(E_{\chi})\] der Grad von $E$.\par
Dabei ist $E_{\chi}=E\times_{\Spec\,k}\mathbb{G}_{m}/G$ das durch den Charakter $\chi$ zu $E$ assoziierte Geradenbündel, wobei die Aktion von $G$ auf $E\times_{\Spec\,k}\mathbb{G}_{m}$ durch \[(e,f)g=(eg,\chi(g^{-1})f)\] für alle $e\in E$, $f\in\mathbb{G}_{m}$ und $g\in G$ gegeben ist.\\
\end{defn}\par

\begin{bem}\label{null}
Es sei weiterhin $X$ eine glatte und projektive Kurve vom Geschlecht $g$ über einem Körper $k$ und $G$ eine zusammenhängende reduktive algebraische Gruppe über $k$. Dann gilt:\par
\begin{enumerate}\item[(i)]Der Grad $d_{G_{X}}$ des trivialen Torseurs $G_{X}$ auf $X$ ist der Nullhomomorphismus.
\item[(ii)] Ist $E$ ein $G$"=Torseur auf $X$ und $f:Y\to X$ ein endlicher Morphismus glatter und projektiver Kurven über $k$, so gilt $d_{f^{*}E}=[K(Y):K(X)]d_{E}$.
\item[(iii)] Ist $E$ ein $G$"=Torseur auf $X$ und $\overline{k}$ ein algebraischer Abschluss von $k$, so gilt $d_{E\otimes_{k}\overline{k}}=d_{E}$. \end{enumerate}
\end{bem}
\bewub{\begin{enumerate} \item[(i)] Es sei $\chi$ ein beliebiger Charakter der Gruppe $G$. Dann erhält man als das zu $\chi$ assoziierte Geradenbündel $(G_{X})_{\chi}$ das folgende Geradenbündel:
\begin{align*}(G_{X})_{\chi} & =(G\times_{\Spec\,k}X\times_{Spec\,k}\mathbb{G}_{m})/<(e,x,f)g=(eg,x,\chi(g^{-1})f)> \\
 & = X\times_{\Spec\,k}(G\times_{\Spec\,k}\mathbb{G}_{m}/<(e,f)g=(eg,\chi(g^{-1})f)>),
 \end{align*} wobei $e,g\in G$, $x\in X$ und $f\in\mathbb{G}_{m}$ sind.\par
Also ist $(G_{X})_{\chi}$ ein triviales Geradenbündel auf $X$ und damit \[d_{G_{X}}(\chi)=\deg((G_{X})_{\chi})=0.\]
Da der Charakter $\chi$ beliebig gewählt war, ist damit der Grad $d_{G_{X}}$ des trivialen Torseurs $G_{X}$ auf $X$ der Nullhomomorphismus.
\item[(ii)] Dies folgt sofort aus der Tatsache, dass die Konstruktion des zu einem Charakter $\chi$ von $G$ assoziierten Geradenbündels mit Basiswechsel vertauscht, und den Rechenregeln für das Verhalten des Grades eines Geradenbündels unter endlichem Pullback.
\item[(iii)] Es gilt $\deg_{k}D=\deg_{\overline{k}}(D\otimes_{k}\overline{k})$ für jeden Divisor $D$ auf $X$ nach \cite[Proposition 7.3.7]{L1}. Aufgrund der Definition des Grades eines Geradenbündels und der Tatsache, dass alle zu betrachtenden Konstruktionen mit Basiswechsel vertauschen, folgt daraus die Behauptung.\hfill$\Box$\end{enumerate}}\par

\text{}\\

Für uns wichtig ist der Begriff des Grades eines Torseurs vor allem wegen des folgenden Resultats von Holla und Narasimham:\\\par

\begin{sat}[\protect{\cite[Theorem 1.2]{HN}}]
Es seien $G$ eine zusammenhängende reduktive algebraische Gruppe und $X$ eine glatte und projektive Kurve jeweils über einem algebraisch abgeschlossenen Körper $k$ beliebiger Charakteristik. Dann ist die Menge der Isomorphieklassen der semistabilen $G$"=Torseure auf $X$ mit fixiertem Grad beschränkt.\par
Dabei heiße eine Menge $\mathcal{S}$ von $G$"=Torseuren beschränkt, falls es ein Schema $S$ von endlichem Typ über $k$ und eine Familie von $G$"=Torseuren gibt, die durch $S$ parametrisiert wird, so dass jedes Element von $\mathcal{S}$ auf $X$ isomorph zu dem $G$"=Torseur auf $X$ ist, den man durch Einschränken der gegebenen Familie auf einen geeigneten abgeschlossenen Punkt von $S$ erhält.
\end{sat}

Weil nach einem Theorem von Lang (vgl. \cite[Theorem 1.3]{La} bzw. \cite{Sp}) \mbox{$H^{1}(k,G_{\overline{k}})=0$} für eine zusammenhängende reduktive algebraische Gruppe $G$ über einem endlichen Körper $k$ gilt, wobei $\overline{k}$ ein algebraischer Abschluss von $k$ ist, und es damit keine nichttrivialen Formen gibt, folgt daraus mittels noetherschem Descent:\\\par

\begin{kor}\label{isoklassen}
Es seien $G$ eine zusammenhängende reduktive algebraische Gruppe und $X$ eine glatte und projektive Kurve jeweils über einem endlichen Körper $k$. Dann ist die Menge der Isomorphieklassen der semistabilen $G$"=Torseure auf $X$ mit fixiertem Grad beschränkt ist.\\
\end{kor}\par

\begin{kor}\label{endliso}
Es seien $G$ eine zusammenhängende reduktive algebraische Gruppe und $X$ eine glatte und projektive Kurve jeweils über einem endlichen Körper $k$. Dann gibt es nur endlich viele Isomorphieklassen semistabiler $G$"=Torseuren auf $X$ mit fixiertem Grad.
\end{kor}
\bewub{Nach Korollar \ref{isoklassen} ist die Menge $\mathcal{S}$ der Isomorphieklassen der semistabilen $G$"=Torseure auf $X$ mit fixiertem Grad beschränkt. Per Definition gibt es also ein Schema $S$ von endlichem Typ über $k$ und eine Familie von $G$"=Torseuren auf $X$, die durch $S$ parametrisiert wird, so dass jedes Element der Menge der Isomorphieklassen der semistabilen $G$"=Torseure auf $X$ isomorph zu dem $G$"=Torseur auf $X$ ist, den man durch Einschränken der gegebenen Familie auf einen geeigneten abgeschlossenen Punkt von $S$ erhält. Ist $S$ affin, so ist $S$ ein Unterschema der affinen Geraden über $k$ und besitzt daher nur endlich viele abgeschlossene Punkte, so dass die Behauptung folgt. Ist $S$ nicht affin, so wähle man eine endliche affine Überdeckung von $S$ (eine solche existiert, da $S$ quasikompakt ist) und folgere dann wie im affinen Fall. \hfill$\Box$}\\\par
\newpage

\section{Zusammenhänge zwischen Semistabilität und Trivialisierbarkeit}

Ähnlich wie im von Lange und Stuhler in \cite{LS} behandelten Fall der Vektorbündel ergeben sich auch im Fall von $G$"=Torseuren Zusammenhänge zwischen Semistabilität und Trivialisierbarkeit der Torseure. Dazu definieren wir:\\\par

\begin{defn}
Es sei $G$ eine zusammenhängende reduktive algebraische Gruppe von endlicher Präsentation über einem vollkommenen Körper $k$ der Charakteristik $p>0$ und $X$ eine glatte und projektive Kurve über $k$. Dann heißt ein $G$"=Torseur $E$ auf $X$ \emph{streng semistabil}, falls $E$ semistabil ist und zudem für alle $r\in\nat$ das Pullback $(F^{r})^{*}E$ unter dem $r$"=fachen Produkt des absoluten Frobenius $F$ auf $X$ ebenfalls semistabil ist.\\
\end{defn}\par

Diese Definition erlaubt es nun, das folgende Resultat zu zeigen, das ein Analogon zu dem bekannten Resultat von Lange und Stuhler im Falle von Vektorbündeln anstelle von Prinzipalbündeln darstellt (siehe \cite[Satz 1.9.]{LS}):\\\par

\begin{sat}\label{ls19} 
Es seien $X$ eine glatte und projektive Kurve vom Geschlecht $g$ über einem endlichen Körper $k$ und $G$ eine zusammenhängende reduktive algebraische Gruppe von endlicher Präsentation über $k$.
Dann ist $E$, aufgefasst als $G$"=Torseur für die étale Topologie, genau dann trivialisierbar, wenn $E$ streng semistabil vom Grad Null ist.
\end{sat}  
\bewub{Ist $E$ durch eine endliche Überlagerung \mbox{$f\colon Y\to X$} trivialisierbar, so ist auch für jedes $r\in\nat$ das Pullback $(F^{r})^{*}E$ unter der $r$"=ten Potenz des absoluten Frobenius auf $X$ durch $f$ trivialisierbar. Nach Lemma \ref{endlichsemist} ist damit $(F^{r})^{*}E$ semistabil, da nach Korollar \ref{torsortrivial} der triviale Torseur $G_{Y}$ auf $Y$ semistabil ist. Also ist $E$ streng semistabil. Aus Bemerkung \ref{null} ergibt sich ferner, dass $d_{E}\equiv0$ ist.\\
Ist umgekehrt $E$ streng semistabil vom Grad Null, so ist \mbox{$(F^{r})^{*}(E)$} semistabil vom Grad Null für alle $r\in\nat$. Da nach Korollar \ref{endliso} die Anzahl der Isomorphieklassen semistabiler $G$"=Torseure $\widetilde{E}$ auf $X$ mit fest gewähltem Grad $d_{\widetilde{E}}$ endlich ist, gibt es dann ein $s\in\nat$, so dass man Isomorphismen von $G$"=Torseuren \mbox{$(F^{r})^{*}(E)\cong(F^{s})^{*}(E)$} für alle $r\geq s$ hat. Nach Korollar \ref{potenz} ist also $E$ trivialisierbar. \hfill$\Box$}\\\par

\section{Anwendung auf die Charakterisierung der Kategorie $\mathscr{B}_{\mathfrak{X}_{\mathfrak{o}},D}(G)$}

Es sei nun im folgenden wieder $G$ stets ein zusammenhängendes reduktives Gruppenschema von endlicher Präsentation über $\mathfrak{o}$ und $k=\overline{\mathbb{F}_{p}}$ der Restklassenkörper von $\overline{\mathbb{Z}}_{p}$.\par
Dann können wir nun das erste Hauptresultat dieser Arbeit zeigen:\\\par

\begin{thm}\label{thm20}
Es sei $\mathfrak{X}$ ein glattes Modell über $\overline{\mathbb{Z}}_{p}$ einer glatten und projektiven Kurve $X$ über $\overline{\mathbb{Q}}_{p}$ von von Null verschiedenem Geschlecht. Dann liegt ein $G$"=Torseur $P$ auf $\mathfrak{X}$ genau dann in der Kategorie $\mathscr{B}_{\mathfrak{X}_{\mathfrak{o}}}(G)$, wenn $P_{k}$ streng semistabil vom Grad Null auf der glatten und projektiven Kurve $\mathfrak{X}_{k}$ über $k$ ist. (Für die Tatsache, dass $\mathfrak{X}_{k}$ eine glatte und projektive Kurve ist, siehe u.\,a. Korollar \ref{irred}, \cite[Lemma 4.3.7]{L1}, \cite[Corollary 4.3.14]{L1}.) 
\end{thm}
\bewub{\begin{itemize}
\item[(i)] Liegt $P$ in $\mathscr{B}_{\mathfrak{X}_{\mathfrak{o}}}(G)$, so existiert nach Theorem \ref{thm16} eine Überdeckung $(\pi\colon\mathcal{Y}\to\mathfrak{X})\in S_{\mathfrak{X}}^{ss}$, so dass $\pi_{k}^{*}P_{k}$ trivial ist. Da $\pi_{k}$ endlich ist, ist nach Korollar \ref{torsortrivial} und Lemma \ref{endlichsemist} also $P_{k}$ ein semistabiler $G$"=Torseur auf $\mathfrak{X}_{k}$. Wegen $\pi\circ F_{\mathfrak{X}_{k}}=F_{\mathcal{Y}_{k}}\circ\pi$ folgt genauso, dass für alle \mbox{$r\in\nat$} auch $(F_{\mathfrak{X}_{k}}^{r})^{*}P_{k}$ semistabil ist, so dass also $P_{k}$ streng semistabil ist. Da der Grad des trivialen Torseurs der Nullmorphismus ist und nach Bemerkung \ref{null}$(ii)$ $d_{\pi_{k}^{*}P_{k}}=[K(\mathcal{Y}_{k}):K(\mathfrak{X}_{k})] d_{P_{k}}$ gilt, ist zudem der Grad von $P_{k}$ ebenfalls der Nullmorphismus.\par
\item[(ii)] Sei nun umgekehrt $P_{k}$ streng semistabil vom Grad Null. Da unter den gegebenen Voraussetzungen jeder $G$"=Torseur auch ein Torseur bezüglich der étalen Topologie ist, können wir im folgenden $P$ und $P_{k}$ stets als Torseure für die étale Topologie auffassen.\par
Mittels noetherschen Descents sieht man, dass die Familie $(X,\mathfrak{X},G,P_{k})$ zu einer Familie $(X_{K},\mathfrak{X}_{\mathfrak{o}_{K}},G_{0},P_{0})$ über einer endlichen Erweiterung $K$ von $\overline{\mathbb{Q}}_{p}$ mit Restklassenkörper $\kappa\cong\mathbb{F}_{q},q=p^{r},$ absteigt. Dabei ist $X_{K}$ eine glatte und projektive Kurve über $K$, $\mathfrak{X}_{\mathfrak{o}_{K}}$ ein glattes Modell von $X_{K}$ über $\mathfrak{o}_{K}$, $G_{0}$ eine zusammenhängende reduktive Gruppe von endlicher Präsentation über einer normalen, endlich erzeugten $\mathfrak{o}_{K}$"=Algebra $A$ und $P_{0}$ ein étaler $G$"=Torseur auf der speziellen Faser $\mathfrak{X}_{0}$ von $\mathfrak{X}_{\mathfrak{o}_{K}}$. Dies zeigt man genauso wie im Beweis von Theorem \ref{thm16}, bis auf die Änderung, dass hier $P_{k}$ als étaler $G$"=Torseur aufgefasst wird, und man beachte, dass die Eigenschaft eines Morphismus, "`étale"' zu sein, nach \cite[Proposition 17.7.8]{EGA IV} stabil (d.\,h. nach Übergang zu einer geeigneten Erweiterung) unter noetherschem Descent ist, ebenso wie nach \cite[Corollaire 8.4.5]{EGA IV} und \cite[Exposé XIV, Remarque 2.9]{SGA3} die Eigenschaften "`zusammenhängend"' und "`reduktiv"' des gegebenen Gruppenschemas $G$.\par
Man überlegt nun leicht, dass der abgestiegene $G_{0}$"=Torseur streng semistabil vom Grad Null ist: Nach Lemma \ref{noether} ist $P_{k}$ semistabil und wegen $(F^{r}_{\mathcal{X}_{k}})^{*}P_{k}=((F_{\mathfrak{X}_{0}}^{r})^{*}P_{0})\otimes_{k}\overline{k}$ folgt aus demselben Resultat, dass auch $(F_{\mathfrak{X}_{0}}^{r})^{*}P_{0}$ für alle $r\in\nat$ semistabil ist, so dass damit $P_{0}$ streng semistabil ist. Dass der Grad von $P_{0}$ Null ist, folgt ferner aus \ref{null}$(iii)$.\par
Also existiert nach Satz \ref{ls19} ein endlicher surjektiver Morphismus \linebreak \mbox{$\varphi\colon\mathcal{Y}_{0}\to\mathfrak{X}_{0}$}, so dass $\varphi^{*}P_{0}$ trivial ist. Nach der Konstruktion des Morphismus $\varphi$ im Beweis von Satz \ref{ls19} kann man dabei $\varphi$ als die Komposition eines endlichen étalen Morphismus $\phi_{0}\colon\mathcal{Y}_{0}\to\mathfrak{X}_{0}$, wobei $\mathcal{Y}_{0}$ eine glatte und projektive Kurve über $\kappa$ ist, und der Frobenius"=Potenz $F^{s}$ für ein geeignetes $s\geq 0$ wählen.\par 
Ab dieser Stelle verläuft der weitere Beweis wortwörtlich wie der von \cite[Theorem 20]{DW1}, man ersetze lediglich $\mathscr{E}_{k}$ durch $P_{k}$ und den Verweis auf \cite[Theorem 16]{DW1} durch den auf Theorem \ref{thm16}. \hfill$\Box$\\\end{itemize}}\par

Für den allgemeinen Fall definieren wir ähnlich wie in \cite{DW1}:\\\par

\begin{defn}
Es sei $R$ ein Bewertungsring mit Quotientenkörper $Q$ und Restklassenkörper $k$. Man betrachte ein Modell $\mathfrak{X}$ über $R$ einer glatten und projektiven Kurve $X$ über $Q$ und einen $G$"=Torseur $E$ auf $\mathfrak{X}$, wobei $G$ ein zusammenhängendes reduktives Grupppenschema von endlicher Präsentation über $R$ ist.\par
Dann sagt man, dass $E$ streng semistabile Reduktion vom Grad Null habe, falls das Pullback von $E_{k}$ zu der Normalisierung $\widetilde{C}$ jeder irreduziblen Komponente $C$ (versehen mit der reduzierten Struktur) von $\mathfrak{X}_{k}$ streng semistabil vom Grad Null ist. Man beachte dabei, dass jedes $\widetilde{C}$ eine glatte und projektive Kurve über $k$ ist. 
\end{defn}

Außerdem benötigen wir die folgende Verallgemeinerung von \cite[Theorem 18]{DW1}:\\\par

\begin{sat}\label{thm18}
Es sei $G$ ein zusammenhängendes reduktives Gruppenschema von endlicher Präsentation über $\mathbb{F}_{q}$, $X$ ein rein eindimensionales eigentliches Schema über $\mathbb{F}_{q}$ und $E$ ein $G$"=Torseur auf $X$.\par
Dann sind äquivalent:\par
\begin{itemize}
\item[(i)] Das Pullback von $E$ zu der Normalisierung jeder irreduziblen Komponente von $X$ ist streng semistabil vom Grad Null.
\item[(ii)] Es gibt einen endlichen surjektiven Morphismus $\varphi\colon Y\to X$, wobei $Y$ ein rein eindimensionales eigentliches Schema über $\mathbb{F}_{q}$ ist, so dass $\varphi^{*}E$ ein trivialer $G$"=Torseur ist.
\item[(iii)] Genauso wie in $(ii)$, aber mit $\varphi$ als Komposition $\varphi\colon Y\xrightarrow{F^{s}}Y\xrightarrow{\pi}X$ für ein $s\geq 0$, wobei $\pi$ endlich, étale und surjektiv ist und $F=Fr_{q}=Fr_{p}^{r}$ der $q=p^{r}$"=lineare Frobenius auf $Y$ ist. 
\end{itemize}
\end{sat}
\bewub{\begin{itemize}{
\item[(iii)$\Rightarrow$(ii)] Dies ist trivial.
\item[(ii)$\Rightarrow$(i)] Sei $(ii)$ gegeben. Dann wird jede irreduzible Komponente $C$ von $X$ endlich dominiert durch eine irreduzible Komponente $D$ von $Y$. Es folgt, dass das Pullback von $E$ zu $\widetilde{C}$ durch den endlichen Morphismus $\widetilde{D}\to\widetilde{C}$ trivialisiert wird. Da man nach Lemma \ref{endlichsemist} Semistabilität nach Pullback entlang einer endlichen Überdeckung prüfen kann und da der absolute Frobenius funktoriell ist, folgt $(i)$.
\item[(i)$\Rightarrow$(iii)] Nach Korollar \ref{endliso} gibt es nur endlich viele Isomorphieklassen semistabiler $G$"=Torseure mit fixiertem Grad jeweils auf den endlich vielen einzelnen irreduziblen Komponenten von $X$. Als ersten Schritt folgern wir daraus, dass es nur endlich viele Isomorphieklassen étaler $G$"=Torseure auf $X$ gibt, deren Pullbacks zu den Normalisierungen der irreduziblen Komponenten von $X$ jeweils semistabil vom Grad Null sind.\par
Dazu nehmen wir zunächst an, dass $X$ reduziert sei. Dann ist der Normalisierungsmorphismus $\pi\colon\widetilde{X}=\coprod\widetilde{C}_{v}\to X$ ein endlicher Morphismus, wobei $X=\cup_{v}C_{v}$ die Zerlegung von $X$ in seine irreduziblen Komponenten sei. \par
Ferner ist der natürliche Morphismus $\alpha\colon G_{X}\to\pi_{*}G_{\widetilde{X}}$ injektiv: Der Normalisierungsmorphismus $\pi$ ist endlich und surjektiv und daher eine Überdeckung für die $h$"=Topologie auf $X$ (vgl. \cite[Definition 3.1.2]{V}), insbesondere also ein universeller topologischer Epimorphismus. Da jedes Schema $U$, das étale über $X$ ist, reduziert ist, ist also nach \cite[Lemma 3.2.1]{V} der Morphismus $\pi_{U}\colon U\times_{X}\widetilde{X}\to U$ ein kategorieller Epimorphismus in der Kategorie der Schemata, so dass die Abbildung \[\alpha(U)\colon\Gamma(U,G_{U})\to\Gamma(U\times_{X}\widetilde{X},G_{U\times_{X}\widetilde{X}})\] injektiv ist.\par
Da der Träger des Kokerns von $\alpha\colon G_{X}\to\pi_{*}G_{\widetilde{X}}$ höchstens aus den singulären Punkten von $X$ besteht, ist dieser Kokern eine Wolkenkratzergarbe von Mengen, also nach \cite[Exercise 2.2.9]{L1} von der Gestalt $\prod_{x\in X^{sing}}i_{x*}S_{x}$. Wir behaupten nun, dass jede der Mengen $S_{x}$ endlich ist:\par
Dazu zeigen wir, dass wir die folgende exakte Sequenz von étalen Garben auf $X$ haben:
\[0\to G_{X} \to \pi_{*}G_{\widetilde{X}} \to \bigoplus_{x\in X^{sing}}i_{x*}(\bigoplus_{y\mathop{|}x}G(\kappa(y))/G(\kappa(x))). \] Es sei $i\colon x\to X$ ein abgeschlossener Punkt. Dann hat man das kartesische Diagramm \begin{equation*}
        \xymatrix@=3em{%
         x\times_{X}\widetilde{X} \ar[r]^{\widetilde{i_{x}}} \ar[d]^{\pi_{x}} & \widetilde{X} \ar[d]^{\pi} \\
         x \ar@{^{(}->}[r]_{i_{x}} & X}
           \end{equation*} und nach \cite[Satz II.6.4.2]{Ta} existiert ein kanonischer Basiswechselmorphismus \[i_{x}^{*}\pi_{*}G_{\widetilde{X}}\stackrel{\sim}{\longrightarrow} \pi_{x*}\widetilde{i_{x}}^{*}G_{\widetilde{X}},\] der ein Isomorphismus ist.\par 
           Daraus folgt \[(\pi_{*}G_{\widetilde{X}})_{x}=(\pi_{x*}\widetilde{i_{x}}^{*}G_{\widetilde{X}})_{x}\] und nach \cite[Corollary II.3.5]{M} \[(\pi_{*}G_{\widetilde{X}})_{x}=\bigoplus_{y\mathop{|}x}G_{\widetilde{X},y}.\] Nach \cite[Remark II.2.9(d)]{M} ist ferner $G_{X,x}=G(\kappa(x))$ und \mbox{$G_{\widetilde{X},y}=G(\kappa(y))$}. Da sich der Halm von $\bigoplus_{x\in X^{sing}}i_{x*}(\bigoplus_{y\mathop{|}x}G(\kappa(y))/G(\kappa(x)))$ im Punkt $x$ zu $\prod_{\widetilde{x}\mathop{|}x}(G(\kappa(\widetilde{x})/G(\kappa(x))$ berechnet, ist also die Sequenz \[0\to G_{X} \to \pi_{*}G_{\widetilde{X}} \to \bigoplus_{x\in X^{sing}}i_{x*}(\bigoplus_{y\mathop{|}x}G(\kappa(y))/G(\kappa(x)))\] auf jedem Halm exakt und damit auch selbst exakt. Insbesondere folgt damit $S_{x}\subset\prod_{\widetilde{x}\mathop{|}x}(G(\kappa(\widetilde{x}))/G(\kappa(x)))$ für alle $x\in X^{sing}$. Aus der Endlichkeit von $\kappa(\widetilde{x})$ und der Endlichkeit des Produkts folgt also, dass $S_{x}$ für alle $x\in X^{Sing}$ endlich ist.\par
Wie im Fall der Vektorbündel folgt also aus \cite[Corollaire III.3.2.4]{G}, dass es nur endlich viele Isomorphieklassen von $G$"=Torseuren auf $X$ gibt, die bestimmte vorgegebene Isomorphieklassen von $G$"=Torseuren auf den Kurven $\widetilde{C}_{v}$ induzieren. Insbesondere gibt es also nur endlich viele Isomorphieklassen von $G$"=Torseuren auf $X$, deren Pullback zu den Normalisierungen der irreduziblen Komponenten von $X$ semistabil vom Grad Null ist.\par           
Für den allgemeinen Fall, dass $X$ nicht notwendigerweise reduziert ist, genügt es daher zu zeigen, dass die kanonische Abbildung \[\varphi\colon H^{1}_{\acute{e}t}(X,G)\longrightarrow H^{1}_{\acute{e}t}(X_{red},G)\] endliche Fasern besitzt. Mittels Devissage sieht man, dass es zu zeigen genügt, dass für jedes Ideal $\mathcal{J}\subset\mathcal{O}_{X}$ mit $\mathcal{J}^{2}=0$ die Abbildung \[i^{*}\colon H^{1}_{\acute{e}t}(X,G)\longrightarrow H^{1}_{\acute{e}t}(X',G)\] endliche Fasern besitzt, wobei $i\colon X'\into X$ das durch $\mathcal{J}$ definierte abgeschlossene Unterschema von $X$ ist.\par
Ähnlich wie im Beweis von \ref{thm16} hat man nach \cite[Lemme VII.1.3.5]{G} die folgende exakte Sequenz von Garben auf $X_{\acute{e}t}$: \begin{equation*}
        \xymatrix@=1em{%
         1 \ar[r] & i_{*}(Lie(G_{X'}/X')\otimes_{\mathcal{O}_{X}}\mathcal{J}) \ar[rrrr]^-{h} & & & & G \ar[rrr]^-{adj}& & & i_{*}G_{X'} \ar[r] & 1.}
           \end{equation*}
Man beachte dabei, dass $i$ ein universeller Homöomorphismus ist, so dass nach \cite[Remark II.3.17]{M} $\widetilde{X'}_{\acute{e}t}\cong\widetilde{X}_{\acute{e}t}$ gilt und man $\mathcal{J}$ kanonisch mit einem $\mathcal{O}_{X'}$"=Modul identifizieren kann. Desweiteren bezeichnet $\adj$ die Komposition des Adjunktionsmorphismus $G_{X}\to i_{*}i^{*}G_{X}$ mit dem Morphismus $i_{*}i^{*}G_{X}\to i_{*}G_{X'}$, der von dem kanonischen Morphismus $i^{*}G_{X}\to G_{X'}$ induziert wird. 
Der Morphismus $h$ schließlich wird durch die Bijektion \[i_{*}(Lie(G_{X'}/X')\otimes_{\mathcal{O}_{X'}}\mathcal{J})(S)\stackrel{\sim}{\longrightarrow}\ker(G_{X}\to i_{*}G_{X'})(S),\alpha\to e\alpha\] induziert, wobei $e$ der Einsschnitt von $G_{X'}$ ist und $S$ étale über $X$ ist.\par        
Aus dieser exakten Sequenz erhält man die exakte Kohomologiesequenz \[H^{1}_{\acute{e}t}(X,i_{*}(Lie(G_{X'}/X')\otimes_{\mathcal{O}_{X'}}\mathcal{J}))\xrightarrow{h} H^{1}_{\acute{e}t}(X,G_{X})\to H^{1}_{\acute{e}t}(X,i_{*}G_{X'}).\]
Nach \cite[Corollaire VIII.5.6]{SGA 4} ist nun $R^{1}i_{*}G_{X'}=0$, da $i$ als abgeschlossene Immersion insbesondere ganz ist.\par
Nach \cite[Proposition V.3.1.3]{G} ist damit $H^{1}_{\acute{e}t}(X,i_{*}G_{X'})$ isomorph zu $H^{1}_{\acute{e}t}(X',G_{X'})$, so dass also die Sequenz  \[H^{1}_{\acute{e}t}(X,i_{*}(Lie(G_{X'}/X')\otimes_{\mathcal{O}_{X'}}\mathcal{J}))\xrightarrow{h} H^{1}_{\acute{e}t}(X,G_{X})\xrightarrow{i^{*}} H^{1}_{\acute{e}t}(X',G_{X'})\] exakt ist.\par
Wie oben gilt nun nach \cite[Corollaire VIII.5.6]{SGA 4} \[R^{1}i_{*}(Lie(G_{X'}/X')\otimes_{\mathcal{O}_{X'}}\mathcal{J})=0.\] Nach \cite[Proposition V.3.1.3]{G} hat man somit einen Isomorphismus \begin{equation*}
 \xymatrix@=3em{%
H^{1}_{\acute{e}t}(X',Lie(G_{X'}/X')\otimes_{\mathcal{O}_{X'}}\mathcal{J})\ar[r]^{j}_{\sim} & H^{1}_{\acute{e}t}(X,i_{*}(Lie(G_{X'}/X')\otimes_{\mathcal{O}_{X'}}\mathcal{J})),}\end{equation*} so dass also die Sequenz  \[H^{1}_{\acute{e}t}(X',Lie(G_{X'}/X')\otimes_{\mathcal{O}_{X'}}\mathcal{J})\stackrel{h\circ j}{\longrightarrow}H^{1}_{\acute{e}t}(X,G_{X})\xrightarrow{i^{*}} H^{1}_{\acute{e}t}(X',G_{X'})\] exakt ist.\par
Die kanonische $\mathcal{O}_{X'}$"=Modulstruktur von $Lie(G_{X'}/X')$ induziert ferner die folgende kanonische exakte Sequenz von Garben auf $X'_{\acute{e}t}$:
\[0\longrightarrow\ker(\alpha)\longrightarrow\mathcal{J}\stackrel{\alpha}{\longrightarrow}Lie(G_{X'}/X')\otimes_{\mathcal{O}_{X'}}\mathcal{J}\longrightarrow0,\] wobei $\alpha$ der kanonische Morphismus ist. Dies liefert eine Surjektion \begin{equation*}\xymatrix@=2em{%
H^{1}_{\acute{e}t}(X',\mathcal{J})\ar@{>>}[rrr]^-{\alpha} & & & H^{1}_{\acute{e}t}(X',Lie(G_{X'}/X')\otimes_{\mathcal{O}_{X'}}\mathcal{J}),}\end{equation*} da $X'$ höchstens eindimensional ist und nach \cite[Proposition III.3.3]{M} \[H^{2}_{\acute{e}t}(X',\mathscr{F})\cong H^{2}(X',\mathscr{F}_{Zar})\] für jede étale Garbe $\mathscr{F}$ auf $X'$ gilt, wobei $\mathscr{F}_{Zar}$ die Einschränkung auf die Zariski-Topologie ist.
Damit erhält man die folgende nicht"=abelsche Kohomologiesequenz: \[H^{1}_{\acute{e}t}(X',\mathcal{J})\xrightarrow{h\circ j\circ\alpha} H^{1}_{\acute{e}t}(X,G_{X})\xrightarrow{i^{*}} H^{1}_{\acute{e}t}(X',G_{X'}).\]
Da $H^{1}_{\acute{e}t}(X',\mathcal{J})$ ein endlich"=dimensionaler $\mathbb{F}_{q}$"=Vektorraum und damit endlich ist, folgt daraus, dass $\varphi$ endliche Fasern besitzt.\par
Damit haben wir gesehen, dass es nur endlich viele Isomorphieklassen von $G$"=Torseuren auf $X$ gibt, deren Pullbacks zu den Normalisierungen der irreduziblen Komponenten von $X$ semistabil vom Grad Null sind.\par
Es sei nun ein $G$"=Torseur $E$ wie in $(i)$ gegeben. Dann sind die Pullbacks zu $\widetilde{C}_{v}$ aller $G$"=Torseure $(F_{X}^{n})^{*}E$ auf $X$ jeweils semistabil vom Grad Null. Da es nur endlich viele Isomorphieklassen derartiger $G$"=Torseure auf $X$ gibt, existiert ein $s\geq 0$, so dass $(F_{X}^{t})^{*}E\cong(F^{s}_{X})^{*}E$ für alle $t\geq s$ gilt. Für den $G$"=Torseur $E':=(F_{X}^{s})^{*}E$ gilt damit $(F^{r}_{X})^{*}E'=E'$ für alle $r=t-s\geq1$. Fassen wir nun $E$ und $E'$ als Torseure für die étale Topologie auf, so existiert folglich nach Korollar \ref{pullback2} ein endlicher, étaler und surjektiver Morphismus $\pi\colon Y\to X$, so dass \[\pi^{*}E'=\pi^{*}(F_{X}^{s})^{*}E\] ein trivialer $G$"=Torseur auf $Y$ ist. Es folgt daraus, dass \[(\pi\circ F_{Y}^{s})^{*}E=(F^{s}_{X}\circ\pi)^{*}E=\pi^{*}(F^{s}_{X})^{*}E\] ein trivialer étaler $G$"=Torseur auf $Y$ ist.\par
Weil ferner $X$ ein rein eindimensionales eigentliches $\mathbb{F}_{q}$"=Schema ist, so dass damit auch $Y$ ein rein eindimensionales eigentliches $\mathbb{F}_{q}$"=Schema ist, ist damit $(iii)$ gezeigt.\hfill$\Box$}\\\end{itemize}}\par

Als Verallgemeinerung von \cite[Theorem 17]{DW1} können wir damit zeigen:\\\par

\begin{thm}\label{thm17}
Es sei $\mathfrak{X}$ ein Modell über $\overline{\mathbb{Z}}_{p}$ einer glatten und projektiven Kurve $X$ über $\overline{\mathbb{Q}}_{p}$ sowie $G$ ein zusammenhängendes reduktives Gruppenschema von endlicher Präsentation über $\mathfrak{o}$.\par Dann liegt ein $G$"=Torseur $P$ auf $\mathfrak{X}_{\mathfrak{o}}$ genau dann in der Kategorie $\mathscr{B}_{\mathfrak{X}_{\mathfrak{o}},D}(G)$ für einen gewissen Divisor $D$, wenn $P$ streng semistabile Reduktion vom Grad Null hat. Ist dies der Fall, so gibt es zwei Divisoren $D$ und $\widetilde{D}$ auf $X$ mit disjunktem Träger, so dass $P$ sowohl in der Kategorie $\mathscr{B}_{\mathfrak{X}_{\mathfrak{o}},D}(G)$ als auch in der Kategorie $\mathscr{B}_{\mathfrak{X}_{\mathfrak{o}},\widetilde{D}}(G)$ liegt.
\end{thm}
\bewub{Liegt der $G$"=Torseur $P$ auf $\mathfrak{X}_{\mathfrak{o}}$ in der Kategorie $\mathscr{B}_{\mathfrak{X}_{\mathfrak{o}},D}(G)$ für einen gewissen Divisor $D$ auf $X$, so gibt es nach Theorem \ref{thm16} eine Überdeckung $\pi\colon\mathcal{Y}\to\mathfrak{X}$ aus $S_{\mathfrak{X},D}^{good}$, so dass $\pi_{k}^{*}P_{k}$ ein trivialer $G$"=Torseur ist.\par
Es sei $\mathfrak{X}_{k}=\cup_{v}C_{v}$ die Zerlegung von $\mathfrak{X}_{k}$ in seine irreduziblen Komponenten. Da $\mathfrak{X}$ irreduzibel ist (vgl. die Argumentation in §1.2), $\pi(\mathcal{Y})$ nach \cite[Proposition 6.1.10]{EGA II} abgeschlossen ist und den generischen Punkt von $\mathfrak{X}$ erhält, ist $\pi$ surjektiv. Also wird jedes $C_{v}$ durch eine irreduzible Komponente $\widetilde{C}_{v}$ von $\mathcal{Y}_{k}$ endlich dominiert. Daraus folgt mit Korollar \ref{torsortrivial} und Lemma \ref{endlichsemist}, dass die Pullbacks von $P_{k}$ zu allen $\widetilde{C}_{v}$ jeweils streng semistabil vom Grad Null sind.\par
Es habe nun umgekehrt der $G$"=Torseur $P$ auf $\mathfrak{X}_{\mathfrak{o}}$ streng semistabile Reduktion vom Grad Null und er werde im folgenden stets als étaler $G$"=Torseur aufgefasst. Dann sieht man ähnlich wie z.B. in den Beweisen der Theoreme \ref{thm16} und \ref{thm20} mittels noetherschen Descents, dass es eine endliche Erweiterung $K$ von $\mathbb{Q}_{p}$ mit Ganzheitsring $\mathfrak{o}_{K}$ und Restklassenkörper $\kappa\cong\mathbb{F}_{q}$ gibt, so dass die Familie $(X,\mathfrak{X},C_{v},G,P_{k})$ zu einer Familie $(X_{K},\mathfrak{X}_{\mathfrak{o}_{K}},C_{v0},G_{0},P_{0})$ mit den entsprechenden Eigenschaften absteigt. Insbesondere ist $P_{0}$ ein étaler $G_{0}$"=Torseur auf der speziellen Faser $\mathfrak{X}_{0}=\mathfrak{X}_{\mathfrak{o}_{K}}\otimes\kappa$, dessen Pullbacks zu den Normalisierungen $\widetilde{C}_{v0}$ der irreduziblen Komponenten $C_{v0}$ von $\mathfrak{X}_{0}$ streng semistabil vom Grad Null sind. 
Nach Satz \ref{thm18} erhält man dann einen endlichen étalen Morphismus $\widetilde{\pi_{0}}\colon\widetilde{\mathcal{Y}_{0}}\to\mathfrak{X}_{0}$, so dass für die Komposition $\widetilde{\varphi_{0}}\colon\widetilde{\mathcal{Y}_{0}} \xrightarrow{F^{s}}\widetilde{\mathcal{Y}_{0}}\xrightarrow{\widetilde{\pi_{0}}}\mathfrak{X}_{0}$ das Pullback $\widetilde{\varphi_{0}}^{*}P_{0}$ trivial ist. Man beachte dabei, dass man bei dieser Aussage $s$ durch jede beliebige größere ganze Zahl $s'$ und dann $F$ durch jede beliebige Potenz von $F$ ersetzen darf. Wie in \cite{DW1} erlaubt \cite[Exposé IX, Théorème 1.10]{SGA 1} es, den Morphismus $\widetilde{\pi_{0}}$ zu einem endlichen étalen Morphismus $\widetilde{\pi_{\mathfrak{o}_{K}}}\colon\widetilde{\mathcal{Y}_{\mathfrak{o}_{K}}}\to\mathfrak{X}_{0}$ mit $\widetilde{\pi_{0}}$ als spezieller Faser zu liften. Ersetzt man nun $K$ wie in \cite{DW1} durch eine geeignete endliche Erweiterung $K_{1}$ von $K$, so kann man $\widetilde{\pi_{\mathfrak{o}_{K_{1}}}}$ durch ein Objekt $\pi_{\mathfrak{o}_{K_{1}}}\colon\mathcal{Y}_{\mathfrak{o}_{K_{1}}}\to\mathfrak{X}_{\mathfrak{o}_{K_{1}}}$ aus der Kategorie $S_{\mathfrak{X}_{\mathfrak{o}_{K_{1}}},\emptyset}$ dominieren, wobei man wie in \cite{DW1} annehmen darf, dass $\mathcal{Y}_{\mathfrak{o}_{K_{1}}}$ nicht nur semistabil, sondern auch regulär ist. Schreibt man zur Vereinfachung der Notation wieder $K$ statt $K_{1}$, $P_{0}$ statt $P_{1}$, $\mathfrak{o}_{K}$ statt $\mathfrak{o}_{K_{1}}$, $\kappa$ statt $\kappa_{1}$ etc., so folgt, dass das Pullback von $P_{0}$ unter der Komposition $\varphi\colon\mathcal{Y}_{0}\xrightarrow{F^{s}}\mathcal{Y}_{0}\xrightarrow{\pi_{0}}\mathfrak{X}_{0}$ ein trivialer $G_{0}$"=Torseur ist.\par
Ab hier überträgt sich der weitere Beweis in \cite[Theorem 17]{DW1} nun wortwörtlich, wobei man lediglich $\mathcal{E}_{\kappa}$ durch $P_{k}$, $\mathcal{E}$ durch $P$ und $\mathcal{E}_{0}$ durch $P_{0}$ sowie den Verweis auf \cite[Theorem 16]{DW1} durch den auf Theorem \ref{thm16} ersetze. \hfill$\Box$}\\\par
\newpage
\thispagestyle{empty}

\chapter{Prinzipalbündel auf $X_{\mathbb{C}_{p}}$ und Paralleltransport}
\section{Potentiell streng semistabile Reduktion und Paralleltransport}

Es sei weiterhin $X$ eine glatte und projektive Kurve über $\overline{\mathbb{Q}_{p}}$, $\mathfrak{X}$ ein Modell von $X$ über $\overline{\mathbb{Z}_{p}}$ und $G$ ein zusammenhängendes reduktives Gruppenschema von endlicher Präsentation über $\mathfrak{o}$.\par
Um die Resultate der vorherigen Kapitel auf den Fall von Prinzipalbündeln über $X_{\mathbb{C}_{p}}$ auszudehnen, benötigen wir die folgende Definition in Anlehnung an die entsprechende Definition im Falle von Vektorbündeln (siehe \cite[§0]{DW1}):\\\par

\begin{defn}
Man sagt, dass ein $G$"=Torseur $E$ auf $X_{\mathbb{C}_{p}}$ streng semistabile Reduktion vom Grad Null habe, wenn er die folgende Eigenschaft besitzt:\par
$E$ lässt sich zu einem $G$"=Torseur $\widetilde{E}$ auf $\mathfrak{X}_{\mathfrak{o}}$ für ein geeignetes Modell $\mathfrak{X}$ der Kurve $X$ ausdehnen und das Pullback der speziellen Faser $\widetilde{E}_{k}$ dieses $G$"=Torseurs $\widetilde{E}$ zu der Normalisierung jeder irreduziblen Komponente von $\mathfrak{X}_{k}$ ist streng semistabil vom Grad Null.\par
Ferner sagt man, dass $E$ potentiell streng semistabile Reduktion vom Grad Null besitze, wenn es einen endlichen étalen Morphismus $\alpha\colon Y\to X$ von glatten und projektiven Kurven über $\overline{\mathbb{Q}_{p}}$ gibt, so dass $\alpha^{*}_{\mathbb{C}_{p}}E$ streng semistabile Reduktion vom Grad Null hat.\\   
\end{defn}\par

Es sei nun $\mathcal{B}_{X_{\mathbb{C}_{p}}}^{s}(G)$ die Kategorie der $G$"=Torseure auf $X_{\mathbb{C}_{p}}$ mit streng semistabiler Reduktion vom Grad Null und $\mathcal{B}_{X_{\mathbb{C}_{p}}}^{ps}(G)$ die Kategorie der $G$"=Torseure auf $X_{\mathbb{C}_{p}}$ mit potentiell streng semistabiler Reduktion vom Grad Null.\par

Dann können wir nun die beiden Hauptresultate dieser Arbeit zeigen:\\\par

\begin{thm}\label{s}
Es gilt $\mathcal{B}^{s}_{X_{\mathbb{C}_{p}}}(G)=\bigcup_{D}\mathcal{B}_{X_{\mathbb{C}_{p}}, D}(G)$.\par
Jeder $G$"=Torseur in $\mathcal{B}^{s}_{X_{\mathbb{C}_{p}}}(G)$ liegt sowohl in $\mathcal{B}_{X_{\mathbb{C}_{p}}, D}(G)$ als auch in \linebreak $\mathcal{B}_{X_{\mathbb{C}_{p}}, \widetilde{D}}(G)$ für geeignete Divisoren $D$ und $\widetilde{D}$ mit disjunkten Träger. \par
Ferner gibt es für jeden Torseur $E\in\mathcal{B}^{s}_{X_{\mathbb{C}_{p}}, D}(G)$ einen eindeutig bestimmten stetigen Funktor $\rho_{E}$ von $\Pi_{1}(X)$ nach $\mathcal{P}(G(\mathbb{C}_{p}))$ mit $\rho_{E}(x)=E_{x}$ für alle \mbox{$x\in X(\mathbb{C}_{p})$}, so dass $\rho_{E}$ verträglich mit den Funktoren $\rho_{E}$ von $\Pi_{1}(X-D)$ nach $\mathcal{P}(G(\mathbb{C}_{p}))$ ist, die für die Divisoren $D$ mit $E\in\mathcal{B}_{X_{\mathbb{C}_{p}},D}(G)$ im zweiten Kapitel konstruiert wurden. \par
Insgesamt erhält man so einen Funktor $\rho\colon\mathcal{B}_{X_{\mathbb{C}_{p}}}^{s}(G)\to\Rep_{\Pi_{1}(X)}(G(\mathbb{C}_{p}))$, der sich funktoriell bezüglich Morphismen von Kurven über $\overline{\mathbb{Q}_{p}}$, Morphismen von zusammenhängenden reduktiven Gruppenschemata von endlicher Präsentation über $\mathfrak{o}$ und $\mathbb{Q}_{p}$"=Automorphismen von $\overline{\mathbb{Q}_{p}}$ verhält und mittels der in §2.4. dargestellten Funktoren verträglich mit den entsprechenden Funktoren im Vektorbündelfall ist.\par
Für jeden Punkt $x_{0}\in X(\mathbb{C}_{p})$ ist der "`Faserfunktor in $x_{0}$"' \[\mathcal{B}_{X_{\mathbb{C}_{p}}}^{s}(G)\to \mathcal{P}(G(\mathbb{C}_{p})), E\mapsto E_{x_{0}}, f\mapsto f_{x_{0}}\] ein treuer Funktor.
\end{thm}
\bewub{Ist $E\in\mathcal{B}_{X_{\mathbb{C}_{p}},D}(G)$ für einen Divisor $D$, so dehnt sich nach Definition der Kategorie $\mathcal{B}_{X_{\mathbb{C}_{p}},D}(G)$ und nach Theorem \ref{thm17} $E$ zu einem $G$"=Torseur $\widetilde{E}$ auf $\mathfrak{X}_{\mathfrak{o}}$ für ein geeignetes Modell $\mathfrak{X}$ von $X$ aus, der streng semistabile Reduktion vom Grad Null hat. Nach Definition ist daher $E\in\mathcal{B}^{s}_{X_{\mathbb{C}_{p}}}(G)$.\par Ist umgekehrt $E'\in\mathcal{B}^{s}_{X_{\mathbb{C}_{p}}}(G)$, so folgt ebenfalls aus Theorem \ref{thm17}, dass $E'\in\mathcal{B}_{X_{\mathbb{C}_{p}},D}(G)$ für einen geeigneten Divisor $D$ ist und es genauer Divisoren $D$ und $\widetilde{D}$ mit disjunktem Träger gibt, für die $E$ sowohl in $\mathcal{B}_{X_{\mathbb{C}_{p}}, D}(G)$ als auch in $\mathcal{B}_{X_{\mathbb{C}_{p}}, \widetilde{D}}(G)$ liegt. Insgesamt folgt also $\mathcal{B}^{s}_{X_{\mathbb{C}_{p}}}(G)=\bigcup_{D}\mathcal{B}_{X_{\mathbb{C}_{p}}, D}(G)$.\par
Ist $E$ ein Torseur in $\mathcal{B}_{X_{\mathbb{C}_{p}}}^{s}(G)$, so existiert für alle Divisoren $D$ mit \linebreak\mbox{$E\in\mathcal{B}_{X_{\mathbb{C}_{p}},D}(G)$}, wie in §2.2 gezeigt, ein stetiger Funktor \[\rho_{E,D}\colon\mathcal{B}_{X_{\mathbb{C}_{p}},D}(G)\to\Rep_{\Pi_{1}(X-D)}(G(\mathbb{C}_{p}))\] mit $\rho_{E}(x)=E_{x}$ für alle $x\in (X-D)(\mathbb{C}_{p})$. Nach Korollar \ref{seifert} folgt dann, dass es genau einen Funktor \[\rho_ {E}\colon\Pi_{1}(X)\to\mathcal{P}(G(\mathbb{C}_{p}))\] gibt, der auf $\Pi_{1}(X-D)$ jeweils die Funktoren $\rho_{E,D}$ induziert.\par
Es ist klar, dass durch die Zuordnung $E\mapsto\rho_{E}$ ein Funktor \[\rho\colon\mathcal{B}_{X_{\mathbb{C}_{p}}}^{s}(G)\to\Rep_{\Pi_{1}(X)}(G(\mathbb{C}_{p}))\] definiert wird.\par
Die Funktorialitätseigenschaften folgen direkt aus den Sätzen \ref{galoiskonj1}, \ref{unabh}, \ref{gruppen1}, \ref{darstellung1} und \ref{vb1}.\par
Es sei nun $x_{0}$ ein fest gewählter $\mathbb{C}_{p}$"=wertiger Punkt von $X$. Dann ist zu zeigen, dass die Abbildung \[\Hom_{\mathcal{B}_{X_{\mathbb{C}_{p}}}^{s}(G)}(E,E')\longrightarrow\Hom_{\mathcal{P}(G(\mathbb{C}_{p}))}(E_{x_{0}},E'_{x_{0}})\] für alle $G$"=Torseure $E$ und $E'$ in $\mathcal{B}_{X_{\mathbb{C}_{p}}}^{s}(G)$ injektiv ist.\par
Es seien also $E$ und $E'$ zwei beliebig gewählte $G$"=Torseure in $\mathcal{B}_{X_{\mathbb{C}_{p}}}^{s}(G)$ und \mbox{$f,g\colon E\longrightarrow E'$} zwei Morphismen von $G$"=Torseuren auf $X_{\mathbb{C}_{p}}$, so dass \mbox{$f_{x_{0}}=g_{x_{0}}$} gilt.\par
Ist nun \mbox{$x\in X(\mathbb{C}_{p})$} ein beliebiger $\mathbb{C}_{p}$"=wertiger Punkt von $X$, so wähle man einen Weg $\gamma$ von $x_{0}$ nach $x$. (Ein solcher existiert nach der Definition von $\Pi_{1}(X)$ (vgl. \cite{SGA 1}).) Dann existiert jeweils ein Paralleltransport $\rho_{E}(\gamma)\colon E_{x_{0}}\longrightarrow E_{x}$ bzw.  \mbox{$\rho_{E'}(\gamma)\colon E'_{x_{0}}\longrightarrow E'_{x}$}, so dass man aufgrund der Funktorialität der Konstruktion des Paralleltransports ein kommutatives Diagramm \begin{equation*}
        \xymatrix@=3em{%
         E_{x_{0}} \ar[r]^{\rho_{E}(\gamma)} \ar[d] & E_{x}\ar[d] \\ 
       E'_{x_{0}} \ar[r]^{\rho_{E'}(\gamma)} & E'_{x} }
           \end{equation*} sowohl dann erhält, wenn man als die senkrechten Pfeile $f_{x_{0}}$ und $f_{x}$ wählt, als auch, wenn man $g_{x_{0}}$ und $g_{x}$ wählt. Aufgrund der Kommutativität des Diagramms und der Tatache, dass alle horizontalen Abbildungen Isomorphismen sind, folgt dann, dass $f_{x}=g_{x}$ gilt. Da $x\in X(\mathbb{C}_{p})$ beliebig gewählt war, ist also $f_{x}=g_{x}$ für alle $\mathbb{C}_{p}$"=wertigen Punkte $x$ von $X$ und also insbesondere für alle $\mathbb{C}_{p}$"=wertigen Punkte von $E$. Damit stimmen auch die von $f$ bzw. $g$ induzierten Abbildungen $\widetilde{f},\widetilde{g}\colon E(\mathbb{C}_{p})\longrightarrow E'(\mathbb{C}_{p})$ überein.
           Weil $G_{X_{\mathbb{C}_{p}}}$ eine zusammenhängende reduktive algebraische Gruppe von endlicher Präsentation über $X_{\mathbb{C}_{p}}$ ist, sind $E$ und $E'$ affin, glatt, reduktiv und von endlicher Präsentation über $X_{\mathbb{C}_{p}}$ (dies folgt jeweils mit fpqc"=Descent wie in Lemma \ref{descent}) und damit auch nach \cite[Proposition 5.3.4]{EGA II} projektiv über $\mathbb{C}_{p}$. Daher kommen $E$ und $E'$ nach \cite[S. 50]{EH} von Varietäten $\mathcal{E}$ bzw. $\mathcal{E'}$ über $\mathbb{C}_{p}$ im Sinne der klassischen algebraischen Geometrie à la Weil. Da ein Morphismus $\mathcal{E}\to\mathcal{E'}$ eindeutig durch den Morphismus $\mathcal{E}(\mathbb{C}_{p})\to\mathcal{E'}(\mathbb{C}_{p})$ bestimmt ist, stimmen daher $f$ und $g$ auf der Menge der abgeschlossenen Punkte von $E$ überein, die nach \cite[Exercise II.3.14]{Ha} eine offene dichte Teilmenge von $E$ ist. Nach \cite[Lemma 7.2.2.1]{EGA I} stimmen also $f$ und $g$ auf ganz $E$ überein, so dass damit insgesamt der Faserfunktor in $x_{0}$ ein treuer Funktor ist.\hfill$\Box$}\\\par

Damit können wir nun das in der Einleitung genannte Theorem beweisen, wofür wir zunächst noch an \cite[Proposition 31]{DW1} erinnern, die im Beweis verwendet wird:\\\par

\begin{sat}[\protect{\cite[Proposition 31]{DW1}}]\label{31}
Es sei $\alpha\colon Y\to X$ eine Galoisüberlagerung von Varietäten über $\overline{\mathbb{Q}_{p}}$. Dann faktorisiert ein stetiger Funktor $W\colon\Pi_{1}(Y)\to\mathcal{C}$ in eine topologische Kategorie $\mathcal{C}$ genau dann in $W=V\circ\alpha_{*}$ für einen stetigen Funktor $V\colon\Pi_{1}(X)\to\mathcal{C}$, wenn $W\circ\sigma_{*}=W$ für jedes $\sigma\in\Gal(Y/X)$ gilt. Ist $\alpha$ nur endlich und étale und nicht notwendeigerweise galoissch, so bestimmt die Relation $W=V\circ\alpha_{*}$ bereits $V$ in eindeutiger Weise.\\
\end{sat}\par

\begin{thm}
Es sei $E$ ein $G$"=Torseur in $\mathcal{B}_{X_{\mathbb{C}_{p}}}^{ps}(G)$, also ein $G$"=Torseur auf $X_{\mathbb{C}_{p}}$ mit potentiell streng semistabiler Reduktion vom Grad Null.\par
Dann gibt es funktorielle Isomorphismen von "`Paralleltransport"' entlang étaler Wege der Fasern von $E_{\mathbb{C}_{p}}$ auf $X_{\mathbb{C}_{p}}$. Insbesondere gibt es für jeden Torseur $E\in\mathcal{B}^{s}_{X_{\mathbb{C}_{p}}, D}(G)$ einen eindeutig bestimmten stetigen Funktor $\rho_{E}$ von $\Pi_{1}(X)$ nach $\mathcal{P}(G(\mathbb{C}_{p}))$ mit $\rho_{E}(x)=E_{x}$ für alle $x\in X(\mathbb{C}_{p})$, so dass $\rho_{E}$ verträglich mit den oben definierten Funktoren $\rho_{\alpha^{*}E}$ von $\Pi_{1}(Y)$ nach $\mathcal{P}(G(\mathbb{C}_{p}))$ ist, falls $\alpha^{*}E\in\mathcal{B}_{Y_{\mathbb{C}_{p}}}^{s}$ für einen endlichen étalen Morphismus $\alpha\colon Y\to X$ von glatten und projektiven Kurven über $\overline{\mathbb{Q}_{p}}$ gilt.\par
Insgesamt erhält man so einen Funktor $\rho\colon\mathcal{B}_{X_{\mathbb{C}_{p}}}^{ps}(G)\to\Rep_{\Pi_{1}(X)}(G(\mathbb{C}_{p}))$, der sich funktoriell bezüglich Morphismen von glatten und projektiven Kurven über $\overline{\mathbb{Q}_{p}}$, Morphismen von zusammenhängenden reduktiven Gruppenschemata von endlicher Präsentation über $\mathfrak{o}$ und $\mathbb{Q}_{p}$"=Automorphismen von $\overline{\mathbb{Q}_{p}}$ verhält und mittels der in §2.4. dargestellten Funktoren verträglich mit den entsprechenden Funktoren im Vektorbündelfall ist.\par
Für jeden Punkt $x_{0}\in X(\mathbb{C}_{p})$ ist der "`Faserfunktor in $x_{0}$"' \[\mathcal{B}_{X_{\mathbb{C}_{p}}}^{ps}\to \mathcal{P}(G(\mathbb{C}_{p})), E\mapsto E_{x_{0}}, f\mapsto f_{x_{0}}\] ein treuer Funktor.
\end{thm}
\bewub{Es sei $E$ ein $G$"=Torseur auf $X_{\mathbb{C}_{p}}$ mit potentiell streng semistabiler Reduktion vom Grad Null. Dann existiert ein endlicher étaler Morphismus $\alpha\colon Y\to X$ von glatten, projektiven Kurven über $\overline{\mathbb{Q}_{p}}$, so dass $E\in\mathcal{B}_{Y_{\mathbb{C}_{p}}}^{s}$ ist. Ohne Beschränkung der Allgemeinheit darf man annehmen, dass $\alpha$ galoissch ist. Nach Theorem \ref{s} gilt somit \[\rho_{\alpha^{*}E}\circ\sigma_{*}=\rho_{\sigma^{*}(\alpha^{*}E)}=\rho_{\alpha^{*}E}\] für alle \mbox{$\sigma\in\Gal(Y/X)$}. Nach Satz \ref{31} existiert daher ein eindeutig bestimmter stetiger Funktor \[\rho(E)=\rho_{E}\colon\Pi_{1}(X)\longrightarrow\mathcal{P}(G(\mathbb{C}_{p})),\] so dass $\rho_{\alpha{*}E}=\rho_{E}\circ\alpha_{*}$ gilt. Insbesondere gilt $\rho_{E}(x)=E_{x}$ für alle $x\in X(\mathbb{C}_{p})$. Für einen étalen Weg $\gamma$ von $x_{1}$ nach $x_{2}$ in $X$ hat man ferner den Morphismus \[\rho_{E}(\gamma)=\rho_{\alpha^{*}E}(\gamma')\colon E_{x_{1}}=(\alpha^{*}E)_{y_{1}}\longrightarrow(\alpha^{*}E)_{y_{2}}=E_{x_{2}}.\] Dabei liegt \mbox{$y_{1}\in Y(\mathbb{C}_{p})$} über $x_{1}$ und $\gamma'$ ist der eindeutig bestimmmte Weg in $Y$ von $y_{1}$ zu dem Punkt $y_{2}$ über $x_{2}$, für den \mbox{$\alpha_{*}\gamma'=\gamma$} gilt. Für einen Morphismus $f\colon E\to E'$ von $G$"=Torseuren in $\mathcal{B}^{ps}_{X_{\mathbb{C}_{p}}}(G)$ ist der Morphismus \mbox{$\rho(f)=\rho_{f}\colon\rho_{E}\to\rho_{E'}$} durch die Familie der Abbildungen $f_{x}\colon E_{x}\to E'_{x}$ für alle $x\in X(\mathbb{C}_{p})$ definiert.\par
Wir behaupten nun, dass diese Konstruktion einen wohldefinierten Funktor \[\rho\colon\mathcal{B}^{ps}_{X_{\mathbb{C}_{p}}}(G)\longrightarrow\Rep_{\Pi_{1}(X)}(G(\mathbb{C}_{p}))\] definiert, der die vorher definierten Funktoren $\rho$ auf $\mathcal{B}_{X_{\mathbb{C}_{p}}}^{s}(G)$ ausdehnt. Dazu zeigen wir als erstes, dass die Definition von $\rho_{E}$ unabhängig von der Wahl von $\alpha$ ist.\par
Sind also $\alpha_{1}\colon Y_{1}\to X$ und $\alpha_{2}\colon Y_{2}\to X$ zwei endliche étale Galoisüberlagerungen von glatten und projektiven Kurven über $\overline{\mathbb{Q}_{p}}$, so existiert eine dritte $\alpha_{3}\colon Y_{3}\to X$, die $\alpha_{1}$ und $\alpha_{2}$ dominiert, d.\,h. es gibt Morphismen $\pi_{i}\colon Y_{3}\to Y_{i}$ für $i=1,2$, so dass jeweils $\alpha_{3}=\alpha_{i}\circ\pi_{i}$ gilt. Angenommen, es sei $\alpha_{i}^{*}E\in\mathcal{B}_{Y_{\mathbb{C}_{p}}}^{s}(G)$ für $i=1,2$. Nach der obigen Konstruktion gibt es dann für $i=1,2$ Funktoren $\rho_{i}\colon\Pi_{1}(X)\to\mathcal{P}(G(\mathbb{C}_{p}))$, so dass jeweils $\rho_{\alpha_{i}^{*}E}=\rho_{i}\circ\alpha_{i*}$ gilt. Es ist nun zu zeigen, dass $\rho_{1}=\rho_{2}$ gilt. Nach Theorem \ref{s} gilt \[\rho_{\alpha_{3}^{*}E}=\rho_{\pi_{i}^{*}(\alpha^{*}E)}=\rho_{\alpha_{i}^{*}E}\circ\pi_{i*}=\rho_{i}\circ\alpha_{i*}\circ\pi_{i*}=\rho_{i}\circ\alpha_{3*}\] für $i=1,2$. Die Eindeutigkeitsaussage in Satz \ref{31} impliziert dann $\rho_{1}=\rho_{2}$.\par
Als nächstes wird gezeigt, dass für jeden Morphismus $f\colon E\longrightarrow E'$ in $\mathcal{B}_{X_{\mathbb{C}_{p}}}^{ps}(G)$ die Familie von Abbildungen $f_{x}\colon E_{x}\longrightarrow E'_{x}$ einen Morphismus in $\Rep_{\Pi_{1}(X)}(G(\mathbb{C}_{p}))$ definiert. Ohne Beschränkung der Allgemeinheit sei angenommen, dass sowohl $\alpha^{*}E$ als auch $\alpha^{*}E'$ in $\mathcal{B}_{Y_{\mathbb{C}_{p}}}^{s}(G)$ liegen. Dann definiert $\rho_{\alpha^{*}f}$, d.\,h. die Familie der Abbildungen $(\alpha^{*}f)_{y}\colon(\alpha^{*}E)_{y}\longrightarrow(\alpha^{*}E')_{y}$, einen Morphismus in $\Rep_{\Pi_{1}(Y)}(G(\mathbb{C}_{p}))$. Nach Definition von $\rho_{E}(\gamma)$ für einen étalen Weg in $X$ kommutiert dann das Diagramm \begin{equation*}
        \xymatrix@=3em{%
         E_{x_{1}} \ar[r]^{f_{x_{1}}} \ar[d]^{\rho_{E}(\gamma)} & E'_{x_{1}}\ar[d]^{\rho_{E'}(\gamma)} \\ 
       E_{x_{2}} \ar[r]^{f_{x_{2}}} & E'_{x_{2}} }
           \end{equation*} für jeden étalen Weg $\gamma$ in $X$, so dass also die Familie $(f_{x})_{x\in X(\mathbb{C}_{p})}$ einen Morphismus in $\Rep_{\Pi_{1}(X)}(G(\mathbb{C}_{p}))$ definiert.\par
 Also ist $\rho$ wohldefiniert und nach Konstruktion ist klar, dass $\rho$ ein Funktor ist und den auf $\mathcal{B}_{X_{\mathbb{C}_{p}}}^{s}(G)$ gegebenen Funktor fortsetzt.\par
Für den Beweis der behaupteten Funktorialitäten betrachte man stellvertretend für den der übrigen den der Tatsache, dass $\rho$ verträglich mit Morphismen $X\to X'$ von glatten und projektiven Kurven über $\overline{\mathbb{Q}_{p}}$ ist:\par
Es sei also $f\colon X\to X'$ ein Morphismus glatter und projektiver Kurven über $\overline{\mathbb{Q}_{p}}$. Dann ist auf Ebene der Objekte zu zeigen, dass $\rho_{f^{*}E}=\rho_{E}\circ f_{*}$ für alle $E\in\mathcal{B}_{X'_{\mathbb{C}_{p}}}^{ps}(G)$ gilt.\par
Ist also $E\in\mathcal{B}_{X'_{\mathbb{C}_{p}}}^{ps}(G)$, so existiert per Definiton eine endliche étale Galoisüberlagerung $\alpha'\colon Y'\to X'$, so dass $\alpha'^{*}E\in\mathcal{B}_{Y'_{\mathbb{C}_{p}}}^{s}(G)$ gilt. Ist $Y$ die Normalisierung einer irreduziblen Komponente von $f^{-1}(Y')$, so kommutiert das Diagramm  \begin{equation*}
        \xymatrix@=3em{%
       Y \ar[r] \ar[d]^{\alpha} & f^{-1}(Y')=X\times_{X'}Y' \ar[r] & Y' \ar[d]^{\alpha'}\\ 
       X \ar[rr]^{f} & & X',}
           \end{equation*} wobei $\alpha\colon Y\to X$ der kanonische Morphismus ist, den man aus der Wahl von $Y$ erhält.\par
          Bezeichnet $g\colon Y\to Y'$ die obere horizontale Abbildung des Diagramms, so liegt dann $\alpha^{*}f^{*}E=g^{*}\alpha'^{*}E$ in $\mathcal{B}_{Y_{\mathbb{C}_{p}}}^{s}(G)$. Nach Theorem \ref{s} und der Definition von $\rho_{E}$ gilt ferner \[\rho_{\alpha^{*}f^{*}E}=\rho_{g^{*}\alpha'^{*}E}=\rho_{\alpha'^{*}E}\circ g_{*}=\rho_{E}\circ\alpha'_{*}\circ g_{*}=\rho_{E}\circ f_{*}\circ\alpha_{*}.\]
          Andererseits ist $\alpha\colon Y\to X$ eine endliche étale Galoisüberlagerung, so dass $f^{*}E$ in $\mathcal{B}_{X_{\mathbb{C}_{p}}}^{ps}(G)$ liegt. Also ist per Definition $\rho_{f^{*}E}$ der eindeutig bestimmte Funktor mit $\rho_{\alpha^{*}f^{*}E}=\rho_{f^{*}E}\circ\alpha_{*}$. Damit folgt $\rho_{E}\circ f_{*}=\rho_{f^{*}E}$. Es ist außerdem klar, dass das Diagramm \begin{equation*}
        \xymatrix@=3em{%
         \mathcal{B}_{X'_{\mathbb{C}_{p}}}^{ps}(G) \ar[rr]^{\rho} \ar[d]^{f^{*}} & & \Rep_{\Pi_{1}(X')}(G(\mathbb{C}_{p})) \ar[d]^{A(f)} \\ 
       \mathcal{B}_{X_{\mathbb{C}_{p}}}^{ps}(G) \ar[rr]^{\rho} &  & \Rep_{\Pi_{1}(X)}(G(\mathbb{C}_{p}))}
           \end{equation*} auf dem Niveau der Morphismen kommutiert.\par
Dass die Konstruktion des Funktors $\rho_{E}$ mit Morphismen $G\to G'$ zusammenhängender reduktiver Gruppenschemata von endlicher Präsentation über $\mathfrak{o}$ sowie mit Galoiskonjugation durch $\mathbb{Q}_{p}$"=Automorphismen von $\overline{\mathbb{Q}_{p}}$ verträglich ist sowie ähnlich wie im Fall der Kategorien $\mathcal{B}_{\mathfrak{X}_{\mathfrak{o}},D}(G)$ bzw. $\mathcal{B}_{X_{\mathbb{C}_{p}},D}(G)$ oder $\mathcal{B}_{X_{\mathbb{C}_{p}}}^{s}(G)$ mit den entsprechenden Funktoren für den Fall der Vektorbündel (siehe \cite{DW1}) kompatibel ist, folgt in analoger Weise ebenfalls aus den bekannten Funktorialitäts"= und Verträglichkeitsaussagen für $\rho_{\alpha^{*}E}$ und der Eindeutigkeitsaussage in Satz \ref{31}.\par
Dass der Faserfunktor \[\mathcal{B}_{X_{\mathbb{C}_{p}}}^{ps}(G)\to\mathcal{P}(G(\mathbb{C}_{p}))\] ein treuer Funktor ist, zeigt man genauso wie in Theorem \ref{s}. \hfill$\Box$}\\\par

\section{Offene Fragen}

\begin{itemize}
\item Wie groß ist die Kategorie $\mathcal{B}_{X_{\mathbb{C}_{p}}}^{ps}(G)$? Liegt insbesondere jeder semistabile $G$"=Torseur vom Grad Null auf $X_{\mathbb{C}_{p}}$ in dieser Kategorie? Die analoge Frage im Fall der Vektorbündel, d.\,h. ob jedes semistabile Vektorbündel vom Grad Null auf $X_{\mathbb{C}_{p}}$ in $\mathcal{B}^{ps}_{X_{\mathbb{C}_{p}}}$ liegt, ist bislang nur für Kurven $X$ vom Geschlecht $g=0$ oder $g=1$ positiv beantwortet worden, für Kurven $X$ von Geschlecht $g\geq 2$ aber auch dort weiterhin unbeantwortet.
\item Was ist das essentielle Bild des Funktors $\rho\colon\mathcal{B}_{X_{\mathbb{C}_{p}}}^{ps}(G)\to\Rep_{\Pi_{1}(X)}(G(\mathbb{C}_{p}))$? Auch dies ist ebenfalls eine offene Frage im Fall der Vektorbündel.
\item Ist der Faserfunktor in einem fest gewählten Punkt $x_{0}\in X(\mathbb{C}_{p})$, \[\mathcal{B}_{X_{\mathbb{C}_{p}}}^{ps}(G)\to\mathcal{P}(G(\mathbb{C}_{p})),\] nicht nur treu, sondern sogar auch volltreu?
\end{itemize} \newpage\thispagestyle{empty}

\markboth{LITERATUR}{}\thispagestyle{plain}
\addcontentsline{toc}{chapter}{\numberline{}Literatur}

\end{document}